\documentclass{article}
\usepackage{amssymb}
\usepackage{amsmath}
\usepackage{color}

\def\thebibliography#1{\section*{References}\list
  {[\arabic{enumi}]}{\settowidth\labelwidth{[#1]}\leftmargin\labelwidth
    \advance\leftmargin\labelsep
    \usecounter{enumi}}
    \def\newblock{\hskip .11em plus .33em minus -.07em}
    \sloppy
    \sfcode`\.=1000\relax}

\newcommand{\refbook}[3]{{\sc #1}{\em\ #2}{\ #3}}
\newcommand{\refer}[5]{{\sc #1}{\ #2}{\em\ #3}{\bf\ #4}{\ #5}}

\newtheorem{lem}{Lemma}[section]
\newtheorem{cor}[lem]{Corollary}
\newtheorem{teo}[lem]{Theorem}
\newtheorem{os}[lem]{Remark}
\newtheorem{defi}[lem]{Definition}
\newtheorem{prop}[lem]{Proposition}
\newtheorem{esem}[lem]{Example}

\newcommand\sign{\mathop{\rm sign}}
\newcommand{\qed}{\thinspace\null\nobreak\hfill\hbox{\vbox{\kern-.2pt\hrule
 height.2pt depth.2pt\kern-.2pt\kern-.2pt \hbox to2.5mm{\kern-.2pt\vrule
 width.4pt \kern-.2pt\raise2.5mm\vbox to.2pt{}\lower0pt\vtop
 to.2pt{}\hfil\kern-.2pt \vrule
 width.4pt \kern-.2pt}\kern-.2pt\kern-.2pt\hrule height.2pt depth.2pt
 \kern-.2pt}}\par\medbreak}

\newcommand{\R}{\mathbb{R}}

\newcommand{\C}{\mathbb{C}}

\newcommand{\N}{\mathbb{N}}

\newcommand{\Rp}{\textrm{\emph{Re}\,}}

\newcommand{\eps}{\varepsilon}

\newcommand{\Om}{\Omega}

\newcommand{\ov}{\overline}

\date{}
\begin{document}

\title{Rellich inequalities in   bounded domains}
\author{{G. Metafune \thanks{Dipartimento di Matematica ``Ennio De Giorgi'', Universit\`a del Salento, C.P.193, 73100, Lecce, Italy.
email:  giorgio.metafune@unisalento.it}  \qquad L. Negro  \thanks{Dipartimento di Matematica ``Ennio De Giorgi'', Universit\`a del Salento, C.P.193, 73100, Lecce, Italy.
email:  luigi.negro@unisalento.it}  \qquad M. Sobajima \thanks{Department of Mathematics, Tokyo University of Science, Japan. email: msobajima1984@gmail.com}\qquad C. Spina} \thanks{Dipartimento di Matematica ``Ennio De Giorgi'', Universit\`a del Salento, C.P.193, 73100, Lecce, Italy.
email:  chiara.spina@unisalento.it}   }

\maketitle
\begin{abstract}
We find necessary and sufficient conditions for the validity of weighted Rellich 
inequalities in $L^p$, $1\le p \le \infty$, for functions in bounded domains vanishing at the boundary. General operators like $L=\Delta+c\frac{x}{|x|^2}\cdot\nabla-\frac{b}{|x|^2}$ are considered. Critical cases and remainder terms are also investigated.

\par

\end{abstract}

\tableofcontents

\section{Introduction}
In this paper we consider the operator 

\begin{equation} \label{L} Lu=\Delta u+c\frac{x}{|x|^2}\cdot\nabla u-\frac{b}{|x|^2} u,\quad  c,\ b\in\R
\end{equation} 
acting in the space $L^p(\Omega)$, for $1\le p \le \infty$, endowed with Dirichlet boundary conditions and we determine all $\alpha's$ (depending on $N,p,c,b$) for which  the following  weighted Rellich 
inequalities hold 
\begin{equation} \label{Intr 1}
\||x|^\alpha Lu\|_p \geq C\||x|^{\alpha-2} u\|_p.
\end{equation}
Note that, when $c=0$, $L$ becomes a Schr\"odinger operator with inverse square potential.
When best constants can be computed, we prove that they are not attained by adding remainder terms. Finally,  when Rellich inequalities above fail, we prove  modified inequalities which include logarithmic terms. \\
The first  results in this direction have been obtained for the Laplacian in unweighted $L^p$-spaces and when $\Omega=\R^N$. In 1956, Rellich proved the inequalities 
$$\left(\frac{N(N-4)}{4}\right)^2\int_{\R^N}|x|^{-4}|u|^2\, dx\leq \int_{\R^N}|\Delta u|^2\, dx$$
for $N\not =2$ and for every $u\in C_c^\infty (\R^N\setminus\{0\})$, see \cite{rellichF}.
These inequalities have been then extended to $L^p$-norms:  in 1996, Okazawa proved in \cite{okazawa} the validity of
$$\left(\frac{N}{p}-2\right)^p\left(\frac{N}{p'}\right)^p\int_{\R^N}|x|^{-2p}|u|^p\, dx\leq \int_{\R^N}|\Delta u|^p\, dx$$ for $1<p<\frac{N}{2}$,  showing also the optimality of the constants. 

Weighted Rellich inequalities have also been studied in \cite{davi-hinz} and later by 
Mitidieri  who proved for $N\geq 3$ and for $2-\frac{N}{p}<\alpha<\frac{N}{p'}$ 
\begin{equation} \label{WRI}
C^p(N,p,\alpha)\int_{\R^N}|x|^{(\alpha-2)p}|u|^p\, dx\leq \int_{\R^N}|x|^{\alpha p}|\Delta u|^p\, dx
\end{equation}
with the optimal constants $C^p(N,p,\alpha)=\left(\frac{N}{p}-2+\alpha\right)^p\left(\frac{N}{p'}-\alpha\right)^p$,
see  \cite[Theorem 3.1]{mitidieri}.\\
In the   recent paper  \cite{caldiroli},  Caldiroli and Musina 
improved  weighted Rellich inequalities for $p=2$ by giving necessary and sufficient conditions on $\alpha$ for the validity of (\ref{WRI}) 
and finding also the optimal constants $C^2(N,2,\alpha)$. 
In particular they  proved that  (\ref{WRI}) is verified for $p=2$ 
if and only if $\alpha\neq N/2+n$, $\alpha \neq -N/2+2-n$ for every $n\in\N_0$.

In \cite{rellich} 
the results in \cite{caldiroli} are extended to $1 \le p \le \infty$, computing also best constants in some cases. It is shown that (\ref{WRI}) holds if and only if $\alpha \neq N/p'+n$, $\alpha \neq -N/p+2-n$ for every $n \in \N_0$.
 Moreover, Rellich inequalities are employed to find necessary and sufficient conditions for the validity of weighted Calder\'on-Zygmund estimates when $1<p<\infty$.
These methods can be  applied 
to general operators as in (\ref{L}),  thus providing
a complete solution  to problem  (\ref{Intr 1}) with $\Omega=\R^N$. 

\smallskip

Let us now consider bounded open sets $\Omega$ containing the origin and spaces of functions vanishing at the boundary.  In contrast with Hardy inequality, where many results in bounded domains improving those in the whole space are known, Rellich inequalities do not seem to have been studied intensively. We quote however  \cite{musina} for $L=\Delta$, where the author discovers a range of parameters $\alpha$ where Rellich inequalities hold in the whole space but not in a bounded $\Omega$, due to the boundary conditions. 

In this paper we find all parameters $\alpha$ for which (\ref{Intr 1}) hold for a general $L$ as in (\ref{L}), assuming that $\Omega$ has a smooth boundary and the condition $D:=b+(N-2+c)^2/4\geq 0$ on the coefficients of $L$, which guarantees the solvability of related elliptic problems. When $\Omega$ is a ball, however, this restriction on $D$ is not necessary.  
Our method is based on the spectral analysis of the auxiliary operator $A=|x|^2 \Delta +cx\cdot \nabla$, as explained in Section 2.
In particular, we show that, setting $\lambda_n=n(N-2+n)$, (\ref{L}) holds if and only if
\begin{align} \label{eq1}
\nonumber \alpha&<N\Bigl (\frac12-\frac1p\Bigr)+1+\frac{c}{2}+ \sqrt{ D} \quad\text{and}\;\\[1ex]
\alpha&\neq N\Bigl (\frac12-\frac1p\Bigr)+1+\frac{c}{2}- \sqrt {D+\lambda_n}, \quad\forall\, n\in \N_{0}.
\end{align} 
When $\Omega$ is a ball centered at the origin, the above characterization holds also when $D<0$ (changing the square roots with their  real parts) and in the extreme cases $p=1, \infty$. However, when $\Omega=\R^N$ the results in \cite{rellich} say that Rellich inequalities hold if and only if
\begin{equation} \label{eq3}
\alpha\neq N\Bigl (\frac12-\frac1p\Bigr)+1+\frac{c}{2}\pm \Rp\sqrt {D+\lambda_n}, \quad\forall\, n\in \N_{0}.
\end{equation}
The reason for the difference between \eqref{eq2} and \eqref{eq3} is explained in Section 2 in an elementary way in the case of the ball, by showing explicit counterexamples due to the boundary. 

Rellich inequalities  can be proved by using integration by parts and applying Hardy-type inequalities only when
\begin{equation} \label{eq2}
N\Bigl (\frac12-\frac1p\Bigr)+1+\frac{c}{2}-\sqrt {D}<\alpha < N\Bigl (\frac12-\frac1p\Bigr)+1+\frac{c}{2}+\sqrt {D}.
\end{equation}
This proof allows also to compute the best constant $C:=b+\Bigl(\frac{N}{p}-2+\alpha\Bigr)\Bigl(\frac{N}{p'}-\alpha+c\Bigr)$. For the other values of $\alpha$ appearing in \eqref{eq1}, the best constant is unknown unless $p=2$, see \cite{caldiroli}, \cite{rellich}, or when $p$ is generic but special subspaces of $L^p$ are considered, see \cite{rellich}. 

In the range \eqref{eq2}, Rellich inequalities have essentially a  one dimensional structure, since the (approximate) extremants are radial functions and best constants can therefore be computed. Outside of this range, however, the problem loses its rotational symmetry and the extremants, in special subspaces,  involve spherical harmonics, see \cite{rellich}, again. This explains also why symmetrization arguments based on spherical rearrangements do not work and a spectral analysis appears. Similarly, best constants can be computed on subspaces of $L^p$ which allow a one-dimensional reduction and then on the whole $L^2$, by orthogonal expansions.

Remainder terms are known for the Laplacian in the unweighted case. We quote  \cite{tertikas-zogra} where the authors obtained in particular
\begin{align*}
\int_{\Omega}
|\Delta u|^2
\,dx &
\geq \left(\frac{N(N-4)}{4}\right)^2\int_{\Omega}
\frac{|u|^{2}}{|x|^{4}}
\,dx\\[1ex]&+\left(1+\frac{N(N-4)}{8}\right)\sum_{i=1}^\infty\int_{\Omega}
\frac{|u|^{2}}{|x|^{4}}
X_1^2X_2^2\cdots X_i^2\,dx,
\end{align*}  
for bounded domains $\Omega$ in $\R^N$, $N\geq 5$,  $u\in C_c^\infty(\Omega\setminus\{0\})$, where $X_k=X_k\left(\frac{|x|}{R(\Omega)}\right)$, $R(\Omega)=\sup_{x\in \Omega}|x|$, are iterated radial logarithmic functions.
The result has been extended to $L^p$ norms in  \cite{barbatis-tertikas} under the restriction  $p<\frac{N}{2}$, according to \eqref{eq2} when $\alpha=b=c=0$. A different proof which uses symmetrization and covers also the case $p=\frac{N}{2}$ is given in  \cite{ando}. Rellich inequalites with remainder terms in the whole space have been investigated  in \cite{sano}, where the remainder is given in terms of weighted $L^q$ norms of the Schwartz symmetrization of the functions. 

We prove a similar result for our operator $L$ in weighted $L^p$ norms, considering only one remainder term. When $\alpha$ satisfies \eqref{eq2} we obtain with $C$ above  
\begin{equation*} 
\Big\||x|^\alpha Lu\Big\|_p^p -C^p \Big\||x|^{\alpha-2} u\Big\|_p^p \ge c \Big\||x|^{\alpha-2}\left|\log |R^{-1}x|\right|^{-\frac{2}{p}} u\Big\|_p^p
\end{equation*}
for $u\in C^2_c(B_{R/2}\setminus\{0\})$.
\smallskip
Some explanation on the class of functions here considered is necessary. Since \eqref{eq2} is satisfied, Rellich inequalities hold for both $\Omega$ bounded or $\Omega=\R^N$ but we choose to formulate the above result with reference to the whole space, that is for functions having compact support. A similar formulation for functions only vanishing at $\partial \Omega$, when $\Omega$ is a ball, is also possible but we prefer to point out only the role of the singularity at $0$, since the weight $|x|^\alpha$ has no effect on the boundary.


In the critical cases, when Rellich inequalities  do not hold, we prove that modified inequalities with  logarithmic correction terms are still valid. Again we focus on the singularity at $0$ and consider functions with compact support in $\R^N$. If $$\alpha= N\Bigl (\frac12-\frac1p\Bigr)+1+\frac{c}{2}\pm \Rp \sqrt {D+\lambda_n}$$ for some $n \in \N_0$, $1<p\le \infty$, then 
\begin{equation*}  
\||x|^\alpha Lu\|_p \geq C\Big\||x|^{\alpha-2}\left|\log |R^{-1}x|\right|^{-2} |u|\Big\|_p \quad {\rm when }\ D+\lambda_n \le 0
\end{equation*}
\begin{equation*}  
\||x|^\alpha Lu\|_p \geq C\Big\||x|^{\alpha-2}\left|\log |R^{-1}x|\right|^{-1} |u|\Big\|_p \quad {\rm when }\ D+\lambda_n >0
\end{equation*}
for $u\in C^2_c(B_{R/2}\setminus\{0\})$. When $p=1$, the previous inequalities  hold with $|\log |R^{-1}x||^{-2}$ and  $|\log |R^{-1}x||^{-1}$ replaced by $|\log |R^{-1}x||^{-2-\eps}$ and  $|\log |R^{-1}x||^{-1-\eps}$, respectively.\\
In this way we extend the results already proved in \cite{adimurthi} for the Laplace operator  under the more restrictive conditions $\alpha=0$, $p=\frac{N}{2}$, $N\geq 3$. We also refer to \cite{gazzola} where Rellich inequalities for the Laplacian have been proved with different remainder terms for $\alpha=0$, $p\leq \frac{N}{2}$.  

\smallskip

The treatment of the critical case does not rely on rearrangements, as already explained, but a reduction to the one-dimensional case is still possible via a spectral analysis. In fact we show that Rellich inequalities are true, even in the critical cases, if we consider subspaces of $L^p(\R^N)$ spanned by functions like $f(r)P(\omega)$, where $P$ is a spherical harmonic of degree different from $n$ and the problem is then reduced to find the right inequalities for (linear combinations of) functions $g(r)Q(\omega)$ where $Q$ is a spherical harmonic of degree $n$, hence to a finite number of one-dimensional problems.

\smallskip

Let us explain why semigroups of linear operators appear often in the paper. When $p=2$, Rellich inequalities can be reduced to  a countable set of one-dimensional inequalities, by an orthogonal expansion in spherical harmonics, see for example \cite{rellich}. Moreover, it turns out that is more convenient to work with the operator $A=|x|^2L$ instead of $L$, so that  the radial and the angular parts decouple. When $p \neq 2$ the one-dimensional analysis can be still performed but one needs a substitute for orthogonal expansions. This role is played by the semigroup $e^{tA}$ which allows to compute the spectrum of $A$, by tensor product arguments, since the radial and the angular parts commute. Rellich inequalities are equivalent to spectral inequalities for $A$ and, moreover, the description of the domain of $A$ allows us to identify precise classes where Rellich inequalities hold.

\smallskip
Let us briefly describe the content of the sections. 
In Section 2 we present the basic ideas and some explicit counterexamples which serve as a guide for the rest of the paper. We reduce Rellich inequalities to a spectral problem for an operator with singular coefficients $A=|x|^2 \Delta +cx\cdot \nabla$ which is therefore analysed in detail in Section 3, which is the core of the paper. Rellich inequalities for the ball and for the whole space are easily deduced in Section 4 from the analysis of Section 3. The case of general domains, without any rotational symmetry, is studied in Section 5: here we need $1<p<\infty$ and $D \ge 0$, a condition  which is known to be equivalent to the existence of positive solutions for elliptic and parabolic problems related to $L$. When $L=\Delta-b|x|^{-2}$, this condition reduces to the classical one $b+(N-2)^2/4 \ge 0$. The main tool to pass from the ball to a general $\Omega$ is a pointwise estimate of the Green function of $-L$ which follows from precise bounds of the  heat kernel. Rellich inequalities in exterior domains not containing the origin are easily treated via the Kelvin transform. 
In Section 6 we show that, when Rellich inequalities fail,  modified inequalities which include logarithmic terms are still valid. The situation is similar to Hardy inequality, when the classical one fails.
In Section 7, we analyse the remainder term in Rellich inequalities when \eqref{eq2} is satisfied.

\bigskip
\noindent\textbf{Notation.}
We denote by $\N_0=\N\cup\{0\}$ the natural numbers including 0.  If $\Omega$ is  an open  subset of $\R^N$,  $C_b(\Omega)$ is the Banach space of all continuous and bounded functions in $\Omega$, endowed with the sup-norm, $C_0(\overline{\Omega})$ its subspace consisting of functions vanishing at the boundary and  $C_0^0(\overline{\Omega})$ its subspace consisting of functions vanishing at the origin and at the boundary, when $0 \in \Omega$.  $C_c^\infty(\Omega)$ denotes the space of infinitely continuously differentiable functions with compact support in $\Omega$. 
The unit sphere $\{\|x\|=1\}$ in $\R^N$ is denoted by $S^{N-1}$; $\Delta_0$ is its Laplace-Beltrami operator. We adopt standard notation for $L^p$ and Sobolev spaces when $1 \le p<\infty$ but we use $L^\infty (\Omega)$ for $C_b(\Omega)$ to unify the notation.
$B_r$ is the ball of center $0$ and radius $r$, $B_r^c=\R^N\setminus B$. We write $B$ for $B_1$. For $V\subseteq\R^N$, we denote by  $\overset{\mathrm{o}}{V}$ the interior part of $V$. When $L$ is a closed operator $\sigma (L)$, $P\sigma (L)$, $A\sigma (L)$, $R\sigma (L)$, denote the spectrum, the point-spectrum, the approximate point spectrum and the residual spectrum, respectively. Definitions and the relevant properties are listed in the Appendix.

\section{Basic results and methods}\label{Preliminaries}
Let $L$ be as in (\ref{L}) and let $\Omega$ be an open, bounded,  connected subset of $\R^N$ containing the origin and with a smooth boundary, or $\Omega=\R^N$.
 For $1\leq p\leq\infty$,  $\alpha\in\R$ we define
\begin{align*}
D_{p,\alpha}(\Omega):&=\left\{u:\ |x|^{\alpha-2}u,\ |x|^{\alpha}Lu\in L^p\left(\Omega\right),\ u=0 \text{ on } \partial \Omega\right\}
\end{align*}
 $Lu$ is understood as a distribution in $\Omega \setminus \{0\}$. Since the coefficients of $L$ are $C^\infty$ away from  the origin,  by local elliptic regularity it follows that, if $u \in D_{p, \alpha}(\Omega)$, then   $u \in W^{2,p}_{loc}(\R^N \setminus \{0\})$ when $\Omega=\R^N$ and  $u \in W^{2,p}(\Omega \setminus B_\eps)$ for every $\eps >0$, when $\Omega$ is bounded. This clearly holds for $1<p<\infty$; when $p=\infty$, the same is true for any $q <\infty$.

Note that, 
when $\Omega$ is bounded, also the class
\begin{align*}
 D_{p,\alpha,0}(\Omega):=\{u \in  D_{p,\alpha}(\Omega) , u=0\  {\rm in\ a\ neighborhood\ of }\  \partial \Omega \}
 \end{align*}
could be considered.  However, since every function $u\in D_{p,\alpha,0}(\Omega)$, extended by $0$ to $\R^N$, belongs to $D_{p,\alpha}(\R^N)$, the problem is then reduced to the case of the whole space. A scaling argument, moreover,  shows that  Rellich inequalities  (\ref{Intr 1}) hold in $D_{p,\alpha,0}(\Omega)$ if and only if they hold in $D_{p,\alpha}(\R^N)$.

Defining
\begin{align*}
v(x)=|x|^{\alpha-2}u(x),
 \end{align*}
  it is straightforward to compute that  $|x|^\alpha L u=Av-\mu v$,  where
\begin{align} \label{A}
  A=|x|^2\Delta +(c+4-2\alpha)x\cdot\nabla  \  \  {\rm and } \ \  \mu=b-(2-\alpha)(N-\alpha+c).
\end{align}
Then 
 Rellich inequalities (\ref{Intr 1})
are equivalent to the spectral estimates
\begin{equation} \label{disv}
\|\mu v-Av\|_p \geq C\|v\|_p, \quad v\in D_{p,max}(\Omega) 
\end{equation} 
where
\begin{align*}
D_{p,max}(\Omega):&=\left\{u\in L^p(\Omega): Au\in L^p\left(\Omega\right),\ u=0 \text{ on } \partial \Omega\right\}
\end{align*}
and $Au$ is understood as a distribution as above. Moreover, the constants $C$ in (\ref{Intr 1}) and (\ref{disv}) are the same.

Inequalities (\ref{disv})
hold precisely when $\mu$ does not belong to the approximate point spectrum of $A$. This explains why  a large part of this paper is devoted to the  study of the operator $A$ and of the  fine structure of its spectrum.

In the next proposition we state the above reduction, for further reference, and prove a density result using the same method. We refer to Section \ref{spectral} for basic definitions and results from spectral theory.

\begin{prop} \label{reddensity}
Let $L$ be as in (\ref{L}) and let $\Omega$ be an open, bounded,  connected subset of $\R^N$ containing the origin and with a $C^{2, \beta}$ boundary, or $\Omega=\R^N$. Then 
\begin{itemize}
\item[(i)] Rellich inequalities  (\ref{Intr 1}) hold if and only if $\mu=b-(2-\alpha)(N-\alpha+c) $ does not belong to the approximate point spectrum of $(A,D_{p,max}(\Omega))$.
\item[(ii)] Rellich inequalities  (\ref{Intr 1})
hold for functions in $D_{p,\alpha}(\Omega)$ if and only if they hold for $C^2$ - functions vanishing in a neighbourhood of the origin and on $\partial \Omega$, when $ \Omega$ is bounded, or also in a neighbourhood of infinity, when $\Omega=\R^N$.
\end{itemize}
\end{prop}
{ Proof. } The discussion above shows that Rellich inequalities hold if and only the spectral inequalities (\ref{disv}) are valid in $D_{p,max}(\Omega)$, hence when $\mu$ does not belong to the approximate point spectrum of $A$, by Proposition \ref{Rellich-spectrum}. This proves (i). To prove (ii) it is sufficient to note that the transformation $v(x)=|x|^{\alpha-2}u(x)$ preserves the class of functions defined in (ii) and that, by Lemma \ref{coreOmega} and Proposition \ref{L1}, these functions constitute a core of $(A, D_{p,max}(\Omega))$.
\qed

The interplay between the operators $A$ and $L$ allows to give simple proofs of   Rellich inequalities in special cases where best constants can be computed.

\begin{prop}  \label{easy}
Let $\Omega$ be an open, bounded,  connected subset of $\R^N$  with a $C^{1}$ boundary, or $\Omega=\R^N$. Assume that $1 \le p \le \infty$, that  
$D:=b+\left(\frac{N-2+c}{2}\right)^2>0$ and that 
\begin{equation} \label{range}
 N\Bigl (\frac12-\frac1p\Bigr)+1+\frac{c}{2}-\sqrt {D}<\alpha < N\Bigl (\frac12-\frac1p\Bigr)+1+\frac{c}{2}+\sqrt {D}. 
\end{equation}
Then Rellich inequalities (\ref{Intr 1}) hold in $D_{p, \alpha}(\Omega)$ with $C:=b+\Bigl(\frac{N}{p}-2+\alpha\Bigr)\Bigl(\frac{N}{p'}-\alpha+c\Bigr)$. The constant $C$ is optimal when $\Omega$ contains the origin.
\end{prop}
{ Proof. } We have to show that (\ref{disv}) holds, with the constant $C$ above, for $A$ and $\mu$ defined in (\ref{A}). This is proved in Theorem  \ref{dissipativity}, using only integration by parts and Hardy inequality (change $c$ with $c+4-2\alpha$ and $\lambda-\omega_p$ with $\mu$, therein). We note that $C>0$ is equivalent to  (\ref{range}).

To prove the optimality of $C$, when $ 0 \in \Omega$, we observe that Rellich inequalities are invariant under dilations. If $C_\Omega$ is the best constant in $\Omega$, then $C_{r\Omega}=C_\Omega$ for any $r>0$. Letting $r \to \infty$ we see that $C_{\R^N} \le C_\Omega $. However, $C_{\R^N}=b+\Bigl(\frac{N}{p}-2+\alpha\Bigr)\Bigl(\frac{N}{p'}-\alpha+c\Bigr)$, by \cite[Theorem 3.1]{rellich}.
\qed

Note that when $L=\Delta$, then $D=(N-2)^2/4$ and (\ref{range}) reduces to $2-N/p <\alpha <N/p'$ and $C=\Bigl(\frac{N}{p}-2+\alpha\Bigr)\Bigl(\frac{N}{p'}-\alpha\Bigr)$. If $\Omega$ does not contain the origin the constant $C$ above is not optimal, in general, see again \cite[Section 6]{rellich} for the case of the half space.

Next, we show explicit counterexamples to Rellich inequalities already appeared in \cite{musina} when $L=\Delta$. We distinguish between free counterexamples depending on the singularity at zero, which appear in any  set $\Omega$ containing the origin and counterexamples where the boundary $\partial \Omega$ is involved, appearing only when $\Omega$ is bounded in addition to the preceding ones. We confine here only to the case of the unit ball $B$; the general case will be treated in Section \ref{Rellich Bounded domain}.

We employ spherical coordinates on $\R^N\setminus\{0\}$ and write $x=r\omega$, where $r:=|x|$, $\omega:=x/|x|\in {\mathbb{S}}^{N-1}$. Then
\begin{align*}
L=D_{rr}+\frac{N-1+c}{r} D_r-\frac{b-\Delta_0}{r^2},
\end{align*}
 where $D_{rr}$, $D_r$ denote radial derivatives and $\Delta_0$ is  the Laplace-Beltrami operator on the unit sphere $\mathbb{S}^{N-1}$. Let $P$ be a spherical harmonics of order $n\in N_0$, with $\Delta_0 P=-\lambda_n P$, $\lambda_n=n(N+n-2)$. If $u(r\omega)=v(r)P(\omega)$ then
\begin{align*}
Lu=\left[v_{rr}+\frac{N-1+c}{r} v_r-\frac{b+\lambda_n}{r^2}v\right]P.
\end{align*}
The equation $Lu=0$ has  solutions $|x|^{-s_1^n}P$, $|x|^{-s_2^n}P$ where the function  $r^{-s_1^n}$, $r^{-s_2^n}$  solve 
$$v_{rr}+\frac{N-1+c}{r} v_r-\frac{b+\lambda_n}{r^2}v=0.$$  $s_1^n,s_2^n$ are the roots of the indicial equation $f(s)=-s^2+(N-2+c)s+b+\lambda_n=0$ given by

\begin{equation} \label{defs}
s_1^n:=\frac{N-2+c}{2}-\sqrt{D+\lambda_n},
\quad
s_2^n:=\frac{N-2+c}{2}+\sqrt{D+\lambda_n}
\end{equation}
where
\bigskip
\begin{equation} \label{defD}
D:=
b+\left(\frac{N-2+c}{2}\right)^2.
\end{equation}

The following Examples shows that, due to the singularity of $L$ at $0$,   Rellich inequalities  always fail when $\alpha$ equals one of the values
\begin{align*} 
\alpha_n^{\pm}:= N\Bigl (\frac12-\frac1p\Bigr)+1+\frac{c}{2}\pm\Rp \sqrt {D+\lambda_n}, \quad \,n\in\N_0,
\end{align*}
\begin{esem} \label{esem1}
Let $1\leq p\leq\infty$ and let $\Omega\subseteq \R^N$ be an open subset of $\R^N$ such that $0\in\Omega$. If $\alpha=\alpha_n^{\pm}$, then  Rellich inequalities \eqref{Intr 1} do not hold in $D_{p,\alpha}(\Omega)$.
\end{esem}
{\sc{Proof.}} Suppose, for example, that $\alpha=\alpha_n^-$. Let $s_1^n$ be  defined in \eqref{defs} and $\gamma=-\Rp s_1^n$. We fix $R>0$ such that $B_R\subseteq\Omega$ and  take $P$ a spherical harmonics of order $n$. The function 
\begin{align*}
u(r\omega):=r^\gamma P(\omega), \quad x=r\omega\in B_R
\end{align*}
satisfies $Lu=0$ but  $|x|^{\alpha_n^--2} u\notin L^p\left(B_r\right)$ 
  since
\begin{align}\label{Counterexample 0-1}
\alpha_n^--2+\gamma=-\frac{N}{p},\quad 1\leq p\leq\infty.
\end{align}
Let $\varphi\in C^\infty(\R)$ such that $\mbox{supp}\,\varphi\subseteq [\frac 1 4,\frac 1 2]$ and $\varphi_\epsilon(r):=\varphi(r^\epsilon)$. By construction $u_\epsilon:=u\varphi_\epsilon$ has support in $ [\left(\frac 1 4\right)^{\frac 1 \epsilon},\left(\frac 1 2\right)^{\frac 1 \epsilon}]$, lies in $D_{p,\alpha}\left(\Omega\right)$ and satisfies
\begin{align*}
Lu_\epsilon(r\omega)=P(\omega)\left[r^\gamma\varphi_\epsilon''+(2\gamma+N-1+c)r^{\gamma-1}\varphi_\epsilon'\right].
\end{align*}
If $1\leq p<\infty$ and  $\bar r>0$ such that $\mbox{supp}\,\varphi_\epsilon\subseteq B_{\bar r}$ we get
\begin{align*}
\int_{\Omega}|x|^{(\alpha_n^--2)p} |u_\epsilon|^p\,dx&=\int_{B_{\bar r}}|x|^{(\alpha_n^--2+\gamma)p} |P(\omega)|^p|\varphi_\epsilon|^p\,dx=C\int_0^{\bar r}\frac{|\varphi(r^\epsilon)|^p}{r}\,dr=\frac{C}{\epsilon}\int_{\frac 1 4}^{\frac 1 2}\frac{|\varphi(s)|^p}{s}\,ds,
\end{align*}
where  $C=\int_{\mathbb{S}^{N-1}} |P(\omega)|^p\,d\omega$. On the other hand
\begin{align*}
\int_{\Omega}|x|^{\alpha_n^-p} |Lu_\epsilon|^p\,dx=C\,\epsilon^{p-1}\int_{\frac 1 4}^{\frac 1 2}s^{p-1}\left|\epsilon s\varphi''(s)+(2\gamma+N-2+c+\epsilon)\varphi'(s)\right|^p\,ds.
\end{align*}
It follows, from the previous equalities, that
\begin{align*}
\frac{\int_{\Omega}|x|^{\alpha_n^-p} |Lu_\epsilon|^p\,dx}{\int_{\Omega}|x|^{(\alpha_n^--2)p} |u_\epsilon|^p\,dx}=\epsilon^p\,\frac{\int_{\frac 1 4}^{\frac 1 2}s^{p-1}\left|\epsilon s\varphi''(s)+(2\gamma+N-2+c+\epsilon)\varphi'(s)\right|^p\,ds}{\int_{\frac 1 4}^{\frac 1 2}\frac{|\varphi(s)|^p}{s}\,ds}
\end{align*}
which tends to $0$ as $\epsilon \to 0$, hence  Rellich inequalities do not hold in $D_{p,\alpha}(\Omega)$ for $1\leq p<\infty$.\\
If  $p=\infty$, then $\alpha_n^--2+\gamma=0$ and an analogous computation   yields
\begin{align*}
|x|^{\alpha_n^--2} u_\epsilon(x)&=P(\omega)\varphi(r^\epsilon),\\[1ex]
|x|^{\alpha_n^-p} Lu_\epsilon(x)
&=P(\omega)\left[r^{2\epsilon}\epsilon^2\varphi''(r^\epsilon)+\epsilon(2\gamma+N-2+c+\epsilon)r^\epsilon\varphi'(r^\epsilon)\right]
.
\end{align*}
This implies
\begin{align*}
\frac{\||x|^{\alpha_n^-}Lu_\epsilon\|_\infty}{\||x|^{\alpha_n^--2} u_\epsilon\|_\infty}=\frac{\epsilon\sup_{s\in[\frac 1 4,\frac 1 2]}\left|\epsilon s^2\varphi''(s)+(2\gamma+N-2+c+\epsilon)s\varphi'(s)\right|}{\sup_{s\in[\frac 1 4,\frac 1 2]}|\varphi(s)|}
\end{align*}
which tends to $0$ as $\epsilon \to 0$.
The proof for $\alpha=\alpha_n^+$ is similar,  choosing  $\gamma=-\Rp s_2^n$.
\\\qed

Next we consider the case where $\Omega=B$ and show that, due to the Dirichlet boundary condition at $\partial B$,  new counterexamples appear, in addition to the previous ones. 
The same result is proved in Section \ref{Rellich Bounded domain} for general bounded domains. 
\begin{prop}\label{Counterexample in B}
 If $ 1 \le p \le \infty$ and $\alpha> N\Bigl (\frac12-\frac1p\Bigr)+1+\frac{c}{2}+ \Rp\sqrt{D}$, then the  Rellich inequalities \eqref{Intr 1} cannot hold in $D_{p,\alpha}(B)$.
\end{prop}
{\sc{Proof.}} Let $\alpha> N\Bigl (\frac12-\frac1p\Bigr)+1+\frac{c}{2}+ \Rp\sqrt{D}$ and let $s_{1,2}$ be  defined in \eqref{defs} with $n=0$. The function
\begin{align*}
u(x):=|x|^{-s_2}-|x|^{-s_1}
\end{align*}
satisfies $Lu=0$ and $|x|^{\alpha-2} u\in L^p (B)$,
  since
$
\alpha-2+{\rm Re\ }s_{1,2}>-N/p$.
Furthermore  $u=0$ on $\partial B$, hence  $u\in D_{p,\alpha}(B)$ and, since $Lu=0$,   Rellich inequalities fail.\\\qed

\section{The operator $A=|x|^2\Delta+cx\cdot \nabla$}\label{Section A}


Let  $c \in \R$ and
\begin{equation*} \label{defA}
A=|x|^2\Delta+cx\cdot \nabla.
\end{equation*}
This section is devoted to the  analysis of  $A$  acting on  $L^p(\Omega)$ for $1\leq p\leq \infty$, where $\Omega=\R^N$ or a bounded  domain, endowed with Dirichlet boundary conditions in this last case. The operator is degenerate both at $0$ and at $\infty$.  Employing spherical coordinates on $\R^N\setminus\{0\}$ we write $x=r\omega$, where $r:=|x|$, $\omega:=x/|x|\in {\mathbb{S}}^{N-1}$ and $$\Delta= D_{rr}+\frac{N-1}{r}D_r+\frac{1}{r^2}\Delta_0,$$ where $D_{rr}$, $D_r$ denote radial derivatives and $\Delta_0$ is  the Laplace-Beltrami operator on the unit sphere $\mathbb{S}^{N-1}$. Thus  we obtain
\begin{align*}
A=r^2D_{rr}+(N-1+c)r D_r+\Delta_0.
\end{align*}
Defining 
\begin{equation*} 
\Gamma=r^2D_{rr}+(N-1+c)r D_r,
\end{equation*}
the operators $\Gamma$ and $\Delta_0$ act on independent variables and therefore, when $\Omega$ is spherically symmetric,  generation and spectral properties  of $A$ can be proved through tensor products methods.

We start by  analysing $\Gamma$ and $\Delta_0$ separately and then we  deduce properties of  $A$ on $L^p(\Omega)$ when $\Omega=\R^N$ and $\Omega=B$. This method has the advantage to apply  also on more general  subspaces defined as tensor products of radial functions and spherical harmonics. Finally, we study $A$ in a general open set $\Omega$.

\subsection{The Laplace-Beltrami operator $\Delta_0$ on  $L^p_{J}(S^{N-1})$} 

We summarize in the next proposition some well known results about $\Delta_0$ referring, for example, to \cite{Grigoryan, Mor, SW} for further details. We recall that a spherical harmonic $P^n$ of order $n$ is the restriction to $\mathbb{S}^{N-1}$ of a homogeneous harmonic polynomial of degree $n$. We write $L^\infty (S^{N-1})$ for $C(S^{N-1})$.

\begin{prop}\label{lemma Spherical Harmonics}
The Laplace-Beltrami operator $\Delta_0$ generates an analytic semigroup $(T_{S^{N-1}}(t))_{t \ge 0}$ in $L^p(S^{N-1})$ (with respect to the surface measure $d\sigma$) for every $1\leq  p \le \infty$ . If $1<p<\infty$, its domain $D_p(\Delta_0)$ coincides with $W^{2,p}(S^{N-1}, d\sigma)$ . The spectrum of the operator $(\Delta_0,D_p(\Delta_0))$ is independent of $1 \le p \le \infty$ and consists of eigenvalues $-\lambda_n:=-n(n+N-2)$, $n \in \N_0$. 
The eigenspace corresponding to $-\lambda_n$  consist of all spherical harmonics of degree $n$ and has dimension $a_n$ where $a_0=1$, $a_1=N$ and for $n \ge 2$
$$
a_n= \binom{N+n-1}{n}-\binom{N+n-3}{n-2}.
$$
The linear span of spherical harmonics  coincides with the set of all polynomials and it is  dense in $C(\mathbb{S}^{N-1})$, hence in $L^p(\mathbb{S}^{N-1})$ for every $1\leq p<\infty$.
\end{prop}
{\sc{Proof.}} The generation and spectral properties of the Laplace Beltrami operator $\Delta_0$ are classic result about Heat operators on compact  manifolds. If $1<p<\infty$,  $D_p(\Delta_0)=W^{2,p}(S^{N-1}, d\sigma)$ by elliptic regularity. The analyticity of the semigroup as well as the invariance of the spectrum follows, for example,  from the Gaussian estimates of the  heat kernel of $\Delta_0$ (see e.g.  \cite[Theorem 5.2.1, Theorem 5.5.1]{davies}) using \cite[Corollary 7.5, Theorem 7.10]{ou}. The main properties of spherical harmonics can be found in \cite[Chapter II]{Mor} and \cite[Chapter IV.2]{SW}.\\\qed

Accordingly to the latter proposition let 
$$\sigma(S^{N-1})=\{\lambda_n=n(n+N-2):\ n\in\N_0\}$$
 be the spectrum of $(-\Delta_0,D_p(\Delta_0))$ and let us write  $\{P_j,\}_{j\in\N_0}$ and $\{\lambda(P_j)\}_{j\in\N_0}$  to denote  the sequences of the ($L^2$-orthonormal) eigenfunctions  and their  respectively eigenvalues repeated according to the relative  multiplicity. With this notation  $P_j$ is a spherical harmonics  whose eigenvalue is $\lambda(P_j)=n(n+N-2)$  and $n=\mbox{deg}(P_j)$.

We extend the analysis of $\Delta_0$ on more general  subspaces defined by  spherical harmonics.
\begin{defi} \label{Fjp}
 For a given $J \subseteq \N_0$ we define 
\begin{equation*}  
L^p_{J}(S^{N-1}) =\overline{span\{P_j: j \in J\}},
\end{equation*}
where the closure is taken in $L^p(S^{N-1})$, $1\leq p \le	\infty$.
\end{defi}
It is clear that $L^p_{J}(S^{N-1})$ is $\Delta_0$-invariant and that the domain of 
${\Delta_0}_{|L^p_{J}(S^{N-1})}$ is given by $D_{p}(\Delta_0)\cap L^p_{J}(S^{N-1})$. 
The following lemma is elementary and proved in \cite[Lemma 5.8]{rellich}.

\begin{lem}\label{spectrum Delta0}
Let $1\leq p\leq \infty$ and $J \subseteq \N_0$. Then ${\Delta_0}_{|L^p_{J}(S^{N-1})}$ generates  in $L^p_{J}(S^{N-1})$ the analytic semigroup $$\left(T_{S^{N-1}}(t)_{|L^p_{J}(S^{N-1})}\right)_{t \ge 0}.$$ 
Moreover ${span}\{P_j: j \in J\}$ is a core for $\Delta_0$ in $L^p_J(S^{N-1})$ and 
\begin{equation*}
\sigma({-\Delta_0}_{|L^p_{J}(S^{N-1})})=\{\lambda(P_j):\ j\in J \}
\end{equation*}
where $\lambda (P_j)$ is the eigenvalue whose eigenfunction is $P_j$.
\end{lem}

Note that, since each eigenvalue can have more than one eigenfunction, different set of indexes leads to different spaces but not necessarily to different spectra.

The asymptotic behaviour of $\left(T_{S^{N-1}}(t)_{|L^p_{J}(S^{N-1})}\right)_{t \ge 0}$  in $L^p_{J}(S^{N-1})$ is determined by the first eigenvalue. However we need a better estimate near $t=0$ which relies on  a
 Poincar\'e-type inequality.

\begin{lem}(\em{\cite[Lemma 2.7]{met-soba-spi3}})\label{Poincare-Lpn}
Let $1  <  p < \infty$ and $J \subseteq \N_0$ such that $n:=\mbox{min}J\geq 1$. Let  $\widetilde{C}_{p,n}$ be the best constant for which 
\[
\int_{S^{N-1}}|v|^p\,d\omega
\leq C\int_{S^{N-1}}|\nabla_{\!\tau} v|^2|v|^{p-2}\,d\omega, 
\quad v\in 
C^\infty(S^{N-1})\cap L^p_{J}(S^{N-1}).
\]
Then
$\widetilde{C}_{p,n}$ are finite, decreasing  and satisfy
$\widetilde{C}_{p,n}\to 0$ as $n\to \infty$.
\end{lem}

In the next Proposition we assume that the numbers $\lambda(P_j)$ are listed in the increasing order. 
\begin{prop} \label{asymptotic}
Let $J \subseteq \N_0$ and let $n$ be the smallest integer in $J$. There exists $M$ (depending on $n$ but not on $p$) such that for every $1 \le p \le \infty$
\begin{equation} \label{expdecay}
\|T_{S^{N-1}}(t)_{|L^p_{J}(S^{N-1})}\|_p \le M^{\big|1-\frac{2}{p}\big|}e^{-\lambda(P_n)\, t}.
\end{equation}
Furthermore $M=1$ when $n=0$. If $1<p<\infty$ then
\begin{equation} \label{expdecay 2}
\|T_{S^{N-1}}(t)_{|L^p_{J}(S^{N-1})}\|_p \le e^{-\frac{p-1}{\tilde C_{p,n}}\, t},
\end{equation}
where $\tilde C_{p,n}$ is the best constant of Lemma \ref{Poincare-Lpn}.
\end{prop}
{\sc{Proof.}} The first statement is proved in \cite[Lemma 5.9]{rellich}. To prove the second it is enough to show the dissipativity of $\Delta_0+\frac{p-1}{\tilde C_{p,n}}$ on $L^p_{J}(S^{N-1})$ or equivalently that, for every $u\in C^\infty(S^{N-1})\cap L^p_{J}(S^{N-1})$,
$$-\int_{S^{N-1}}\Delta_0u |u|^{p-2}u 
d\sigma\geq \frac{p-1}{\tilde C_{p,n}} \int_{S^{N-1}}|u|^p\;
d\sigma.$$
 Consider first the case $2\leq p<\infty$.
Setting
$u^\star=u|u|^{p-2}$ we multiply  $\Delta_0 u$  by
$u^\star$ and integrate over $S^{N-1}$. Integrating by parts and using Lemma \ref{Poincare-Lpn} we get
\begin{align*}
-\int_{S^{N-1}}\Delta_0 u\, u^\star\;
d\sigma&=(p-1)\int_{S^{N-1}}|u|^{p-2}|\nabla_\tau
u|^2\;d\sigma\geq \frac{p-1}{\tilde C_{p,n}}\int_{S^{N-1}}|u|^p\; d\sigma.
\end{align*}
For $1<p<2$ it is sufficient to replace $u^\star$ by $u(u^2+\delta)^{\frac{p}{2}-1}$,  $\delta>0$;  and then let $\delta$ to $0$ to obtain  the same inequality. \\\qed

\subsection{The operator $\Gamma$ on $L^p(I,r^{N-1}\,dr)$}\label{section Gamma}
In this section we summarize the main results about generation and spectral properties for the  operator 
$$\Gamma=r^2D_{rr}+(N-1+c)r D_r,$$
 acting, for $1\leq p<\infty$, on   $L^p(I,r^{N-1}dr)$, where $I=]0,\infty[$ or $I=]0,1[$. When  $p=\infty$,   $L^\infty(I,r^{N-1}\,dr)$  stands for the space $C_0^0(I)$ of all the continuous functions defined on $I$ vanishing at both endpoints.\\
For $1\leq p\leq \infty$  we define $\Gamma_p$ as the operator  $\Gamma$ endowed with the domain $D(\Gamma_p)$ defined, when $I=]0,\infty[$, as
\begin{align} \label{dp gamma 0,infty}
D(\Gamma_p)=\{u\in L^p(]0,\infty[,r^{N-1}\,dr),\ r\frac{\partial u}{\partial r},\ r^2\frac{\partial^2 u}{\partial r^2}\in L^p(]0,\infty[,r^{N-1}\,dr)\}
\end{align}
and 
for $I=]0,1[$
\begin{align} \label{dp gamma 0,1}
D(\Gamma_p)=\{u\in L^p(]0,1[,r^{N-1}\,dr),\ r\frac{\partial u}{\partial r},\ r^2\frac{\partial^2 u}{\partial r^2}\in L^p(]0,1[,r^{N-1}\,dr),\ u(1)=0\}.
\end{align}

In the next Theorem we show that $\Gamma_p$ always generates an analytic semigroup  in $L^p(I,r^{N-1}\,dr)$; the spectral analysis is more subtle since the spectrum  and the approximate point spectrum of $\Gamma_p$ drastically change accordingly to $I$ being bounded or not and to the sign of $N\left(1-\frac{2}{p}\right)-2+c$.

Let us introduce some notation: for $1\leq p\leq\infty$ (limiting values are taken for $p=\infty$), let us set
\begin{equation} \label{spettrogamma}
{\cal Q}_p:=\left\{\lambda\in \C\ \textrm{such that}\  {\rm Re}\lambda\leq -\frac{({\rm Im} \lambda)^2}{\left (N\left(1-\frac{2}{p}\right)-2+c\right)^2}-\omega_p\right\}
\end{equation}  and 
\begin{equation} \label{spettrogamma1}
{\cal P}_p:=\left\{\lambda=-\xi^2+i\xi\left (N(1-\frac{2}{p})-2+c\right)-\omega_p,\,  \xi\in\R\right\},
\end{equation}
where
\begin{equation}  \label{omegap}
\omega_p:=\frac{N}{p^2}\left[p(N-2+c)-N\right].
\end{equation}
${\cal P}_p$ is a parabola having vertex $-\omega_p$, symmetric with respect to the $x$ axis  
whereas ${\cal Q}_p$ is the  region enclosed inside  ${\cal P}_p$.  Obviously ${\cal P}_p$ coincides with the boundary of ${\cal Q}_p$ and, when $N\left(1-\frac{2}{p}\right)-2+c=0$, both  reduce to the half line $(-\infty, -\omega_p]$.

\begin{teo} \label{spec-rad}
Let $1\leq p\leq \infty$. Then the operator $\Gamma_p$  generates a strongly continuous  analytic semigroup $(S(t))_{t \ge 0}$  in $L^p(I,r^{N-1}dr)$ which satisfies the estimate 
$$\|S(t)\|_p \le e^{-\omega_p t},\quad  \text{ for }\ t \ge 0.$$
If $I=]0,\infty[$ we have
\begin{equation*} 
\sigma(\Gamma_p)=A\sigma(\Gamma_p)={\cal P}_p
.
\end{equation*}
If $I=]0,1[$, then
 \begin{equation*} 
\sigma(\Gamma_p)={\cal Q}_p.
\end{equation*}
Moreover
\begin{itemize}
\item[(i)] if $N\left(1-\frac{2}{p}\right)-2+c<0$, then \quad $\sigma(\Gamma_p)=A\sigma(\Gamma_p)={\cal Q}_p$,\; $P\sigma(\Gamma_p) \supset\overset{\mathrm{o}}{\cal Q}_p$;
\item[(ii)] if $N\left(1-\frac{2}{p}\right)-2+c=0$, then \quad  $\sigma(\Gamma_p)=A\sigma(\Gamma_p)=(-\infty, -\omega_p]$;
\item[(iii)]  if $N\left(1-\frac{2}{p}\right)-2+c>0$, then \quad $ A\sigma(\Gamma_p)={\cal P}_p$,\; $\overset{\mathrm{o}}{\cal Q}_p=R\sigma(\Gamma_p)\setminus A\sigma(\Gamma_p)$.
\end{itemize} 
\end{teo}
{\sc Proof.} Assume first that $I=]0,1[$. Let $J=]-\infty,0[$ and 
consider the isometry $S$ defined, for $1\leq p<\infty$, by 
$$S:L^p(J, ds)\to L^p( ]0,1[,{r}^{N-1}\,d{r}), \quad (Su)({r})={r}^{-\frac{N}{p}}u(\log {r}),$$
  and, for $p=\infty$, by
    $$S:C_0^0\left(J\right)\to C_0^0\left(]0,1[\right), \quad Su({r})=u(\log {r}).
$$ 
It follows that
$$S^{-1}\Gamma S u=u''+\left(N\left(1-\frac{2}{p}\right)-2+c\right)u'-\omega_p u.$$ 
By classical results, $S^{-1}\Gamma S$, endowed with  domain $D_p(S^{-1}\Gamma S)$
\begin{align*}
W^{2,p}(J)\cap W_0^{1,p}(J)\  ( p<\infty), \qquad
\left\{u\in C_0^0(J)\cap  C^2 \left(J\right):\ S^{-1}\Gamma Su\in C_0^0\left(J\right) \right\}\ ( p=\infty),
\end{align*}
 generates a strongly continuous analytic semigroup in $L^p\left(J\right)$ whose norm is bounded by $e^{-\omega_p t}$.

  It is elementary to check that
$$D(\Gamma_p)=\{Su:\ u\in D_p\left(S^{-1}\Gamma S\right)\}.$$
It follows that $\Gamma_p$  generates a strongly continuous and analytic semigroup $(S(t))_{t \ge 0}$  in the space $L^p(]0,1[, {r}^{N-1}d{r})$ which satisfies   $\|S(t)\|_p \le e^{-\omega_p t}$. The case $I=]0,\infty[$ is  similar and  proved in \cite[Proposition 5.1]{rellich} by considering $S$ with  $J=\R$.\\
 Concerning the second part of the statement we observe that the spectra of $\Gamma_p$ and $S^{-1}\Gamma_p S$ coincide.
 
  When $I=]0,\infty[$, the operator $S^{-1}\Gamma_p S$  is uniformly elliptic in $L^p(\R, ds)$, hence  its spectrum is independent of $p$ and coincides with the spectrum in $L^2(\R, ds)$ which is  ${\cal P}_p$, using the Fourier transform. Furthermore, since ${\cal P}_p$ coincides with its boundary, it follows, from Proposition \ref{boundary-spectrum},  that $\sigma(\Gamma_p)=A\sigma(\Gamma_p)={\cal P}_p$.
  
When $I=]0,1[$ we use  Lemma \ref{ODE2}  to see  that the spectrum  of  $S^{-1}\Gamma_p S$, hence of $\Gamma_p$,  coincides with the region ${\cal Q}_p$. Moreover, for the same reason,  the approximate point spectrum $A\sigma(\Gamma_p)$ coincides with ${\cal Q}_p$ if $N\left(1-\frac{2}{p}\right)-2+c<0$ (and in this case $P\sigma(\Gamma_p)\supset\overset{\mathrm{o}}{\cal Q}_p$), with the boundary ${\cal P}_p$ if $N\left(1-\frac{2}{p}\right)-2+c>0$ (and in this case $\overset{\mathrm{o}}{\cal Q}_p=R\sigma(\Gamma_p)\setminus A\sigma(\Gamma_p)$) and with  the half line $(-\infty, -\omega_p]$ when $N\left(1-\frac{2}{p}\right)-2+c=0$.\\
\qed

\begin{os}\label{Remark Equality domain N=1,p=inf}
Since the domain  $D_p(S^{-1}\Gamma S)$ coincides with its maximal one $$\{u \in L^p(J, ds): S^{-1}\Gamma Su \in L^p(J,ds)\},$$ as it easily follows from the classical interpolative inequalities
$
\|u'\|_p\leq \epsilon \|u''|_p+\frac{C}{\epsilon}\|u\|_p,
$ it follows that $$D(\Gamma_p)=\{ u \in L^p(I, r^{N-1}\, dr): \Gamma u \in L^p(I, r^{N-1}\, dr)\}.$$
 
\end{os}

\subsection{The operator $A=|x|^2\Delta+cx\cdot \nabla$ on $L^p_J(\R^N)$ and $L^p_J(B)$}
In this section   we use tensor arguments to combine the previous results on $\Gamma$ and $\Delta_0$ and deduce generation and spectral properties of  $$A=|x|^2\Delta+cx\cdot \nabla$$ on $L^p(\Omega)$ when $\Omega=\R^N$ and $\Omega=B$. We extend the analysis  also on more general  subspaces defined by tensor products of radial functions and spherical harmonics.

If $X,Y$ are function spaces over $G_1, G_2$ we denote by $X\otimes Y$ the algebraic tensor product of $X,Y$, that is the set of all functions $u(x,y)=\sum_{i=1}^n f_i(x)g_i(y)$ where $f_i \in X, g_i \in Y$ and $x \in G_1, y\in G_2$.
If $T,S$ are linear operators  on $X,Y$ we denote by $T\otimes S$ the operator on $X\otimes Y$ defined by 
$$
T\otimes S \left (\sum_{i=1}^n f_i(x)g_i(y)\right)=\sum_{i=1}^n T f_i(x)Sg_i(y).
$$
Let us fix a  complete orthonormal system of spherical harmonics $\{P_j,\}_{j\in\N_0}$  $L^2(S^{N-1})$ and let $\{\lambda(P_j)\}_{j\in\N_0}$  be  the sequence of the  corresponding eigenvalues  repeated according to their multiplicity.  With this notation  $-\Delta_0(P_j)=\lambda(P_j)P_j$ and $\lambda(P_j)=n(n+N-2)$, where  $n=\mbox{deg}(P_j)$. \\

Unless otherwise specified     $\Omega$ denotes $\R^N$ or $B$,   $I$ stands for $]0,\infty[$, $]0,1[$, respectively. As usual we write $L^\infty (\Omega)$ for $C_0^0(\Omega)$. 
\begin{defi} \label{Fjp1}
Let $1\leq p\leq \infty$ and let $J \subseteq \N_0$. We define 
\begin{equation*}  
L^p_{J}(\Omega) =\overline{L^p\left(I, r^{n-1}dr\right)\otimes L^p_J(S^{N-1})}=\overline{L^p\left(I, r^{n-1}dr\right)\otimes \mbox{span}\{P_j: j \in J\}},
\end{equation*}
where the closure is taken in $L^p(\Omega)$.   Fixing $n \in \N_0$ we write  $L^p_{\ge n}(\Omega), L^p_n(\Omega), L^p_{<n}(\Omega)$  when $J$ identifies all spherical harmonics of order $\ge n$, $n$ and $< n$ respectively. The spaces $L^p_{>n}(\Omega), L^p_{\le n}(\Omega), L^p_{\neq n}(\Omega)$ are defined similarly.
\end{defi}

Note that  $L^p_{J}(\Omega)=L^p(\Omega)$ if $J=\N_0$.

\smallskip

The next lemma clarifies the structure of the spaces  $L^p_{J}(\Omega)$.

\begin{lem} \label{projection}
Assume that the $L^2$ orthogonal projection $P: L^2(S^{N-1}) \to L^2_J(S^{N-1})$ extends to a bounded projection $P$ in $L^p(S^{N-1})$.  Then 
\begin{equation} \label{complement}
L^p(\Omega)= L^p_{J}(\Omega) \oplus L^p_{\N_0\setminus J}(\Omega)
\end{equation}
and
\begin{equation} \label{caratterizzazione}
 L^p_{J}(\Omega)=\left \{u \in L^p(\Omega): \int_{S^{N-1}} u(r\, \omega)P_j(\omega) \, d\sigma (\omega)=0\ {\rm for} \ r\in I\ {\rm  and}\  j \not \in J\right \}.
\end{equation}
When $J$ is finite  
\begin{equation} \label{Jfinito1}
L^p_J(\Omega)=\Bigl \{ u=\sum_{j \in J}f_j(r)P_j(\omega): f_j \in L^p(I, r^{N-1}dr)\Bigr \}
\end{equation} and 
the projection $I\otimes P :L^p(\Omega) \to L^p_{J}(\Omega)$ is given by
\begin{equation} \label{Jfinito}
(I\otimes P) u=\sum_{j \in J}T_j u (r)\, P_j(\omega),
\end{equation}
where $$T_j u (r):=\int_{S^{N-1}} u(r\, \omega)P_j(\omega) \, d\sigma (\omega), \quad  \forall\ u\in L^p(\Omega).$$
\end{lem}
{\sc{Proof.}} 
When $\Omega=\R^N$ we refer  to \cite[Lemma 5.11]{rellich}. The proof for $\Omega=B$ is identical.
\\\qed

\begin{os} (i) The equality  
$$ L^p_{J}(\Omega) = \left \{u \in L^p(\Omega): \int_\Sigma u(r\, \omega)P_j(\omega) \, d\sigma (\omega)=0\ {\rm for} \ r\in I\ {\rm  and}\  j \not \in J\right \}
$$
holds without assuming the boundedness of the projection $P$ (see \cite[Proposition 2.8]{rellich Disc}).\\[1ex]
(ii)  $L^p_0(\Omega)$ consists of radial functions and  $L^p(\Omega)=L^p_{\le n}(\Omega)\oplus L^p_{>n}(\Omega)$.
\end{os}

The following result follows from well-known and elementary facts about  Tensor Product Semigroups, see \cite[AI, Section 3.7]{nagel}. A proof is provided in \cite[Proposition 5.14]{rellich} when $\Omega=\R^N$, the case of the ball is similar.

\begin{prop} \label{analyt}
For $1\le p \le \infty$, let  $D(\Gamma_p)$ and $D({\Delta_0}_{|L^p_{J}(S^{N-1})})$ be the domains of $\Gamma_p$ and ${\Delta_0}_{|L^p_{J}(S^{N-1})}$ introduced in the previous subsection. Then the closure of the operator 
$$\left(A,\, D(\Gamma_p)\otimes D({\Delta_0}_{|L^p_{J}(S^{N-1})})\right)$$ generates a strongly continuous  analytic semigroup $(T_{p,J}(t))_{t \ge 0}$
in $L^p_J(\Omega)$.  Let $n$ be the smallest integer in $J$. Then there exists $M$ (depending on $n$ but not on $p$) such that for every $1 \le p \le \infty$
\begin{equation} \label{expdecay Omega}
\|T_{p,J}(t)\|_p \le M^{\big|1-\frac{2}{p}\big|}e^{-(\omega_p+\lambda(P_n))\, t},
\end{equation}
where $\omega_p$ is defined in \eqref{omegap} and $M$ is the constant   in \eqref{expdecay} which satisfies $M=1$ when $n=0$.  Moreover, if $1<p<\infty$, then
\begin{equation} \label{expdecay 2 Omega}
\|T_{p,J}(t)\|_p \le e^{-\left(\omega_p+\frac{p-1}{\tilde C_{p,n}}\right)\, t},
\end{equation}
where $\tilde C_{p,n}$ is the best constant of Lemma \ref{Poincare-Lpn}.
\end{prop}
%

\begin{defi}\label{Def A_pJ}
We denote by $A_{p,J}$ the closure of  $(A,D(\Gamma_p) \otimes D({\Delta_0}_{|L^p_{J}(S^{N-1})}))$ 
in $L^p_{J}(B)$. When $J=\N_0$ we write $A_{p}$ for $A_{p,J}$ and $T_{p}(t)$ for $T_{p,J}(t)$.
\end{defi}

 The proof of the following corollary is immediate.

\begin{cor} \label{restriction}
$T_{p,J} (t)$ is the restriction of $T_{p}(t)$ to $L^p_J(B)$ and its generator $A_{p,J}$ is the part of $A_{p}$ in $L^p_J(B)$.
\end{cor}

As in \cite[Proposition 5.16]{rellich}, we prove that the smooth functions are a core for $A_{p,J}$.

\begin{prop}  \label{core Palla R^N}
Let $1\leq p\leq\infty$. The set 
$$
C_{c,0}^2\left(B\right):=\big\{u\in C^2_c\left(\bar B\right): u=0  \text{ on } \partial B\ \text{and on a neighborhood of } 0\big\}.
$$
 is a core for $A_{p,J}$ when $\Omega=B$. When $\Omega=\R^N$, $C^\infty _c\left(\R^N\setminus\{0\}\right)$ is a core for $A_{p,J}$.
\end{prop}
{\sc Proof.} Let us suppose  that $\Omega=B$. Recalling the proof of Theorem \ref{spec-rad}, we observe that,  since by Proposition \ref{Sobolev approximation 1,infty} the set 
$$\big\{u\in C^2_c\left(]-\infty,0]\right):\ u(0)=0 \big\}$$
 is dense  in $D_p(S^{-1}\Gamma_p S)$, then  
 $${\cal F}:=\big\{u\in C^2_c\left(]0,1]\right):\ u(1)=0\big\}$$
  is  dense in $D(\Gamma_p)$. Moreover $ \mbox{span}\{P_j: j \in J\}$ is dense in $D({\Delta_0}_{|L^p_{J}(S^{N-1})})$. Since by construction $D(\Gamma_p)\otimes  D({\Delta_0}_{|L^p_{J}(S^{N-1})})$ is a core for $A_{p,J}$, it follows that 
$$\mathcal{F}\otimes  \mbox{span}\{P_j: j \in J\}$$
is dense in $D(A_{p,J})$. Observing that 
$$\mathcal{F}\otimes \mbox{span}\{P_j: j \in J\}\subseteq C_{c,0}^2\left(\Omega\right)$$
 we get the thesis. The proof for $\Omega=\R^N$ is similar.\\\qed

In order to prove the main result of this section, namely 
$$\sigma(A_{p,J})=\sigma(\Gamma_p)+\sigma({\Delta_0}_{|L^p_{J}(S^{N-1})}),$$
we need two preliminary lemmas. The first provides some regularity properties of the projection defined in \eqref{Jfinito} and is proved in \cite[Lemma 2.15]{met-soba-spi3} when $\Omega=\R^N$.

\begin{lem}\label{regularity of projections}
Let $J\subseteq \N_0$ and let $j_0\in J$.   Let us consider the operator $T_{j_0} :L^p_J(\Omega) \to L^p(I,\ r^{N-1}dr)$ defined by 
$$T_{j_0} u (r):=\int_{S^{N-1}} u(r\, \omega)P_{j_0}(\omega) \, d\sigma (\omega), \quad  \forall\ u\in L^p(\Omega)$$
and the projection 
$$I\otimes P_{j_0} :L^p_J(\Omega) \to L^p_{j_0}(\Omega)=L^p(I,\ r^{N-1}dr)\otimes P_{j_0}$$ given, for   $u\in L^p_J(\Omega)$, $r\in I$, $\omega\in S^{N-1}$, by
\begin{equation*}
(I\otimes P_{j_0})\ u(r\omega)=\,T_{j_0} u (r) P_{j_0}(\omega).
\end{equation*}
Then $T_{j_0}$, $I\otimes P_{j_0}$ are well defined and bounded operator. Furthermore $T_{j_0}$ maps $D(A_{p,J})$ onto $D(\Gamma_p)$ and one has
\begin{align}\label{A Tu}
T_{j_0}Au=\Big(\Gamma-\lambda(P_{j_0})\Big)T_{j_0} u, \quad \forall u\in D(A_{p,J}).
\end{align}
\end{lem}

The next lemma relates the spectra of $\Gamma_p$ and $A_{p,J}$.

\begin{lem}\label{heritability spectrum}
Let $1\leq p\leq \infty$, $J\subseteq\N_0$ and $j_0\in J$. Let $\Omega$ stand for $\R^N$ or $B$ and let $A_{p,J}$ be the operator defined in Definition \ref{Def A_pJ}. The following properties hold.
\begin{itemize}
\item[(i)]If $\lambda\in P\sigma(\Gamma_p)$ then $\lambda-\lambda(P_{j_0})\in P\sigma(A_{p,J})$;
\item[(ii)] If $\lambda\in A\sigma(\Gamma_p)$ then $\lambda-\lambda(P_{j_0})\in A\sigma(A_{p,J})$;
\item[(iii)] If $\lambda\in R\sigma(\Gamma_p)$ then $\lambda-\lambda(P_{j_0})\in R\sigma(A_{p,J})$;
\end{itemize}
\end{lem}
{\sc{Proof.}}
Let $\lambda\in P\sigma(\Gamma_p)$ and let $0 \neq u\in D(\Gamma_p)$ be  such that $\Gamma u=\lambda u$. Then  it is immediate to see  that the function $f=uP_{j_0}$ satisfies $f\in D(A_{p,J})$ and $Af=\left(\lambda-\lambda(P_{j_0})\right)f$. This proves (i).\\
Assertion (ii) follows similarly by using Lemma \ref{Char Aspectrum}.\\
Let us now consider  (iii) and let $\lambda\in R\sigma(\Gamma_p)$. Recalling Definition \ref{defi R spectrum} we have to show that   $\mbox{rg}\left(\lambda-\lambda(P_{j_0})-A_{p,J}\right)$ is not dense in $L^p(A_{p,J})$. 
 Since $\lambda\in R\sigma(\Gamma_p)$,   $\mbox{rg}(\lambda-\Gamma_p)$ is not dense in $L^p(I,r^{N-1}dr)$ and therefore  there exists a  linear form $0\neq G$ in the dual space $\left(L^p(I,r^{N-1}dr)\right)'$ which vanishes over $\mbox{rg}(\lambda-\Gamma_p)$.  Let us consider the projection
\begin{align*}
T_{j_0}:L^p_J(\Omega)\to L^{p}(I,r^{N-1}dr),\quad u\mapsto T_{j_0}u(r)=\int_{S^{N-1}}u(r\omega)P_{j_0}(\omega)\,d\sigma(\omega).
\end{align*}
Using  Lemma  \ref{regularity of projections} we see that $0\neq T=G\circ T_{j_0}$ belongs to  the dual space $\left(L^p_j(\Omega)\right)'$ and satisfies for $u\in D(A_{p,J})$,
\begin{align*}
T\left(\lambda-\lambda(P_{j_0})-A\right)u=G\Big(T_{j_0}\left(\lambda-\lambda(P_{j_0})-A\right)u\Big)=G\Big(\left(\lambda-\Gamma_p \right)T_{j_0}u\Big)=0.
\end{align*}
This implies that $T$  vanishes over $\mbox{rg}(\lambda-\lambda(P_{j_0})-A_{p,J})$ and proves (iii).\\\qed

We can finally describe in detail the spectrum of $A_{p,j}$.  We are mainly interested in the computation of the complement of the approximate point spectrum, that is the set of all $\lambda$ such that the inequality
\begin{align*}
\|u\|\leq C\|\lambda u-Au\|,\quad \forall u\in D(A_{p,J})
\end{align*}
holds,
since it is equivalent to Rellich inequalities. Observe that the situation is more complicate in the case where $N\left(1-\frac{2}{p}\right)-2+c>0$ since residual spectra appear.

We recall that ${\cal P}_p$ and ${\cal Q}_p$ are defined in \eqref{spettrogamma} and \eqref{spettrogamma1}.

\begin{teo} \label{Spectrum main} 
Let $1\leq p\leq \infty$, $J\subseteq\N_0$ and $j_0:=\min\{j\in J\}$. The following properties hold 

\begin{itemize}
\item[1.]
If $\Omega=\R^N$, the  spectrum of $A_{p,J}$ in $L^p_J(\R^N)$ is given by
\begin{equation*} 
\sigma(A_{p,J})=A\sigma(A_{p,J})=\bigcup\limits_{j\in J}({\cal P}_p-\lambda(P_j))
\end{equation*}
and reduces to $]-\infty,-\omega_p-\lambda(P_{j_o})]$ when $N\left(1-\frac{2}{p}\right)-2+c=0$. 
\item[2.] If $\Omega=B$, the  spectrum of $A_{p,J}$ in $L^p_J(B)$
 is given by
 \begin{equation*} 
 \sigma(A_{p,J})={\cal Q}_p-\lambda (P_{j_0})
\end{equation*}
and reduces to $]-\infty,-\omega_p-\lambda(P_{j_o})]$ when $N\left(1-\frac{2}{p}\right)-2+c=0$. In particular we have \medskip
\begin{itemize}

\item[(i)] If $N\left(1-\frac{2}{p}\right)-2+c<0$, then $$A\sigma(A_{p,J})={\cal Q}_p-\lambda (P_{j_0}),\quad P\sigma(A_{p,J}) \supset\overset{\mathrm{o}}{\cal Q}_p-\lambda(P_{j_0}).$$

\item[(ii)] If $N\left(1-\frac{2}{p}\right)-2+c=0$, then  $$A\sigma(A_{p,J})=(-\infty, -\omega_p-\lambda (P_{j_0})].$$

\item[(iii)]  If $N\left(1-\frac{2}{p}\right)-2+c>0$, then 

\begin{align*}
A\sigma(A_{p,J})&=\bigcup\limits_{j\in J}({\cal P}_p-\lambda (P_{j}));\\[1ex]
R\sigma(A_{p,J})\setminus A\sigma(A_{p,J})&=\left (\overset{\mathrm{o}}{\cal Q}_p-\lambda(P_{j_0})\right )\setminus \bigcup\limits_{j\in J}({\cal P}_p-\lambda (P_{j})).
\end{align*}

\end{itemize}
\end{itemize}
\end{teo}
{\sc Proof.} 
We give a proof only when $\Omega=B$, since the case $\Omega=\R^N$ is similar and proved in  \cite[Theorem 5.17]{rellich}.
Let us prove first the inclusion 
\begin{align*}
\sigma(A_{p,J})\subseteq 
\sigma(\Gamma_p)+\sigma({\Delta_0}_{|L^p_{J}(S^{N-1})})
={\cal Q}_p-\lambda (P_{j_0}).
\end{align*}
Let $\lambda \not \in {\cal Q}_p-\lambda (P_{j_0})$ and fix $n \in \N_0$ such that
\begin{equation} \label{boundgamma}
-\omega_p-\lambda(P_k) < {\rm Re}\, \lambda \quad {\rm for\  every\ } k > n.
\end{equation} 
According to Lemma \ref{projection} we write
$L^p_J(B)= L^p_{J_n}(B) \oplus L^p_{J\setminus J_n}(B)$,
where $J_n=J\cap \{0,1,\dots ,n\}$ (note that if $J_n=\emptyset$ then  $L^p_{J_n}(B)={0}$ and $L^p_J(B)\subseteq L^p_{> n}(B)$). Since both  $L^p_{J_n}(B)$  and  $L^p_{J\setminus J_n}(B)$ are $A_{p,J}$ invariant, then $\lambda \in \rho(A_{p,J})$ if and only if  $\lambda \in \rho(A_{p,J_n})$ and  $\lambda \in \rho(A_{p,J \setminus J_n})$. The second inclusion follows immediately from \eqref{expdecay Omega} with $J\setminus J_n$ instead of $J$, since $\rm{Re}\, \lambda$ is greater than the growth bound of $(T_{p, J \setminus J_n})_{t \ge 0}$, by  (\ref{boundgamma}). Concerning the first inclusion let us suppose that $J_n\neq\emptyset$ and, without loss of generality, let us assume $J_n=\{0,1,\dots,n\}$. We note that
$$
 L^p_{J_n}(B) =\oplus_{i=0}^n  L^p_{i}(B) =\oplus_{i=0}^n  L^p\left((0,1),r^{N-1}dr\right)\otimes P_i 
$$
 and that each $ L^p_{i}(B) $ is $A_{p,J}$ invariant. Moreover, $\lambda-A_{p,J}$ coincides with $\left (\lambda+\lambda (P_i)-\Gamma_p\right ) \otimes I$ on $ L^p_{i}(B) $, hence it is invertible on it, since $\lambda+\lambda(P_i) \not \in {\cal Q}_p=\sigma(\Gamma_p)$ by assumption.
This shows that $\lambda \in \rho (A_{p,J})$, hence
\begin{align}\label{Spectrum palla eq 1}
\sigma(A_{p,J})\subseteq {\cal Q}_p-\lambda (P_{j_0}).
\end{align}

Let us prove the opposite inclusion. Using the description of the spectrum of $\Gamma_p$  proved in Theorem \ref{spec-rad} and Lemma \ref{heritability spectrum}, we get immediately the reverse inclusion and (i) and (ii).

  In the case $N\left(1-\frac{2}{p}\right)-2+c>0$, Lemma \ref{heritability spectrum} only shows that 
\begin{align*}
A\sigma(A_{p,J})&\supseteq\bigcup\limits_{j\in J}({\cal P}_p-\lambda (P_{j}));\\[1ex]
R\sigma(A_{p,J})&\supseteq \left (\overset{\mathrm{o}}{\cal Q}_p-\lambda(P_{j_0})\right )\setminus \bigcup\limits_{j\in J}({\cal P}_p-\lambda (P_{j})).
\end{align*}
To end the proof we need to show that, if $\lambda\in \left (\overset{\mathrm{o}}{\cal Q}_p-\lambda(P_{j_0})\right )\setminus \bigcup\limits_{j\in J}({\cal P}_p-\lambda (P_{j}))$, then $\lambda\notin A\sigma(A_{p,J})$. Recalling Proposition \ref{Rellich-spectrum} this is equivalent to the validity, for some $C>0$, of the inequality
\begin{align}\label{Spectrum palla eq 2}
\|\lambda v-Av\|_p \geq C \|v\|_p,\quad \forall v\in D(A_{p,J}).
 \end{align}
 Let us fix  $\lambda\in \left (\overset{\mathrm{o}}{\cal Q}_p-\lambda(P_{j_0})\right )\setminus \bigcup\limits_{j\in N_0}({\cal P}_p-\lambda (P_{j}))$ and  let $\ov{n}\in N_0$ sufficiently large such that   $\lambda\notin {\cal Q}_p-\lambda(P_{\ov n})$.
Then, by \eqref{Spectrum palla eq 1}, $\lambda$ belongs to the resolvent of the operator $A_{p,>\ov{n}}$ in $L^p_{>\ov{n}}(B)$. It follows that \eqref{Spectrum palla eq 2} is true in $L^p_{>\ov{n}}(B)$.\\
Since from \eqref{complement}, $L^p(B)= L^p_{\leq \ov n}(B) \oplus L^p_{>\ov n}(B)$, it remains  to  proves \eqref{Spectrum palla eq 2} for any  $v\in D(A_{p,J})\cap L^p_{\leq \ov n}(B)$. Recalling \eqref{Jfinito1} and Lemma \ref{regularity of projections}, one has  
\begin{align*}
 v(\rho\omega)=\sum_{i=1}^{\ov{n}}c_i(\rho)P_i(\omega),
 \end{align*}
for some $c_i\in D(\Gamma_p)$.  Then
\begin{align*}
\|\lambda v-Av\|_p &=\|\sum_{i=1}^{\ov{n}}P_i\left(\lambda+\lambda(P_i)-\Gamma\right)c_i\|_p\geq C \sum_{i=1}^{\ov{n}}\|P_i\left(\lambda+\lambda(P_i)-\Gamma\right)c_i\|_p\\[1ex]
&=C \sum_{i=1}^{\ov{n}}\|\left(\lambda+\lambda(P_i)-\Gamma\right)c_i\|_{L^p\left((0,1), r^{N-1}dr\right)},
\end{align*}
where in the last equality we have used spherical coordinates to evaluate the integrals.\\
By the assumption on $\lambda$ and recalling (iii) in Theorem \ref{spec-rad}, one has $\lambda+\lambda(P_i)\notin {\cal P}_p=A\sigma(\Gamma_p)$, which implies, for a possibly different constant $C>0$,
\begin{align*}
\nonumber \|\lambda v-Av\|_p &\geq C \sum_{i=1}^{\ov{n}}\|\left(\lambda+\lambda(P_i)-\Gamma\right)c_i\|_{L^p\left((0,1), r^{N-1}dr\right)}\\[1ex]
&\geq C \sum_{i=1}^{\ov{n}}\|c_i\|_{L^p\left((0,1), r^{N-1}dr\right)}\geq C \|v\|_p.
\end{align*}
This proves \eqref{Spectrum palla eq 2} in the remaining case.\\\qed

\begin{os} {\rm 
The inclusion 
\begin{align*}
\sigma(A_{p,J})\subseteq 
\sigma(\Gamma_p)+\sigma({\Delta_0}_{|L^p_{J}(S^{N-1})})=
{\cal Q}_p-\lambda (P_{j_0}).
\end{align*}
follows also from the more general result  \cite[Theorem 7.3]{arendt1} since the semigroups generated by $\Gamma$ and ${\Delta_0}_{|L^p_{J}(S^{N-1})}$ are analytic and commute.
 }
\end{os}

\begin{cor} \label{a-priori}
 Let $\Omega$ be equal to $\R^N$ or $B$ and assume that $\lambda+\omega_p>0$. Then the best constant for which the inequality 
\begin{equation} \label{estA}
\|u\|_p\leq C\|\lambda u-A u\|_p,\quad \forall u\in D(A_p)
\end{equation} 
holds is given by 
\begin{align*}
C=\frac{1}{\lambda+\omega_p}.
\end{align*}
\end{cor}
{\sc{Proof.}}
If $\lambda+\omega_p>0$, then  $\lambda\in\rho(A_p)$, by the preceding theorem, and  then the optimal constant in \eqref{estA} is $
\|R(\lambda,A_p)\|_p$. Recalling \eqref{distance parabola} we have
\begin{align*}
\|R(\lambda,A_p)\|_p\geq \frac{1}{\mbox{dist}(\lambda,\sigma(A_p))}=\frac{1}{\mbox{dist}(\lambda,{\cal{P}}_p))}=\frac{1}{\lambda+\omega_p}.
\end{align*}
Using the contractivity estimates \eqref{expdecay Omega} and writing the resolvent as the Laplace transform of the semigroup we see that also  the reverse inequality
\begin{align*}
\|R(\lambda,A_p)\|_p\leq \frac{1}{\lambda+\omega_p}
\end{align*}
holds.
\qed

\subsection{The operator $A=|x|^2\Delta+cx\cdot \nabla$ on $L^p(\Omega)$}
In this section we complete the study of the operator $A$ in $\R^N$ and $B$ by providing a complete description of the domain. Then we use the results in the whole space to extend our results to bounded sets containing the origin. In particular we prove that the domain of the operator coincides with the maximal one, see Proposition \ref{L1}. This allows to state the precise class of functions where Rellich inequalities hold. Note that $A$ is singular both at $0$ and at $\infty$.

 Let $\beta\in(0,1]$. In what follows we assume $\Omega$ to be $\R^N$ or  a  bounded open connected  subset of $\R^N$ whose boundary $\partial \Omega$ is $C^{2,\beta}$ and such that $0\notin\partial\Omega$. 
 For any $p\in ]1,\infty[$ we define $A_p$ by $A_p=Au$ in
\begin{align} \label{dp}
D_p(\Omega)&=\left\{u\in W^{2,p}(\Omega\setminus B_\epsilon)\cap L^p(\Omega)\ \forall \epsilon>0\,:\,u=0 \text{ on } \partial \Omega,\  |x|\nabla
u,\ |x|^2D^2 u\in L^p(\Omega)\right\};
\end{align}
for $p=1$ we define $\left(A_1,D_1(\Omega)\right)$ as
\begin{align}\label{d1}
D_1(\Omega)&=\left\{u\in L^1(\Omega):\,u=0 \text{ on } \partial \Omega,\  |x|\nabla
u,\ |x|^2\Delta u\in L^1(\Omega)\right\}.
\end{align}
When $\Omega=\R^N$ and, correspondingly,  $\partial\Omega=\emptyset$,  the  requirement "$u=0$ on  $\partial \Omega$" must be disregarded. When $\Omega$ is bounded  the Dirichlet boundary condition $u(x)=0$ for $x\in\partial\Omega$ makes sense in the sense of traces since $u$ has first derivatives in $L^p$  in a neighbourhood of the boundary $\partial \Omega$. The case $0\notin\Omega$ is classical since the term $|x|$ is negligible and, for $1<p<\infty$, $D_p(\Omega)$ becomes $W^{2,p}(\Omega)\cap W_0^{1,p}(\Omega)$. \\
For $p=\infty$, we also consider the operator $A_\infty$ endowed with the domain
\begin{align}\label{def A Omega infty}
D_\infty(\Omega)&=\left\{u\in C_0^0(\overline{\Omega}):\ Au\in C_0^0(\overline{\Omega}),\  |x|\nabla
u,\ |x|^2\Delta u\in C^0(\overline{\Omega})\right\},
\end{align}
where   $C^0(\overline{\Omega})$ denotes the space of bounded and continuous  functions defined in $\overline{\Omega}$ and  vanishing at the origin, if $0\in\Omega$; $C_0^0(\overline{\Omega})$ is its subspace consisting of functions vanishing also at $\infty$ when $\Omega=\R^N$ and at  the boundary $\partial\Omega$, otherwise.
\medskip

When $\Omega$ is bounded we use Proposition \ref{Partition unity} to fix $\delta>0$ such that the subsets
\begin{align*}
K_{\delta}:=\left\{x\in\R^N:\ \mbox{dist}(x,\partial\Omega)<\delta\right\},\quad \Omega_{\delta}:=K_{\delta}\cap\Omega
\end{align*}
have $C^{2,\beta}$ boundary. Furthermore we can write $\overline\Omega=\overline\Omega_\delta\cup\Omega_0$ where $\Omega_0$ is an open subset $\Omega_0\subset\subset\Omega$   and we fix a partition of  unity $\{\eta_\delta^2,\eta_0^2\}$ such that

\begin{align}\nonumber
(i)&\quad \eta_\delta\in C_c^\infty(K_{\delta}), \quad 0\leq\eta_\delta\leq 1,\quad \eta_\delta=1 \;\text{ in }\; \overline\Omega_{\frac \delta 2};\\[1ex]\label{Partiton unity eq}
(ii)&\quad\eta_0\in C_c^\infty(\Omega_0), \quad 0\leq\eta_\delta\leq 1;\\[1ex]\nonumber
(iii)&\quad \eta_\delta^2+\eta_0^2=1\hspace{5ex} \text{ in }\;\overline\Omega.
\end{align}

\bigskip
In order to  identify a core for $A_p$ we define
\begin{align*}
C_{c,0}^2\left(\Omega\right):&=\big\{u\in C^2_c\left(\bar\Omega\setminus\{0\}\right):\ u=0  \text{ on } \partial \Omega\big\}\\
&=\big\{u\in C^2_c\left(\bar\Omega\right):\ u=0  \text{ on } \partial \Omega\ \text{and on a neighborhood of } 0\big\}.
\end{align*}

\smallskip

\begin{lem}\label{coreOmega}
The space $C_{c,0}^2\left(\Omega\right)$ is dense in $D_p(\Omega)$, endowed with
the norm
\begin{align*}
\|u\|_{D_p(\Om)}&=\|u\|_p+\||x|\nabla u\|_p+\||x|^2D^2 u\|_p,\quad (1<p<\infty);\\[1.5ex]
\|u\|_{D_p(\Om)}&=\|u\|_p+\||x|\nabla u\|_p+\||x|^2\Delta u\|_p,\quad\hspace{2ex} (p=1,\infty).
\end{align*}
When $\Omega=\R^N$, $C^\infty _c\left(\R^N\setminus\{0\}\right)$ is is dense in $D_p(\Omega)$.
 \end{lem}
{\sc Proof.} Let us consider, preliminarily,  $\Omega=\R^N$.\\
 Let $u\in D_p(\R^N)$; we approximate $u$ with functions in
$D_p(\R^N)$ having compact support in $\R^N\setminus\{0\}$. Let
$$\Omega_n=\left\{x\in\R^N:\ |x|\geq \frac{1}{n}\right\},\quad \xi_n=\chi_{\Omega_{\frac n 2}}\ast\phi_\frac{1}{n}$$ where $\phi$ is a classical mollifier supported in $B_1$, with $\int_{\R^N}\phi=1$
 and $\phi_{\frac{1}{n}}(x)=n^N\phi\left(nx\right)$. It is easy to check that $\xi_n(x)=1$ for $x\in\Omega_n$, $\xi_n$ is supported in $\R^N\setminus\{0\}$
  and that $0 \le \xi_n \le 1$, $|\nabla \xi_n|\leq Cn$, $|D^2\xi_n|\leq Cn^2$. Consider also a smooth function $\eta$ such that $\chi_{B_1}\leq\eta\leq \chi_{B_2}$
  and, for every $n\in\N$, define $\eta_n(x)=\eta\left(\frac{x}{n}\right)$.
Set $u_n=\xi_n\eta_n u$. It is immediate to check, using Lebesgue's Theorem, that $u_n$ tends
to $u$ in $L^p(\R^N)$. Concerning the gradient term, we have
\begin{align*}
 \||x|(\nabla(\xi_n\eta_n u)-\nabla u)\|_p^p\leq&
 \int_{\R^N}|x|^p|\xi_n\eta_n-1|^p|\nabla u|^p\,dx\\[1ex]
 &+
 \int_{\R^N}|x|^p|\nabla\xi_n|^p|\eta_n|^p|u|^p\,dx+
 \int_{\R^N}|x|^p|\xi_n|^p|\nabla\eta_n|^p |u|^p\,dx\\[1ex]
 \leq&\int_{\R^N}|x|^p|\xi_n\eta_n-1|^p|\nabla u|^p\,dx\\[1ex]
 &+C n^p
 \int_{|x|\leq\frac 1 n}|x|^p |u|^p\,dx+
 C  n^{-p}\int_{\{n \le |x| \le
 2n\}}|x|^p|u|^p\,dx.
\end{align*}
The last inequality implies
\begin{align*}
 \||x|(\nabla(\xi_n\eta_n u)-\nabla u)\|_p^p\,dx\leq&  \int_{\R^N}|x|^p|\xi_n\eta_n-1|^p|\nabla u|^p\,dx\\[1ex]&
 +C\int_{|x|\leq\frac 1 n}|u|^p\,dx+C\int_{\{n \le |x|\le 2n\}} |u|^p\,dx
\end{align*}
which tends to 0 by dominated convergence. Using a similar argument
one shows that, if $1<p<\infty$, $|x|^2D^2u_n$ tends to $|x|^2D^2 u$ in $L^p(\R^N)$ and that, if $p=1,\infty$,   $|x|^2\Delta u_n$ tends to $|x|^2\Delta u$ in $L^p(\R^N)$. This proves that $u_n$ tends to $u$ in $D_p(\R^N)$; we also note that, by construction, $\mbox{supp\,} u_n\subseteq \mbox{supp\,} u$.  Finally
we can use a standard convolution argument to approximate in $D_p(\R^N)$ functions having compact support in $\R^N\setminus\{0\}$ with $C_c^\infty\left(\R^N\setminus\{0\}\right)$
functions.\\
Let us consider, now, a bounded set   $\Omega\subset\R^N$ and let $u\in D_p(\Omega)$. We use the partition of unity defined in \eqref{Partiton unity eq} to write
\begin{align*}
u=\eta_0^2u+\eta_\delta^2 u:=u_0+u_\delta.
\end{align*}
The function  $u_0$ satisfies $\mbox{supp\,}u_0=\Omega_0\subset\subset\Omega$: the same proof as before shows that we can approximate $u_0$ in $D_p(\Omega)$ with $C_c^\infty(\Omega\setminus\{0\})$ functions.\\
On the other hand the function  $u_\delta$ satisfies $u_\delta\in D_p(\Omega_\delta)$ since $u=0$ on $\partial\Omega$ and $\mbox{supp}\, \eta_\delta\subseteq K_{\delta}$. Since no singularity appears in $\Omega_\delta$, the approximation problem is a classical one: Proposition \ref{Sobolev approximation 1,infty} then proves that $u_\delta$ can be approximated in $D_p(\Omega)$ with functions in $ C_{c,0}^2\left(\Omega\right)$.
\qed

The previous Lemma shows that $C_{c,0}^2\left(\Omega\right)$
 is a core for $A_p$. When $\Omega=\R^N$ or $\Omega= B$, Proposition \ref{core Palla R^N} states that $C_{c,0}^2\left(\Omega\right)$ is also a core for the operator $A_{p,J}$ of Definition \ref{Def A_pJ}. We have therefore proved the following result which provides a description of the operators introduced in the previous subsection.
\begin{prop} \label{domainap}
Let $1\leq p\leq\infty$ and  $\Omega=\R^N$ or $\Omega= B$. Then  the operator $A_p$ coincides with that  of Definition \ref{Def A_pJ} for $J=\N_0$.
\end{prop}

%

In the next lemma we state some interpolative and  a-priori  estimates.
\begin{lem} \label{apriori}
Let $1\leq p\leq \infty$. Then there  exist $\eps_0, \ C>0$ depending only on $c,N,\Om$  such that for
every $0<\eps<\eps_0$ and  $u\in D_p(\Om)$ one has
\begin{align} \label{apriori1}
 \||x|\nabla u\|_p&\leq \eps\|A u\|_p+\frac{C}{\eps}\|u\|_p.
\end{align}
Moreover, if $1<p<\infty$,
\begin{align}\label{apriori2}
 \||x|^2D^2 u\|_p&\leq C(\|A u\|_p+\|u\|_p).
\end{align}
\end{lem}
{\sc Proof.}  In view of Lemma \ref{coreOmega}, it is enough to prove these estimates 
for
$u\in C_{c,0}^2\left(\Omega\right)$. The proof of  (\ref{apriori1})  follows as in \cite[Lemma 2.4]{for-lor} with minor modifications (in particular, one intersects the balls $B(x_0, \rho)$ with $\Omega$). To prove (\ref{apriori2}) for $1<p<\infty$, it is sufficient to apply the classical elliptic estimate $\|D^2u\|_p \le C\|\Delta u\|_p$ (which holds both in $\R^N$ as well as in a bounded $\Omega$ if $u$ vanishes at the boundary)  to $|x|^2u$ and then to interpolate the terms containing $\nabla u$, by (\ref{apriori1}).
\qed

\bigskip
\noindent  
In the next Propositions we  prove  dissipativity properties for $A_p$ through Hardy type inequalities. In the spirit of Section \ref{Section Rellich }, this is equivalent to the fact that the Rellich inequalities \eqref{Intr 1} for the operator $L$,  when $b$ is sufficiently large,  can be proved  using  integration by parts and Hardy  inequalities \eqref{WHardy2-RN}. We begin by the recalling the following  result.
\begin{prop}\label{hardy2} (see \cite[Proposition 8.3]{rellich}).
Let $1<p<\infty$, $\beta\in \R$. Then, if $N-2+\beta\neq 0$, 
for every $u\in C_c^\infty(\R^N\setminus\{0\})$, 
\begin{equation}\label{WHardy2-RN}
\int_{\R^N}
|x|^{\beta}|\nabla u|^{2}|u|^{p-2}
\,dx
\geq 
\left(\frac{N-2+\beta}{p}\right)^2
\int_{\R^N}
|x|^{\beta-2}|u|^{p}
\,dx;
\end{equation}
\end{prop}

We prove now that $A_p$ is  quasi-dissipative.
\begin{teo} \label{dissipativity}
Let $1\leq p\leq \infty$ and set $\omega_p=\frac{N}{p^2}\left[p(N-2+c)-N\right]$.  Then, for every $u\in D_p(\Omega)$, $\lambda>0$,
\begin{equation} \label{contractivity}
\lambda \|u\|_p \le \|(\lambda-A-\omega_p)u\|_p.
\end{equation}
\end{teo}
{\sc Proof.}  We consider, preliminarily, $1<p<\infty$ and  prove the inequality 
\begin{equation} \label{ausiliaria}
-\int_{\Omega}Au |u|^{p-2}u\;
dx\geq \omega_p
 \int_{\Omega}|u|^p\;
dx.
\end{equation}
Let $2\leq p<\infty$. 
By Proposition \ref{coreOmega},  we may assume that $u\in C_{c,0}^2\left(\Omega\right)$.
Setting
$u^\star=u|u|^{p-2}$ we multiply  $Au$  by
$u^\star$ and integrate over $\Omega$. Integrating by parts we get
\begin{align*}
-\int_{\Omega}Au\, u^\star\;
dx&=(p-1)\int_{\Omega}|x|^2|u|^{p-2}|\nabla
u|^2\;dx- (c-2)\int_{\Omega}x\cdot\nabla u\, u|u|^{p-2}\; dx\\[1ex]
&=(p-1)\int_{\Omega}|x|^2|u|^{p-2}|\nabla
u|^2\;dx- \left(\frac{c-2}{p}\right)\int_{\Omega}x\cdot\nabla|u|^p\; dx\\[1ex]
&=(p-1)\int_{\Omega}|x|^2|u|^{p-2}|\nabla
u|^2\;dx+ N\left(\frac{c-2}{p}\right)\int_{\Omega}|u|^p\; dx.
\end{align*}
By Hardy inequality (\ref{WHardy2-RN}) with $\beta=2$,
\begin{align*}
-\int_{\Omega}&Au\, u^\star\;
dx\geq \left[(p-1)\frac{N^2}{p^2}+ N\left(\frac{c-2}{p}\right)\right]\int_{\Omega}|u|^p\; dx=\omega_p \int_{\Omega}|u|^p\; dx
\end{align*}
and therefore
$$-\int_{\Omega}Au |u|^{p-2}u 
dx\geq \omega_p \int_{\Omega}|u|^p\;
dx.$$

For $1<p<2$ the integration by parts is not straightforward (but still allowed, see
\cite{met-spi}) since $|u|^{p-2}$ becomes singular near the zeros of $u$. In this case  it is sufficient to replace $u^\star$ by $u(u^2+\delta)^{\frac{p}{2}-1}$  where $\delta$ is a positive parameter and then let $\delta$ to $0$ obtaining  the required estimates also in this case. 

It is clear that (\ref{ausiliaria}) implies (\ref{contractivity}) which is therefore proved for $1<p<\infty$. Letting $p \to 1, \infty$, we see that (\ref{contractivity}) holds in all cases.
\qed

\begin{os} {\rm
\begin{itemize}
 \item[(i)] $\omega_\infty=0$ and $\omega_1=(c-2)N$;
\item[(ii)] $\omega_p\geq 0$ iff $p\geq \frac{N}{N-2+c}$. Moreover $\omega_p$ attains its maximum value at $\ov{p}=\frac{2N}{N-2+c}$ and $\omega_{\ov{p}}=\left(\frac{N-2+c}{2}\right)^2$.
\end{itemize}}
 \end{os}
The previous theorem, combined with Lemma \ref{apriori}, allows us to deduce the following result. 
\begin{cor} \label{aprioriLambda}
Let $1\leq p\leq\infty$. There exist two constants $\Lambda>0$ and $C>0$ such that, for
every $u\in D_p(\Omega)$ and every $\rm{Re}\lambda\geq \Lambda_p$
$$|\lambda|\|u\|_p+|\lambda|^\frac{1}{2}\||x|\nabla u\|_p\leq C\|\lambda u-Au\|_p.$$
If $1<p<\infty$, we have also 
$$
\||x|^2D^2 u\|_p\leq C\|\lambda u-Au\|_p.
$$
\end{cor}
{\sc Proof.} The estimate
$$|\lambda|\|u\|_p\leq C\|\lambda u-Au\|_p$$ is nothing but
sectoriality. The gradient estimate follows from it, using
(\ref{apriori1}) with $\eps=|\lambda|^{-\frac{1}{2}}$. 
The Hessian estimate for $1<p<\infty$ follows from (\ref{apriori2}).\\
\qed

The next theorem shows that $A_p$ is the generator of a contractive analytic semigroup in $L^p(\Omega)$.

\begin{teo} \label{gen-prel}
For any $1\leq p\leq\infty$, the operator $(A_p+\omega_p, D_p(\Omega))$ generates a contractive analytic semigroup in $L^p(\Omega)$.
\end{teo}
{\sc{Proof.}} To distinguish, we write $\tilde{A}_p$ for $A_p$ when $\Omega=\R^N$. 
Observe that, by Proposition \ref{analyt},  $\tilde{A}_p$  generates an analytic semigroup, hence  its resolvent contains a sector 
$$\Sigma_{\theta,\rho}=\{\lambda \in \C: |\lambda| \ge \rho, |{\rm Arg }\lambda |<\theta \},$$ with $\theta >\pi/2$ where the following resolvent estimate holds
$$\|(\lambda- \tilde{A}_p)^{-1}\|_p \le \frac{M}{|\lambda|}.$$
Let $\Omega\subset\R^N$ and define $\eta_\delta$ and $\eta_0$ as in  (\ref{Partiton unity eq}). 
For $\lambda \in \Sigma_{\theta, \rho}$,
$f\in
L^p(\Omega)$, set $R_0(\lambda)f=\eta_0(\lambda-\tilde{A}_p)^{-1}(\eta_0 f)\in D_p(\R^N)$,  $R_\delta(\lambda)f=\eta_\delta(\lambda-A_\delta)^{-1}(\eta_\delta f)\in W^{2,p}(K_\delta)\cap W_0^{1,p}(K_\delta)$ where $A_\delta$ is the operator $A$ in $K_\delta$ with Dirichlet boundary conditions.
We have 
\begin{align*}
(\lambda-A)R_0(\lambda)f&=(\lambda-A)\eta_0(\lambda-\tilde{A}_p)^{-1}(\eta_0 f)\\&= 
\eta_0(\lambda-A)(\lambda-\tilde{A}_p)^{-1}(\eta_0f)+[\eta_0,
A](\lambda-\tilde{A}_p)^{-1}(\eta_0 f)\\&= \eta_0^2 f +[\eta_0,
A](\lambda-\tilde{A}_p)^{-1}(\eta_0f):=\eta_0^2f+S_0(\lambda)f
\end{align*}
where
$$
[\eta_0, A]g=\eta_0(Ag)- A(\eta_0 g)
$$
is a first order operator supported on $K_\delta$.  Using Corollary \ref{aprioriLambda} (and disregarding $|x|$ which is bounded above and below from 0 in $K_\delta$) we see that
$$\|S_0(\lambda)f\|_p \leq c_1\frac{\|f\|_p}{|\lambda|^\frac{1}{2}}$$
for $\lambda \in \Sigma_{\theta, \rho}$ and with $c_1$ depending
only on $\delta$. In similar way we get 
\[
(\lambda-A)R_\delta(\lambda)f=\eta_\delta^2f+S_\delta(\lambda)f
\]
with 
$$\|S_\delta(\lambda)f\|_p \leq c_1\frac{\|f\|_p}{|\lambda|^\frac{1}{2}}$$
for $\lambda \in \Sigma_{\theta, \rho}$ and with $c_1$ depending
only on $\delta$, by classical results, since $A_\delta$ is uniformly elliptic in $K_\delta$.
 Then setting
$$
R(\lambda):=R_0(\lambda)+R_\delta(\lambda), \quad
S(\lambda):=S_0(\lambda)+S_\delta(\lambda),
$$
we have
\[
(\lambda-A)R(\lambda)f=f+S(\lambda)f.
\]
Choosing  $|\lambda| >\rho_1$ large enough, 
we find $\|S(\lambda)\|_p\leq\frac{1}{2}$ and 
then the operator $I+S(\lambda)$ is invertible in $L^p(\Omega)$. 
Setting $V(\lambda)=(I+S(\lambda))^{-1}$
we have
\[
(\lambda-A)R(\lambda)V(\lambda)f=f
\]
and hence the operator $R(\lambda)V(\lambda)$, which maps $L^p(\Omega)$ into $D_p(\Omega)$, 
is a right inverse of $\lambda-A$. Since both $\|R_0(\lambda)\|_p,\| R_\delta (\lambda)\|_p \le M|\lambda|^{-1}$ and $\|V(\lambda)\|_p  \le 2$, then 
\begin{equation} \label{StimaRes}
\|R(\lambda)V(\lambda)\|_p\leq \frac{C}{|\lambda|}
\end{equation}
for $\lambda \in \Sigma_{\theta, \rho_1}$. 
Clearly, 
$R(\lambda)V(\lambda)$ coincides with $(\lambda-A_p)^{-1}$ whenever
this last is injective, in particular  for $\lambda>-\omega_p$, By Proposition \ref{dissipativity}.
Then $(-\omega_p, \infty)\subset \rho (A_p)$, the  a-priori estimates (\ref{StimaRes}) shows that the norm of the resolvent cannot blow up in   $\Sigma_{\theta, \rho_1}$, hence $\Sigma_{\theta, \rho_1} \subset \rho(A_p)$ and the proof is complete.

\qed

 In the next proposition we prove that the domain $D_p(\Omega)$ coincides
with the maximal one. In what follows, $Au$ is understood in the sense of distributions in $\Omega \setminus \{0\}$. Since the coefficients of $A$ are $C^\infty$ away from the origin,  by local elliptic regularity it follows that  $u \in W^{2,p}_{loc}(\R^N \setminus \{0\})$ when $\Omega=\R^N$ and that $u \in W^{2,p}(\Omega \setminus B_\eps)$ for every $\eps >0$, when $\Omega$ is bounded. This clearly holds for $1<p<\infty$; when $p=\infty$, the same is true for any $q <\infty$.
\begin{prop} \label{L1}
Let $1\leq p\leq \infty$. The domain $D_p(\Omega)$ defined in \eqref{dp} coincides
 with the maximal domain
\begin{align}\label{maximal domain}
 D_{p,max}(\Omega)=\{u\in L^p(\Omega): \,u=0 \text{ on } \partial \Omega,\  Au\in L^p(\Omega)\}.
\end{align}

\end{prop}
{\sc Proof.} The inclusion  $D_p(\Omega) \subset D_{p,max}(\Omega)$ is obvious.
Conversely, let  $u\in D_{p,max}(\Omega)$ and $\lambda>0$ be in the
resolvent set of $(A_p, D_p(\Omega))$. Set $f=\lambda u-A_p u$ and
$v=u-R(\lambda, A_p)f$. Then $v$ belongs to $D_{p,max}(\Omega)$ and
satisfies $\lambda v-A_p v=0$. We prove that $v\equiv 0$ if
$\lambda$ is large enough.
Let us consider for large $n$ 
$$\Omega_n=\left\{x\in\Omega:\ |x|\geq \frac{1}{n},\ \mbox{dist}(x,\partial{\Omega})\geq\frac{1}{n}\right\},\quad \xi_n=\chi_{\Omega_{\frac n 2}}\ast\phi_\frac{1}{n}$$ where $\phi$ is a classical mollifier supported in $B_1$, with $\int_{\R^N}\phi=1$
and $\phi_{\frac{1}{n}}(x)=n^N\phi\left(nx\right)$. It is easy to check that $\xi_n(x)=1$ for $x\in\Omega_{\frac{n}{3}}$, $\xi_n$ is supported in $\Omega_n$ and that $0 \le \xi_n \le 1$, $|\nabla \xi_n|\leq Cn$, $|D^2\xi_n|\leq Cn^2$. Consider also a smooth function $\eta$ such that $\chi_{B_1}\leq\eta\leq \chi_{B_2}$ and set $\eta_n(x)=\eta\left(\frac{x}{n}\right)$, 
$\zeta_n=\xi_n \eta_n$. Since $|\nabla \xi_n| \le Cn\chi_{(\Omega_{n} \setminus \Omega_\frac{n}{3})}$ and  $|\nabla \eta_n| \le Cn^{-1} \chi_{(B_{2n}\setminus B_n)}$,  it follows that the function $\nabla \zeta_n$ has support in $F_n:=\left (\Omega_{n} \setminus \Omega_\frac{n}{3}\right) \cup \left (B_{2n} \setminus B_n\right)$ and satisfies $|x|^2|\nabla \zeta_n|^2 \le C$, with $C$ independent of $n$.\\
Let us consider, first,  the case where $p\geq 2$. Integrating by parts the identity
$$
\int_{\Omega}(\lambda v-Av)v|v|^{p-2}\zeta_n^2=0
$$  we obtain
\begin{align*}
 0 =&\lambda \int_{\Omega}|v|^p\zeta_n^2\,
dx+(p-1)\int_{\Omega}|x|^2|\nabla v|^2|v|^{p-2}\zeta_n^2\, dx
\\[1ex]
&+2\int_{\Omega}|x|^2\zeta_n|v|^{p-2}v\nabla v\cdot\nabla\zeta_n\,
dx+(2-c)\int_{\Omega}\zeta_n^2|v|^{p-2}v\,x\cdot\nabla v\, dx.
\end{align*}
Using  H\"{o}lder's inequality  we obtain
\begin{align*}
 \left|\int_{\Omega}|x|^2\zeta_n|v|^{p-2}v\nabla v\cdot\nabla\zeta_n\, dx\right|&\leq \left(\int_{\Omega}|x|^2\zeta_n^2|\nabla v|^2|v|^{p-2}\, dx\right)^\frac{1}{2}\left(\int_{\Omega}|x|^2|v|^p|\nabla\zeta_n|^2\, dx\right)^\frac{1}{2}\\[1ex]
 &\leq C\left(\int_{\Omega}|x|^2\zeta_n^2|\nabla v|^2|v|^{p-2}\, dx\right)^\frac{1}{2}\left(\int_{\Omega\cap F_n}|v|^p\, dx\right)^\frac{1}{2}\\[1ex]
 &\leq
\eps \int_{\Omega}|x|^2\zeta_n^2|\nabla v|^2|v|^{p-2}\,
dx+\frac{C}{\eps}\int_{\Omega\cap F_n} |v|^p\,
dx.
\end{align*}
 Similarly
$$\left|\int_{\Omega}\zeta_n^2|v|^{p-2}v x\cdot\nabla v\, dx\right|\leq \eps\int_{\Omega}|x|^2\zeta_n^2|\nabla v|^2|v|^{p-2}\, dx+\frac{C}{\eps}\int_{\Omega} |v|^p\zeta_n^2\, dx.$$
Combining the
last inequalities we obtain, up to slightly changing the constants, 
%
\begin{equation*}
\left(\lambda-\frac{C_1}{\eps}\right)\int_{\Omega}|v|^p\zeta_n^2\,
dx+(p-1-3\eps)\int_{\Omega}|x|^2|\nabla v|^2|v|^{p-2}\zeta_n^2\
dx-\frac{2C_1}{\eps}\int_{\Omega\cap F_n} |v|^p\,
dx\leq 0.
\end{equation*}
Finally, choosing $3\eps <p-1$ and letting $n$ to infinity, we obtain
\begin{equation*}
\left(\lambda-\frac{C_2}{p-1}\right)\int_{\Omega}|v|^p\, dx\leq 0
\end{equation*}
which implies $v\equiv 0$, if $\lambda$ is large enough.
For $1<p<2$ the integration by parts is not straightforward  since $|v|^{p-2}$ becomes singular near the zeros of $v$, but still allowed ( see
\cite{met-spi}) and one concludes as before (or, more simply,  notice that $v$ is a smooth function, by elliptic regularity, replace $v|v|^{p-2}$ by $v(v^2+\delta)^{\frac{p}{2}-1}$ and then let $\delta \to 0$).

For $p=1$, we notice that $v$ is a smooth function away from the origin, by elliptic regularity, and  consider a sequence of smooth functions $h_n:\R\rightarrow \R$ such
that $|h_n|\leq 1,$ $h_n'(s)\ge 0$ and  $h_n(s)\rightarrow
\sign(s)$ for every $s\in\R$.  Integrating by parts the identity
$$
\int_{\Omega}(\lambda v-Av)h_k(v)\zeta_n^2=0
$$ the proof follows  as before.

For $p=\infty$ we note that $v$ vanishes at $0$ and at $\partial \Omega$ when $\Omega $ is bounded or at $\infty$ if $\Omega=\R^N$. Moreover, by elliptic regularity, $v$ is a smooth  function out of the origin. If $v$ is not identically zero, then it has a positive maximum point (or a negative minimum point ) at some $x_0 \in \Omega$. The classical maximum principle yields $Av(x_0) \le 0$, hence $\lambda v(x_0) \le 0$, which is a contradiction for $\lambda >0$.
\qed

Finally,  we consider the domain of the operator $A_{p,J}$ of Subsection 3.3.

\begin{cor} \label{domainpj}
If  $\Omega=\R^N$ or $\Omega=B$, then the domain $D_{p,J}(\Omega)$ of $A_{p,J}$ is given by
$$
D_{p,J}(\Omega)=D_p(\Omega)\cap L^p_J(\Omega)=D_{p,max}(\Omega) \cap L^p_J(\Omega).
$$
\end{cor}
{\sc Proof. } By Corollary \ref{restriction}, the domain of $A_{p,J}$ is the intersection of the domain of $A_p$ with $L^p_J$ and the thesis follows from Propositions \ref{domainap}, \ref{L1}.
\qed

\section{ Rellich inequalities in $\R^N$ and in $B$ }\label{Section Rellich }
In this section we prove weighted Rellich inequalities for the operator
$$L=\Delta +c\frac{x}{|x|^2}\cdot\nabla -\frac{b}{|x|^2} ,\quad  c,\ b\in\R$$ on $L^p(\Omega)$ when $\Omega=\R^N$ and $\Omega=B$.
 For $1\leq p\leq\infty$,  $\alpha\in\R$ and $J \subset \N_0$ we define
\begin{align*}
D_{p,\alpha, J}(\Omega):&=\left\{u:\ |x|^{\alpha-2}u,\ |x|^{\alpha}Lu\in L^p_J\left(\Omega\right),\ u=0 \text{ on } \partial \Omega\right\}
.
\end{align*}
When $J=\N_0$ we write $D_{p,\alpha}(\Omega)$ in place of $D_{p,\alpha, \N_0}(\Omega)$.
As in the previous section $Lu$ is understood as a distribution in $\Omega \setminus \{0\}$. Since the coefficients of $L$ are $C^\infty$ away from  the origin,  by local elliptic regularity it follows that, if $u \in D_{p, \alpha}(\Omega)$, then   $u \in W^{2,p}_{loc}(\R^N \setminus \{0\})$ when $\Omega=\R^N$ and  $u \in W^{2,p}(\Omega \setminus B_\eps)$ for every $\eps >0$, when $\Omega$ is bounded. This clearly holds for $1<p<\infty$; when $p=\infty$, the same is true for any $q <\infty$.

Defining 
\begin{align*}
 \quad \Phi u=v,\,v(x)=|x|^{\alpha-2}u(x),
 \end{align*}
  we have seen in Section \ref{Preliminaries}  that 
\begin{align*}
|x|^\alpha L u=Av-\mu v,\qquad \mu=b-(2-\alpha)(N-\alpha+c)
\end{align*}
where $A$ is the operator of  Section \ref{Section A} with $c+4-2\alpha$ in place of $c$,
\begin{equation*} 
A=|x|^2\Delta +(c+4-2\alpha)x\cdot\nabla.
\end{equation*}
By construction $\Phi \left (D_{p,\alpha,J}(\Omega)\right )$ coincides with the domain  $D_{p,J}(\Omega)=D_{p,max}(\Omega)\cap L^p_J(\Omega)$, see Corollary \ref{domainpj}.
In particular  Rellich inequalities
\begin{equation*} 
\||x|^\alpha Lu\|_p \geq C\||x|^{\alpha-2} u\|_p,\quad u\in D_{p,\alpha,J}(\Omega)
\end{equation*}
are equivalent to the spectral estimates
\begin{equation*}
\|\mu v-Av\|_p \geq C\|v\|_p, \quad v\in D_{p,J}(\Omega)
\end{equation*} 
which, recalling Proposition \ref{Rellich-spectrum}, hold precisely when $\mu\notin A\sigma(A_{p,J})$.
The results of this section are then immediate consequences of Theorem \ref{Spectrum main} and Corollary  \ref{a-priori}.

Let us define 
\begin{equation*} 
\gamma_p(\alpha,c):=\Bigl(\frac{N}{p}-2+\alpha\Bigr)\Bigl(\frac{N}{p'}-\alpha+c\Bigr)=\Bigl(\frac{N-2+c}{2}\Bigr)^2-\Bigl(N\Bigl (\frac12-\frac1p\Bigr)+1+\frac{c}{2}-\alpha\Bigr)^2.
\end{equation*}
and   
\begin{equation*} 
D:=
b+\left(\frac{N-2+c}{2}\right)^2.
\end{equation*}
In what follows we refer to $D$ as the discriminant of  $L$; in \cite{met-soba-spi3, met-negro-spina 1} the authors show  that $D$ takes a fundamental role in generation   properties of $L$.
We recall that ${\cal Q}_p$, ${\cal P}_p$, $\omega_p$ have  been defined  in \eqref{spettrogamma}, \eqref{spettrogamma1}, \eqref{omegap}. For clarity sake, we rewrite them in the present situation where $c+4-2\alpha$ takes the place of  $c$:
 
\begin{align*}
{\cal Q}_p&=\left\{\lambda\in \C\ \textrm{such that}\  {\rm Re}\lambda\leq -\frac{({\rm Im} \lambda)^2}{\left (N\left(1-\frac{2}{p}\right)+2-2\alpha+c\right)^2}-\omega_p\right\},\\[1ex]
{\cal P}_p&=\left\{\lambda=-\xi^2+i\xi\left (N\left(1-\frac{2}{p}\right)+2-2\alpha+c\right)-\omega_p,\,  \xi\in\R\right\},\\[1ex]
\omega_p&=\frac{N}{p^2}\left[p(N+2-2\alpha+c)-N\right].
\end{align*}
 Note that, when $N\Bigl (\frac12-\frac1p\Bigr)+1-\alpha+\frac{c}{2}=0$, then $${\cal Q}_p=]-\infty,-\omega_p].$$

 In the following  lemma we denote with $\sqrt z$ a complex square root of $z$  having  non negative real part.

\begin{lem}\label{Parameters}
Let  $1\leq p\leq \infty$, $j \in \N_0$ and 
$\mu:=b-(2-\alpha)(N-\alpha+c)$. Then  the following properties are equivalent
\begin{itemize}
\item[(i)] $\mu\notin {\cal Q}_p-\lambda(P_{j})$;\medskip
\item[(ii)] $b+\gamma_p(\alpha,c)+\lambda(P_{j})>0$;\medskip
\item[(iii)] $\left|N\Bigl (\frac12-\frac1p\Bigr)+1+\frac{c}{2}-\alpha\right|< \sqrt{D+\lambda(P_{j})}$ and  $D+\lambda(P_{j})>0$;\medskip
\item[(iv)] $\left|N\Bigl (\frac12-\frac1p\Bigr)+1+\frac{c}{2}-\alpha\right|< \Rp\sqrt{D+\lambda(P_{j})}$.
\end{itemize}
\end{lem}
{\sc{Proof.}} The proof follows from elementary calculations after noticing that  
\begin{align*}
\omega_p&=b+\gamma_p(\alpha,c)-\mu,\\[1ex]
\gamma_p(\alpha,c)&=D-b-\Bigl(N\Bigl (\frac12-\frac1p\Bigr)+1+\frac{c}{2}-\alpha\Bigr)^2.
\end{align*}
Since $\mu\in\R$, the conditions $\mu\notin{\cal{P}}_p-\lambda(P_{j})$, $\mu\notin{\cal{Q}}_p-\lambda(P_{j})$ become  $b+\gamma_p(\alpha,c)+\lambda(P_{j})\neq 0$, $b+\gamma_p(\alpha,c)+\lambda(P_{j})> 0$, respectively.\\
\qed

The following is the main result of this section. Part 1  has been already proved in \cite{rellich}.

\begin{teo} \label{RellichFJ}
Let
 $1\leq p\leq \infty$, $\alpha,\ b,\ c\in\R$ and $J\subseteq N_0$ with $j_0:=\min\{j\in J\}$.
\begin{itemize}
\item[1.]
If $\Omega=\R^N$,  Rellich inequalities 
\begin{align*}
\||x|^\alpha Lu\|_p \geq C\||x|^{\alpha-2} u\|_p,\quad u\in D_{p,\alpha,J}(\R^N)
 \end{align*}
hold if and only if 
\begin{align*}
\alpha\neq N\Bigl (\frac12-\frac1p\Bigr)+1+\frac{c}{2}\pm\Rp \sqrt {D+\lambda(P_j)}, \quad \forall\, j\in J,
\end{align*}
or equivalently when \;$b +\gamma_p(\alpha,c)+\lambda(P_j)\neq 0$ for every $j\in J$.
\item[2.] If $\Omega=B$,  Rellich inequalities 
\begin{align*}
\||x|^\alpha Lu\|_p \geq C\||x|^{\alpha-2} u\|_p,\quad u\in D_{p,\alpha,J}(B)
 \end{align*} hold if and only if
\begin{align*}
\alpha&<N\Bigl (\frac12-\frac1p\Bigr)+1+\frac{c}{2}+ \Rp\sqrt{D+\lambda(P_{j_0}}), \quad\text{and}\;\\[1ex]
\alpha&\neq N\Bigl (\frac12-\frac1p\Bigr)+1+\frac{c}{2}-\Rp \sqrt {D+\lambda(P_j)}, \quad\forall\, j\in J.
\end{align*}
In particular the latter conditions are verified
\begin{itemize}
\item[(i)] when  $\alpha\geq N\left(\frac 1 2-\frac{1}{p}\right)+1+\frac c 2$, if and only if \;$b +\gamma_p(\alpha,c)+\lambda(P_{j_0})>0$,
\item[(ii)]  when  $\alpha< N\left(\frac 1 2-\frac{1}{p}\right)+1+\frac c 2$, if and only if \;$b +\gamma_p(\alpha,c)+\lambda(P_j)\neq 0$ for every $j\in J$.
\end{itemize}
\end{itemize}\medskip
If $J=\N_0$  and  $b +\gamma_p(\alpha,c)>0$, that is
\begin{align*}
 \left|N\Bigl (\frac12-\frac1p\Bigr)+1+\frac{c}{2}-\alpha\right|< \Rp\sqrt{D},
 \end{align*}
  then the optimal constant is given by $C= b +\gamma_p(\alpha,c)$.
\end{teo}
{\sc Proof.} Consider  ${\cal Q}_p$, ${\cal P}_p$ and $\omega_p$ defined before Lemma \ref{Parameters} and let $\mu=b-(2-\alpha)(N-\alpha+c)$. Then Rellich inequalities hold if and only if $\mu\notin A\sigma(A_{p,J})$.
The proof of the required claims  follows then easily by combining  Lemma  \ref{Parameters}, Theorem \ref{Spectrum main} and Corollary  \ref{a-priori}.\\\qed

\begin{os}
For a fixed $\alpha$, Rellich inequalities are always true in $L^p_{\geq n}(\Omega)$, for a sufficiently large $n\in\N_0$,  even though they fail in the whole $L^p(\Omega)$. This phenomenon appears also in the extreme cases $p=1,\infty$.  The failure of Rellich  inequalities  for some values of $\alpha$  is, therefore, always determined by subspaces defined by spherical harmonics of  low order.  
 \end{os}
 
\medskip 
When $b=c=0$, the operator reduces to the  Laplace operator $L =\Delta$. In this case 
\begin{align*}
D=\left(\frac{N-2}{2}\right)^2, \quad D+\lambda_n=\left(\frac{N-2}{2}+n\right)^2.
\end{align*}
 Rellich inequalities in bounded domains for the Laplace operator have already been investigated in \cite{musina} where  their validity is proved  for $N\geq 3$, $1<p<\infty$ and   
\begin{align}\label{range easy}
 -\frac N p+2<\alpha<N\left(1-\frac 1 p\right).
 \end{align}
 
This range coincides with the values of $\alpha$ for which Rellich inequalities can be proved  using  integration by parts and the Hardy  inequalities \eqref{WHardy2-RN} (see Theorem \ref{dissipativity}). The following corollary characterizes their validity in the ball.

\begin{cor} \label{Rellich-Delta}
Let
 $1\leq p\leq \infty$, $\alpha\in\R$. If $\Omega=B$,   Rellich inequalities 
\begin{align*}
\||x|^\alpha \Delta u\|_p \geq C\||x|^{\alpha-2} u\|_p,\quad u\in D_{p,\alpha}(B)
 \end{align*} hold if and only if
\begin{align*}
\alpha<N\left(1-\frac 1 p\right),\quad
\alpha\neq-\frac N p+2-n, \quad\forall\, n\in\N_0.
\end{align*}
\medskip
\end{cor}

\section{Rellich inequalities in general domains}\label{Rellich general}\label{Rellich Bounded domain}
Let $\Omega$ be  an open  bounded and connected subset of $\R^N$ whose boundary $\partial \Omega$ is $C^{2,\beta}$ and such that $0 \in\Omega$. In this section we show that Rellich inequalities for the operator $L$ hold in  $\Omega$ if and only if they hold in the ball $B$. In terms of the auxiliary operator $A$, this means that its approximate point spectrum is independent of the bounded set $\Omega$. We have no direct proof of this fact which does not seem to be evident.
We write $L$ in the symmetric form
\begin{equation}\label{symmetrize}
L=\Delta +c\frac{x\,}{|x|^2}\cdot\nabla -\frac{b\,}{|x|^{2}}=
|x|^{-c}{\rm div}(|x|^{c}\nabla ) - \frac{b}{|x|^2}. 
\end{equation}
and we always assume  $1<p<\infty$ and that 
 \begin{equation} \label{D}
D=b+\left (\frac{N-2+c}{2}\right )^2 \ge 0.
\end{equation}
This condition is crucial for the solvability of some elliptic problems related to $L$ which will be studied in the following subsection in a auxiliary  weighted $L^2$ space.

\subsection{The operator $L$ in $L^2(\Omega, d\mu)$}\label{Section estimate}

We need some preliminary facts concerning the operator $L$ in a weighted space and here we suppose $\Omega$ as above  or $\Omega=\R^N$.  
We consider the weighted space $L^2(\Omega,d\mu)$,  $d\mu=|x|^{c}dx$, and the symmetric form 
\begin{align*}
\mathfrak{a}(u,v)
:=
\int_{\Omega}\left(
 \nabla u\cdot \nabla \overline{v}+ \frac{b}{|x|^2}u\overline{v}\right)\,d\mu, \qquad
u,v \in   C_{c,0}^2(\Omega).
\end{align*}
Using \eqref{symmetrize}, we see that for $u,v \in   C_{c,0}^2(\Omega)$
$$\int_{\Omega} (Lu)\, \overline{v} \, d\mu=\mathfrak{a}(u,v).$$


To prove that $\mathfrak{a}$ is non-negative, we make different change of variables according to $D>0$ or $D=0$.
When $D>0$ we write   $u=u_1|x|^{-\frac{c}{2}}$ and $v=v_1|x|^{-\frac{c}{2}}$ to obtain, after integration by parts
\begin{align}
\label{N-dimensional form change variables} 
\mathfrak{a}(u,v)=
\int_{\Omega}\left(
 \nabla u_1\cdot \nabla \overline{v}_1+\left(D-\frac{(N-2)^2}{4}\right) \frac{u_1\overline{v}_1}{|x|^2}\right)\,dx.
\end{align}
Then we use the classical Hardy inequality. When $D=0$ we are in the critical case of Hardy inequality and it is convenient to use the transformation (which is the basis of the proof of Hardy inequality)
 $u=u_1|x|^{-\frac{N-2+c}{2}}$ and $v=v_1|x|^{-\frac{N-2+c}{2}}$. Integrating by parts we get
\begin{align}
\label{N-dimensional form change variables-crit} 
\mathfrak{a}(u,v)=
\int_{\Omega}\left(
 \nabla u_1\cdot \nabla \overline{v}_1\right)|x|^{2-N}\,dx.
\end{align}
To identify the domain of the closure of $\mathfrak{a}$ we use the classical Sobolev space $H^1_0(\Omega)$ and also $H_0^1\left(\Omega,|x|^{2-N}dx\right)$ defined as the closure of  $C_{c}^2(\Omega)$  with respect to the norm  
$$\left\|v\right\|^2_{H_0^1\left(\Omega,|x|^{2-N}dx\right)}=\int_{\Omega}\left[|\nabla v|^2+ \left|v\right|^2\right]\,|x|^{2-N}dx.$$
Note that we use $C^2_c(\Omega)$ and not $C_{c,0}^2(\Omega)$, that is we do not assume that the functions vanish in a neighbourhood of $0$. However, the above definition would not change  using the smaller space. Let us recall, in fact, that, since $N \ge 2$, $C_{c,0}^2(\Omega)$ is dense in $H^1_0(\Omega)$ and the same is true for $H_0^1\left(\Omega,|x|^{2-N}dx\right)$, as we show below. 

\begin{lem}
$C_{c,0}^2(\Omega)$ is dense in $H_0^1\left(\Omega,|x|^{2-N}dx\right)$.
\end{lem}
{\sc Proof. }
Let us assume, for example that $\Omega=\R^N$ and let $f\in C_c^2(\R^N)$. We approximate $f$ in the norm of  $H_0^1\left(\Omega,|x|^{2-N}dx\right)$ with functions belonging to $C_{c,0}^2 (\R^N)$. \\
Let $\varphi\in C^\infty(\R^+)$ such that $ \varphi(r)=0$ if $0\leq r\leq \frac 1 4$ and  $ \varphi(r)=1$ if $ r\geq \frac 1 2$ and set  $\varphi_\epsilon(x):=\varphi(|x|^\epsilon)$. By construction $f\varphi_\epsilon\in C_{c,0}^2(\R^N)$ and,  as $\epsilon\to 0^+$,  $f\varphi_\epsilon$, $\partial_if\varphi_\epsilon$ converge in $L^2\left(\Omega,|x|^{2-N}dx\right)$ to $f$, $\partial_i f$, respectively,  by dominated convergence.\\ It remains to show that  $f\partial_i\varphi_\epsilon$ converges to $0$ in $L^2\left(\Omega,|x|^{2-N}dx\right)$.
This is true since
\begin{align*}
&\int_{\R^N} |f|^2|\partial_i\varphi_\epsilon|^2\,|x|^{2-N}dx \leq \|f\|^2_\infty \int_{(\frac 1 4)^{\frac{1}{\epsilon}}\leq |x|\leq (\frac 1 2)^{\frac{1}{\epsilon}}}|x|^{2\epsilon-2} \epsilon^2\left|\varphi'(|x|^\epsilon)\right|^2\,|x|^{2-N}dx\\[1ex]
&=\epsilon^2\|f\|^2_\infty|S^{N-1}|\int_{(\frac 1 4)^{\frac{1}{\epsilon}}}^{(\frac 1 2)^{\frac{1}{\epsilon}}}|\varphi'(r^\epsilon)|^2r^{2\epsilon-1}\,dr
=\epsilon\|f\|^2_\infty|S^{N-1}|\int_{\frac 1 4}^{\frac 1 2}|\varphi'(s)|^2s\,ds.
\end{align*}
\qed

To prove the main properties of  $\mathfrak{a}$ we may therefore use $C_{c,0}^2(\Omega)$.

\begin{lem}\label{Close form Nd}
Let $D\geq 0$. The form $\mathfrak{a}$ is  non-negative and symmetric in $L^2(\Omega,d\mu)$. For $u\in C_{c,0}^{2}\left(\Omega\right)$, let  $||u||_{\mathfrak{a}}:=\sqrt{\mathfrak{a}(u,u)+||u||^2_{L_\mu^2}}$. Then  $||u||_{\mathfrak{a}}$ is equivalent to $\||x|^{\frac{c}{2}}\,u\|_{H_0^1\left(\Omega \right)}$,  if $D>0$, and to $\||x|^{\frac{N-2+c}{2}}\,u\|_{H_0^1\left(\Omega,|x|^{2-N}dx\right)}$, if $D=0$. 
\end{lem}
{\sc Proof.} If $D>0$ we set $u=v|x|^{-\frac{c}{2}}$. We choose $\eps$ small enough such that $D-\eps\frac{(N-2)^2}{4}>0$. 
Using \eqref{N-dimensional form change variables} and  Hardy inequality
$$\int_\Omega |\nabla v|^2\, dx  \ge \frac{(N-2)^2}{4} \int_\Omega \frac{|v|^2}{|x|^2}\, dx
$$
 we obtain
\begin{align} \label{injective}
\mathfrak{a}(u,u) \geq \eps\int_{\Omega}
 |\nabla v|^2\,dx+  
 \left(D-\eps\frac{(N-2)^2}{4}\right)\int_{\Omega}\frac{|v|^2}{|x|^2}\,dx\geq \eps\int_{\Omega}
 |\nabla v|^2\,dx.
\end{align}
On the other hand, by Hardy inequality again, 
\begin{align*}
\mathfrak{a}(u,u) \leq C\left(\int_{\Omega}
 |\nabla v|^2\,dx+  
 \int_{\Omega}\frac{|v|^2}{|x|^2}\,dx\right) \leq \tilde{C}\int_{\Omega}
 |\nabla v|^2\,dx.
\end{align*}
This proves that $||u||_{\mathfrak{a}}$ and $||v||_{H_0^1\left(\Omega \right)}$ are equivalent norms.
If $D=0$, setting $u=v|x|^{-\frac{N-2+c}{2}}$, we obtain from (\ref{N-dimensional form change variables-crit}) 
$$
\mathfrak{a}(u,u) =\int_{\Omega} |\nabla v|^2|x|^{2-N}
\, dx.
$$
Since also the norms of $u$ in $L^2(\Omega, d\mu)$ and $v$ in $L^2(\Omega, |x|^{2-N}\, dx)$ coincide, we see
that the norms $||u||_{\mathfrak{a}}$ and $||v||_{H_0^1\left(\Omega,|x|^{2-N}dx\right)}$ are equivalent.\\
\qed


Using the density of  $C_{c,0}^2(\Omega)$  in $H^1_0(\Omega)$ and in  $H_0^1\left(\Omega,|x|^{2-N}dx\right)$, we extend the form ${\mathfrak{a}}$  to the domain
\begin{align*}
D({\mathfrak{a}})&=\left\{u\in L^2(\Omega,d\mu): u|x|^{\frac{c}{2}}\in H_0^1\left(\Omega \right)\right\}, &\text{for}\quad D>0,\\[1.5ex]
D({\mathfrak{a}})&=\left\{u\in L^2(\Omega,d\mu): u|x|^{\frac{N-2+c}{2}}\in H_0^1\left(\Omega, |x|^{2-N}dx \right)\right\}, &\text{for}\quad D=0,
\end{align*}
thus obtaining a closed form. 

Note that both the norms of  $u|x|^{\frac{c}{2}}$ and $u|x|^{\frac{N-2+c}{2}}$ in the corresponding spaces equal the norm of $u$ in $L^2(\Omega, d\mu)$. The transformation $u= u_1|x|^{\frac{N-2+c}{2}}$ can be performed also in the case $D>0$. However it leads to the extra term $D(u_1v_1)/|x|^2$ in the integral \eqref{N-dimensional form change variables-crit} which cannot be dominated by the norm of $H_0^1\left(\Omega, |x|^{2-N}dx \right)$. 

Let $-L$ be the operator associated to  ${\mathfrak{a}}$, that is
\begin{align}\label{Definition operator in L^N}
D(L):=
\left\{u\in D({\mathfrak{a}})
\;;\;
\exists v\in L^2_{\mu}\ \text{s.t.}\ {\mathfrak{a}}(u,w)
=
\int_{\Omega}v\overline{w}\,d\mu
\quad\forall w\in D(\mathfrak{a})
\right\}, 
\quad 
-Lu:=v.
\end{align}
Clearly, $L$ is given by \eqref{symmetrize} when $u \in C^2_{c,0}(\Omega)$. In the next lemma we prove the simplest inequality useful to prove compactness when $D=0$. Note that Hardy inequality fails with respect to the weight $|x|^{2-N}$.

\begin{lem}\label{Compactness D(a)}
Let $\Omega$ be bounded and let $R(\Omega):=\max_{x\in\Omega}|x|$. Then, for every $u\in  C_{c,0}^2(\Omega)$,
\begin{align*}
\int_{\Omega} \frac{|u|^2}{|x|}\, |x|^{2-N}\,dx\leq 4R(\Omega)\int_\Omega |\nabla u|^2\,|x|^{2-N}\,dx.
\end{align*}
In particular the immersion $H_0^1\left(\Omega, |x|^{2-N}dx \right)\hookrightarrow L^2\left(\Omega, |x|^{2-N}dx \right)$ is compact.
\end{lem}
{\sc Proof.} Let us fix $u\in  C_{c,0}^2(\Omega)$. Integrating by parts we have
\begin{align*}
\int_{\Omega} \frac{|u|^2}{|x|}\, |x|^{2-N}\,dx=\int_\Omega |u|^2\mbox{div}(|x|^{1-N}x)\,dx=- 2\int_\Omega u\nabla u\cdot(|x|^{1-N}x)\,dx.
\end{align*}
This implies, using the Cauchy-Schwarz inequality,
\begin{align*}
\int_{\Omega} \frac{|u|^2}{|x|}\, |x|^{2-N}\,dx&\leq 2\int_\Omega |u|\,|\nabla u|\,|x|^{2-N}\,dx\leq 2\sqrt{R(\Omega)}\int_\Omega \frac{|u|}{\sqrt {|x|}}|\nabla u|\,|x|^{2-N}\,dx\\[1ex]
&\leq 2\sqrt{R(\Omega)}\left(\int_\Omega \frac{|u|^2}{|x|}\,|x|^{2-N}\,dx\right)^{\frac 1 2}\left(\int_\Omega|\nabla u|^2\,|x|^{2-N}\,dx\right)^{\frac 1 2}
\end{align*}
and the inequality follows.
To prove the compactness of the embedding, we take $u$  in the unit ball 
 $\mathcal{B}$ of $ H_0^1\left(\Omega, |x|^{2-N}dx \right)$ and fix $\epsilon>0$. Then 
\begin{align*}
\int_{\Omega\cap B_\epsilon}|u|^2\, |x|^{2-N}\,dx\leq \int_{\Omega\cap B_\epsilon}\frac{\epsilon}{|x|}|u|^2\, |x|^{2-N}\,dx \leq 4\epsilon R(\Omega).
\end{align*}
Since $L^2\left(\Omega\cap B^c_\epsilon,\,|x|^{2-N}dx\right)=L^2\left(\Omega\cap B^c_\epsilon,\,dx\right)$, the compactness of $\mathcal{B}_{\vert \Omega\cap B^c_\epsilon}$ in $L^2\left(\Omega\cap B^c_\epsilon,\,|x|^{2-N}dx\right)$ is classical. This fact and and the above estimate show that $\mathcal{B}$ is totally bounded.\\
\qed

In the next Proposition we collect the main properties of $L$ in $L^2(\Omega, d\mu)$..
\begin{prop} \label{LForm}
The operator $-L$ defined in (\ref{Definition operator in L^N}) is non-negative and  self-adjoint in $L^2(\Omega, d\mu)$. The generated semigroup $T_{\Omega}(t)$ is positivity preserving in $L^2(\Omega, d\mu)$.
Moreover,  $C_{c,0}^2(\Omega)\hookrightarrow D(L)$ and for every $u\in C_{c,0}^2(\Omega)$ 
$$Lu=\Delta u
+c\frac{x}{|x|^2}\cdot\nabla u - \frac{b}{|x|^2}u.$$
If $\Omega$ is bounded then $L$ has compact resolvent and is invertible in $L^2(\Omega, d\mu)$.
\end{prop}
{\sc Proof.} Non-negativity and self-adjointness of $-L$ follow by construction. The positivity of $T_{\Omega}(t)$ follows from that of the resolvent which is a consequence of the Beurling-Deny conditions.

Let us suppose, now, $\Omega$  be bounded and let us prove that $D({\mathfrak{a}})$ is compactly embedded in $L^2(\Omega, d\mu)$.  To this aim let  $\mathcal{U}$ be a bounded subset of $D({\mathfrak{a}})$. Assume  $D>0$; then the set $\mathcal{U}'=\{u|x|^\frac{c}{2}:\, u\in \mathcal{U}\}$ is a bounded subset of $H_0^1(\Omega)$, hence
totally bounded in $L^2(\Omega)$, by the compactness of the embedding of $H^1_0(\Omega)$ into $L^2(\Omega)$. It is then immediate to check that 
$\mathcal{U}$ is totally bounded in $L^2(\Omega, d\mu)$, which proves the claim. The case $D=0$ follows similarly from Lemma \ref{Compactness D(a)}.

In both cases $L$ has  compact resolvent;  its spectrum  consists, therefore, of eigenvalues and, being injective  by  \eqref{N-dimensional form change variables-crit}, \eqref{injective},    $L$  is  invertible.
\qed

 Next we need a maximum principle for the solution of an homogeneous problem related to $L$. Note that no singularity appears, since $0 \not \in V$ below. However, comparison is  not obvious since the coefficient $b$ can be negative even though $D\geq 0$. 

\begin{lem}\label{Maximum principle for L}
Let  $V$ be  an open  bounded and connected subset of $\R^N$ whose boundary $\partial V$ is $C^{2,\beta}$ and such that $0\notin V$. For every $\varphi\in C^2(\partial V)$ the problem
  \begin{align*}
\begin{cases}
-Lv=0,\quad & \text{in}\quad V,\\
v=\varphi,\quad & \text{in}\quad \partial V,
\end{cases}
\end{align*}
admits a unique solution $v\in C^2\left(V\right)\cap C\left(\bar V\right)$. Moreover $v$ satisfies $\inf_{\partial V}\varphi\leq v(x)\leq\sup_{\partial V} \varphi$ for every $x\in V$ .
\end{lem}
{\sc Proof.}  The transformation $Su(x)=|x|^{-\frac{N-2+c}{2}}v(x)$ turns  $L$  into 
$$SLS^{-1}= \Delta -(N-2)\frac{x}{|x|^2}\nabla-\frac{D}{|x|^2},$$ 
 which is uniformly elliptic with smooth coefficients and  non-positive potential. Then the proof  follows, immediately, by classical results.\\\qed

In order to prove Rellich inequalities in domains, we  need estimates for the Green function of $-L$ in $\Omega$, that is for the integral kernel expressing $(-L)^{-1}$ with respect to the Lebesgue measure.  We start with the case $D>0$ where can use the results of \cite{met-negro-spina 1} and compare the Green function in $\Omega$ with that in $\R^N$. 

\begin{prop}\label{comparison Green}
Let $D>0$ and let $G(x,y)$, $x,y\in\Omega\times \Omega$ be the Green function of the operator $L$, written with respect to the Lebesgue measure. Then
\begin{equation} \label{G}
 0 \leq G(x,y) \leq C\, G_0(x,y),
\end{equation}
where
if $N>2$
\begin{align}
 \label{eGreen}
|x|^{\frac{c}{2}}|y|^{-\frac{c}{2}}G_0(x,y)=|x-y|^{2-N}\left(1\wedge\frac{|x||y|}{|x-y|^2}\right)^{\sqrt{D}-\frac{N-2}{2}}
\end{align}
and if $N=2$
\begin{align}
 \label{estimates Green N=2}
|x|^{\frac{c}{2}}|y|^{-\frac{c}{2}}G_0(x,y)=
\begin{cases}
\dfrac {\left(|x||y|\right)^{\sqrt{D}}}{|x-y|^{2\sqrt{D}}},\quad&\text{if}\quad \frac{|x-y|^2}{|x||y|}\geq 1;\\[4ex]
 1-\log \left(\dfrac{|x-y|^2}{|x||y|}\right),\quad&\text{if}\quad \frac{|x-y|^2}{|x||y|}\leq 1.
\end{cases}
\end{align}
\end{prop}
{\sc Proof.} Let  $T_{\Omega}(t)$, $T(t)$ be the semigroups generated by $L$ in $L^2(\Omega, d\mu)$ and $L^2(\R^N, d\mu)$, respectively. From \cite[Sections 2.3, 2.6, Proposition 4.23]{ou} it follows that $0 \le T_{\Omega}(t)f \le T(t)f$ whenever $0 \le f \in L^2(\Omega, d\mu)$. Furthermore from \cite[Corollary 4.6]{cal-met-negro-spina } $T(t)$ is an integral operator whose kernel $p(t,x,y)$, expressed with respect to the Lebesgue measure,  satisfies, for every $\epsilon>0$ and some constant $C_\epsilon>0$, 
\begin{align*}
0\leq p(t,x,y)\leq
C_\epsilon t^{-\frac{N}{2}} |x|^{-\frac{c}{2}}|y|^{\frac{c}{2}}&\left [\left (\frac{|x|}{\sqrt t}\wedge 1 \right) \left (\frac{|y|}{\sqrt t}\wedge 1 \right)\right ]^{-\frac{N}{2}+1+\sqrt{D}}\exp\left(- \dfrac{|x-y|^2}{(4+\epsilon)t}\right).
\end{align*}
Using \cite[Theorem 1.5]{Arendt Bukhalov}, it follows that also $T_{\Omega}(t)$ is an integral operator whose kernel $p_\Omega$ satisfies the same estimate above. By \cite[Theorem 7.1]{met-negro-spina 1}, since  $D>0$,  we have 
\begin{equation} \label{Greenspace}
\int_0^\infty p(t,x,y)\,dt\leq  C G_0(x,y)
\end{equation}
hence
\begin{align*}
G(x,y)=\int_0^\infty p_{\Omega}(t,x,y)\,dt\le \int_0^\infty p(t,x,y)\,dt\leq  C G_0(x,y).
\end{align*}
\qed

\begin{os}
The inequality between the semigroups above easily follows from the the corresponding one for the resolvents.  Let $\lambda>0$, $0\leq f\in L^2(\Omega, d\mu)$ and set $u=R(\lambda, L_{\Omega})f$, $w=R(\lambda, L_{\R^N})f$. Then $0\leq u\in D({\mathfrak a}_{\Omega})$ and  $0\leq w\in D({\mathfrak a}_{\R^N})$; furthermore  $\lambda u-Lu=\lambda w -L w$ and, for every $v\in D({\mathfrak a}_{\Omega})$ one has 

\begin{align*}
\lambda\int_{\Omega} (u-w)v\, d\mu=
\int_{\Omega}\left(
 \nabla (w-u)\cdot \nabla v+ \frac{b}{|x|^2}(w-u) v\right)\,d\mu.
\end{align*}
Choosing $v=(u-w)^+\in D( {\mathfrak a}_{\Omega})$  we get
\begin{align*} 
& \lambda\int_{\Omega} \left|(u-w)^+\right|^2\, dx=-{\mathfrak a}_{\Omega}\Big((u-w)^+,(u-w)^+\Big)\leq 0
\end{align*}
which implies $(u-w)^+= 0$ that is $u\leq w$. \qed
\end{os}

\medskip
The case $D=0$ is more involved since, in this case,  the integral in \eqref{Greenspace} is divergent near $\infty$. To overcome this problem, we use the boundedness of $\Omega$ to  improve the decay of $p_\Omega$ as $t\to \infty$. We estimate directly $p_\Omega$ without comparing with the kernel in the whole space, by adapting to our case the arguments of   \cite{cal-met-negro-spina }.

We use the change of variable leading to \eqref{N-dimensional form change variables-crit} to get rid of the potential term $b|x|^{-2}$ 
and introduce  the Hilbert space $L^2(\Omega,|x|^{-2s_1}\,d\mu)=L^2(\Omega,|x|^{2-N}\,dx)$,  where $s_1=\frac{N-2+c}{2}$. Then \eqref{N-dimensional form change variables-crit} reads
\begin{align*}
\mathfrak{b}(u,v)
&:=(\nabla u,\nabla v)_{L^2(\Omega,|x|^{2-N}\,dx)}=a\left(|x|^{-s_1}u, |x|^{-s_1}v\right),\\[1ex]
D(\mathfrak{b})
&:= H_0^1\left(\Omega,|x|^{2-N}\,dx\right).
   \end{align*}
  By construction $\mathfrak{b}$ is the inner product in $ H_0^1\left(\Omega,|x|^{2-N}\,dx\right)$ , and  $u\in L^2(\Omega,|x|^{2-N}\,dx)\mapsto |x|^{-s_1} u\in L^2(\Omega,\,d\mu)$ is an isometry 
which maps $D(\mathfrak{b})$ onto $D(\mathfrak{a})$. The operator $-\tilde L$ associated to $\mathfrak{b}$  then satisfies
\begin{align*} 
D(\tilde{L})=|\cdot|^{s_1}\ D(L),\quad 
\tilde{L}u=|\cdot|^{s_1}L(|\cdot|^{-s_1} u)
\end{align*}
hence
\begin{align}\label{relation kernels}
e^{z\tilde{L}}f=|\cdot|^{s_1}e^{zL}(|\cdot|^{-s_1} f) , \quad f \in L^2(\Omega,|x|^{2-N}\,dx).
\end{align}

Clearly $-\tilde L$ is  non-negative and  self-adjoint in  $L^2\left(\Omega,|x|^{2-N}\,dx\right)$. The semigroup $\left(e^{z\tilde L}\right )_{z \in \C_+}$    is analytic, submarkovian and
satisfies
\begin{align}\label{exp decay}
\|e^{-t\tilde L}\|_{L^2\left(\Omega,|x|^{2-N}\,dx\right)}\leq e^{-\lambda_1 t},
\end{align}
where $\lambda_1>0$ is the first eigenvalue of $- \tilde L$, which is positive since $-\tilde L$ is non-negative and invertible, by the similarity with $-L$.

The following lemma is a special case of 
Caffarelli-Kohn-Nirenberg inequalities and we refer to \cite[Lemma 3.2]{met-soba-spi4} for a short proof. It is used to prove the $L^1$-$L^\infty$ bound of the semigroup.
\begin{lem}\label{CKN}
Let $\sigma\in\R\setminus\{-N\}$. Then for every 
$q\in (2,\infty)$ satisfying $\frac{1}{q}\geq \frac{1}{2}-\frac{1}{N}$,
there exists $C_q>0$ such that 
for every $u\in C^2_{c,0}(\Omega)$,
\begin{align*}
\left(
\int_{\Omega}
  |u(x)|^q
|x|^{\sigma}\,dx
\right)^{\frac{1}{q}}
&\leq 
C_q
\left(
\int_{\Omega}
  |\nabla u(x)|^2
|x|^{(1-\frac{2}{N})\sigma}\,dx
\right)^{\frac{N}{2}\left (\frac{1}{2}-\frac{1}{q} \right )}
\left(
\int_{\Omega}
  |u(x)|^2
|x|^{\sigma}\,dx
\right)^{\frac{1}{2}-\frac{N}{2}\left (\frac{1}{2}-\frac{1}{q} \right )}. 
\end{align*}
In particular, when $\Omega$ is bounded and $\sigma \le 0$, then 
\begin{align*}
\left(
\int_{\Omega}
  |u(x)|^q
|x|^{\sigma}\,dx
\right)^{\frac{1}{q}}
&\leq 
C_{q, \Omega}
\left(
\int_{\Omega}
  |\nabla u(x)|^2
|x|^{\sigma}\,dx
\right)^{\frac{N}{2}\left (\frac{1}{2}-\frac{1}{q} \right )}
\left(
\int_{\Omega}
  |u(x)|^2
|x|^{\sigma}\,dx
\right)^{\frac{1}{2}-\frac{N}{2}\left (\frac{1}{2}-\frac{1}{q} \right )}. 
\end{align*}
\end{lem}

\begin{prop}\label{comparison Green D=0}
Let $D=0$ and  $\Omega$ be bounded. Then the semigroup $T_{\Omega}(t)$ generated by $L$ in $L^2(\Omega, d\mu)$ has an  heat kernel $p(t,x,y)$,  with respect to the Lebesgue measure,  which satisfies, for every $\epsilon>0$ and some constant $C_\epsilon>0$ 
\begin{align}\label{HK est D=0}
 p(t,x,y)\leq C_\epsilon t^{-\frac{N}{2}}e^{-\frac{\lambda_1}{3} t}|x|^{-s_1}|y|^{c-s_1}\exp\left(-\dfrac{|x-y|^2)}{(4+\epsilon)t}\right).
\end{align}
The Green function  $G(x,y)$ of $L$, again written with respect to the Lebesgue measure, satisfies for some constant $C, k>0$, 
\begin{equation} \label{G D=0}
 0 \leq G(x,y) \leq C\, G_0(x,y),
\end{equation}
where
if $N>2$
\begin{align*}
G_0(x,y)=|x|^{-s_1}|y|^{c-s_1}e^{-c|x-y|}\left(1\wedge |x-y|\right)^{2-N}
\end{align*}
and if $N=2$
\begin{align*}
G_0(x,y)=|x|^{-s_1}|y|^{c-s_1}
\begin{cases}
e^{-k|x-y|},\quad&\text{if}\quad |x-y|\geq 1;\\[2ex]
 1-\log \left(|x-y|\right),\quad&\text{if}\quad |x-y|\leq 1.
\end{cases}
\end{align*}
\end{prop}
{\sc{Proof.}} We make use of the results and methods of \cite[Sections 3,4]{cal-met-negro-spina }, pointing out the appropriate changes due to the boundedness of $\Omega$.  The $L^p$ norms used here refer to the measure $|x|^{2-N}\, dx$.

The ultracontractivity estimate for $t \ge 0$
   \begin{align*}
 \|e^{t\tilde L}\|_{1\to\infty}\leq C t^{-\frac N 2}
 \end{align*}
follows from Lemma \ref{CKN} with $\sigma=2-N \le 0$ and any fixed $q$ as in its statement, using \cite[Theorem 6.2]{ou}.

 Since $\tilde L$ is self-adjoint we have also $
 \|e^{it\tilde L}\|_{2\to2}\leq 1$ for $t\in\R$. Using  $\|T^*T\|_{1\to \infty}=\|T\|_{1\to 2}^2$ and recalling \eqref{exp decay}, we obtain  for $t>0, \,s\in\R$,
 \begin{align*}
\|e^{-(t+is)\tilde L}\|_{1\to \infty}&\le \|e^{-\frac t 3 \tilde L}\|_{1 \to 2}\|e^{-\frac t 3 \tilde L}\|_{2 \to 2}
\|e^{-is \tilde L}\|_{2 \to 2}\|e^{-\frac t 3 \tilde L}\|_{2 \to \infty}\\[1ex]
&\leq\|e^{-\frac t 3 \tilde L}\|^2_{1 \to 2}e^{-\frac {\lambda_1}{ 3} t }=\|e^{-t  \tilde L}\|_{1 \to \infty}e^{-\frac {\lambda_1}{ 3} t }\le C {\ t}^{-N/2}e^{-\frac {\lambda_1}{ 3} t }.
 \end{align*}

This  proves 
   \begin{align*}
 \|e^{z\tilde L}\|_{1\to\infty}\leq C \left(\Rp z\right)^{-\frac N 2}e^{-\frac {\lambda_1}{ 3} \Rp z },\quad \forall z\in\C^+.
 \end{align*}
The Dunford-Pettis criterion  yields the existence of a  kernel $\tilde p$ such that, for $z\in\C_+$,
\[
e^{z\tilde{L}}f(x)=
\int_{\Omega}\tilde p(z,x,y)f(y)\,|x|^{2-N}dx,\quad f\in  L^1\left(\Omega,|x|^{2-N}\,dx\right)\cap L^\infty(\Omega)
\]
and 
$$\underset{x,y\in \Omega\setminus\{0\}}{\mbox{ sup}} |\tilde p(z,x,y)| \le C \left(\Rp z\right)^{-\frac N 2}e^{-\frac {\lambda_1}{ 3} \Rp z }.$$
By classic results, see e.g. \cite[Theorem 7.20, page 208]{Grigoryan}, $\tilde p$ is a continuous function of $(z,x,y)\in \C_+\times\Omega\setminus\{0\}\times\Omega\setminus\{0\}$, it is  symmetric  in $x,y$ and it is holomorphic in $z$.\\
Furthermore, the same argument as in \cite[Theorem 4.4]{cal-met-negro-spina }  proves that the family $\{e^{t\tilde L}:\ t \ge 0\}$ satisfies the Davies-Gaffney estimate  in $L^2\left(\Omega,|x|^{2-N}\,dx\right)$ that is 
\begin{align*}
\left|\left(e^{t\tilde L}f_1, f_2\right)_{L^2\left(\Omega,|x|^{2-N}\,dx\right)}\right|\leq\exp\left(-\frac{r^2}{4t}-\frac {\lambda_1}{ 3} t \right)\|f_1\|_{L^2\left(\Omega,|x|^{2-N}\,dx\right)}\|f_2\|_{L^2\left(\Omega,|x|^{2-N}\,dx\right)}
\end{align*}
for all $t>0$,  $U_1$,\ $U_2$ open subsets of $\Omega\setminus\{0\}$, $f_i$ in $L^2\left(U_i,|x|^{2-N}\,dx\right)$ and $r:=d(U_1, U_2)$.  Applying \cite[Theorem 4.1]{CS} to the operator $-\frac{\lambda_1}{3}-\tilde L$ we get, for every $z\in\C_+$, $x, y\in \Omega\setminus\{0\}$ (here the joint continuity of $\tilde{p}(t, \cdot, \cdot)$ is used)
\begin{align*}
|\tilde p(z,x,y)|\leq C_1({\rm Re}\,z)^{-\frac{N}{2}}\left(1+{\rm Re}\,\frac{|x-y|^2}{4z}\right)^\frac{N}{2}\exp\left(-\frac{\lambda_1}{3} {\rm Re}\,z-\Rp\dfrac{|x-y|^2)}{4z}\right).
\end{align*}
Recalling \eqref{relation kernels}, the heat kernel $p$ of $L$, taken with respect the Lebesgue measure, satisfies
\begin{equation*}
p(z,x,y)=|x|^{-s_1}|y|^{-s_1}\tilde p(z,x,y)
\end{equation*}
and \eqref{HK est D=0} follows.

To prove the second statement we observe that
\begin{align}\label{GF D=0 eq 1} 
G(x,y)=\int_0^\infty p(t,x,y)\,dt\leq C |x|^{-s_1}|y|^{c-s_1}\int_0^\infty h(t)\,dt,
\end{align}
where we  put $h(t)=t^{-\frac{N}{2}}e^{-\frac{\lambda_1}{3} t}\exp\left(-\dfrac{|x-y|^2)}{(4+\epsilon)t}\right)$. Using  \cite[Formula (29), page 146]{Erdelyi}, we have

\begin{align*}
\int_0^\infty h(t)\,dt&=2\left(\frac{3|x-y|^2}{\lambda_1(4+\epsilon)}\right)^{-\frac{N-2}{4}}K_{\frac{N-2}{2}}\left(2\frac{|x-y|}{\sqrt{4+\epsilon}}\sqrt{\frac{\lambda_1}{3}}\right)\\[1ex]
&=C|x-y|^{-\frac{N-2}{2}}K_{\frac{N-2}{2}}\left(c|x-y|\right),
\end{align*}
 where the $K_\nu$ is the modified Bessel function and satisfies the following asymptotics, see e.g., \cite[9.6 and 9.7]{AS}. 
\begin{align*}
\text{If\;} \nu>0,\quad K_\nu(r)\approx
\begin{cases}
\sqrt{\frac{\pi}{2}}\,r^{-\frac{1}{2}}e^{-r},\quad&\text{if}\quad r\rightarrow\infty;\\[1ex]
\frac{\Gamma(\nu)}{2}\left(\frac{r}{2}\right)^{-\nu},\quad&\text{if}\quad r\rightarrow 0;
\end{cases}\\[2ex]
K_0(r)\approx
\begin{cases}
\sqrt{\frac{\pi}{2}}\,r^{-\frac{1}{2}}e^{-r},\quad&\text{if}\quad r\rightarrow\infty;\\[1ex]
-\log{r},\quad&\text{if}\quad r\rightarrow 0.
\end{cases}
\end{align*}
Inserting this relations into \eqref{GF D=0 eq 1} we get  if $N>2$
\begin{align*}
G(x,y)\leq C |x|^{-s_1}|y|^{c-s_1}e^{-c|x-y|}\left(1\wedge |x-y|\right)^{2-N}
\end{align*}
and if $N=2$
\begin{align*}
G(x,y)\leq C |x|^{-s_1}|y|^{c-s_1}
\begin{cases}
e^{-c|x-y|},\quad&\text{if}\quad |x-y|\geq 1;\\[2ex]
 1-\log \left(|x-y|\right),\quad&\text{if}\quad |x-y|\leq 1.
\end{cases}
\end{align*}
\qed

\subsection{Main result}
As  in the cases $\Omega=B$ or $\Omega=\R^N$, we define
\begin{align*}
D_{p,\alpha}(\Omega):&=\left\{u:\ |x|^{\alpha-2}u,\ |x|^{\alpha}Lu\in L^p\left(\Omega\right),\ u=0 \text{ on } \partial \Omega\right\}.
\end{align*} Our main result is the following

\begin{teo} \label{RellichOmega}
Let $N\geq 2$, $1< p <\infty$ and assume that (\ref{D}) holds. 
Rellich inequalities 
\begin{align*}
\||x|^\alpha Lu\|_p \geq C\||x|^{\alpha-2} u\|_p,\quad u\in D_{p,\alpha}(\Omega)
 \end{align*} hold if and only if
\begin{align*}
\alpha&<N\Bigl (\frac12-\frac1p\Bigr)+1+\frac{c}{2}+ \sqrt{D} \quad\text{and}\;\\[1ex]
\alpha&\neq N\Bigl (\frac12-\frac1p\Bigr)+1+\frac{c}{2}-\sqrt {D+\lambda_n}, \quad\forall\, n\in\N_0.
\end{align*}
\end{teo}
{\sc Proof.} 
We  first prove that, if $\alpha$ is as in the assumptions,  Rellich inequalities are true. 
Let $B_R\subseteq\R^N$ be  such that $\Omega\subset B_R$, $R>0$. Without loss of generality we may assume that $R=1$. For a sufficiently small $\delta>0$, set 
\begin{align*}
\Omega_{\delta}:=\left\{x\in\Omega:\ \mbox{dist}(x,\partial\Omega)<\delta\right\}.
\end{align*}
We take  a linear extension operators  $E: W^{2,p}(\Omega)\to W^{2,p}_0(B)$  such that 
$$\|Eu\|_{W^{2,p}(B\setminus\Omega)}\leq C\|u\|_{W^{2,p}(\Omega_\delta)}$$  and let $u\in C_{c,0}^2(\Omega)$.
By Theorem \ref{RellichFJ} and since all coefficients are bounded in $B\setminus \Omega$,  we have
\begin{align*}
\int_	\Omega |x|^{(\alpha-2)p}|u|^p\, dx &\leq \int_B |x|^{(\alpha-2)p}|Eu|^p\, dx\leq C \int_B |x|^{\alpha p}|L(Eu)|^p\, dx\\&\leq  C\left( \int_\Omega |x|^{\alpha p}|Lu|^p\, dx+\|Eu\|^p_{W^{2,p}(B\setminus\Omega)}\right)\\&\leq C\left( \int_\Omega |x|^{\alpha p}|Lu|^p\, dx+\|u\|^p_{W^{2,p}(\Omega_\delta)}\right).
\end{align*}
By the interior estimates for elliptic equations (see \cite[Theorem 1, Sec. 4, Ch.9]{Krylov})
\begin{align*}
\|u\|_{W^{2,p}(\Omega_\delta)}&\leq C\left( \|Lu\|_{p,\Omega_{2\delta}}+\|u\|_{p,\Omega_{2\delta}}\right)\\[1ex]
&\leq C\left( \||x|^\alpha Lu\|_{p,\Omega_{2\delta}}+\|u\|_{p,\Omega_{2\delta}}\right)\\[1ex]
&\leq C\left( \||x|^\alpha Lu\|_{p,\Omega}+\|u\|_{p,\Omega_{2\delta}}\right).
\end{align*}
To conclude the proof we  show that $\|u\|_{p,\Omega_{2\delta}}\leq C\||x|^\alpha Lu\|_{p,\Omega}$. 
\smallskip

Set $f=-|x|^\alpha Lu$.  Since $u\in C_{c,0}^2(\Omega) \subset D(L)$ and $L$ is invertible, by Proposition \ref{LForm}, then  $u=(-L)^{-1}\frac{f}{|x|^\alpha}$. Using the estimates proved in Section \ref{Section estimate}, the Green function $G$ of  $L$ in $\Omega$ satisfies
\begin{equation*} 
 0 \leq G(x,y) \leq C\, G_0(x,y),
\end{equation*}
where $G_0$ is defined in  Proposition \ref{comparison Green} when $D>0$ and in Proposition \ref{comparison Green D=0} when $D=0$.

Let us suppose preliminarily that $D>0$. Then, for $x\in\Omega_{2\delta}$,
\begin{align*}
|u(x)|&\leq C\int_\Omega G_0(x,y)|y|^{-\alpha}|f(y)|\, dy\\[1ex]&=
C\int_{\{|x||y|\geq |x-y|^2\}}G_0(x,y)|y|^{-\alpha}|f(y)|\, dy\\[1ex]
&+ C\int_{\{|x||y|\leq |x-y|^2\}}G_0(x,y)|y|^{-\alpha}|f(y)|\, dy=: u_1(x)+u_2(x).
\end{align*}
Since $x\in\Omega_{2\delta}$, there exists $a>0$ such that $|x|\geq a>0$.
Consider first $u_1(x)$.  
\begin{align*}
|u_1(x)|&\leq C\int_{\left\{|x||y|\geq |x-y|^2,\ |y|\geq \frac{a}{2}\right\}}G_0(x,y)|y|^{-\alpha}|f(y)|\, dy\\&+ C\int_{\left\{|x||y|\geq |x-y|^2,\ |y|\leq \frac{a}{2}\right\}}G_0(x,y)|y|^{-\alpha}|f(y)|\, dy=: I_1(x)+I_2(x).
\end{align*}
For the first term $I_1(x)$ we get 
\begin{align*}
I_1(x)\leq C\int_\Omega |x-y|^{2-N}|f(y)|\, dy,\quad \text{if\, }N>2,\\[1ex]
I_1(x)\leq C\int_\Omega \left|\log|x-y|\right||f(y)|\, dy,\quad \text{if\, }N=2,
\end{align*}
which  therefore implies $\|I_1\|_{p,\Omega_{2\delta}}\leq C\|f\|_{p,\Omega}$. For $I_2(x)$, observe that  $|x-y|\geq\frac{a}{2}$, therefore $|x||y|\geq |x-y|^2\geq\frac{a^2}{4}$ and, recalling that $\Omega\subset B$, $|y|\geq \frac{a^2}{4|x|^2}\geq \frac{a^2}{4}$. It follows that 
$$I_2(x)\leq C\int_\Omega |f(y)|\, dy,$$ and $\|I_2\|_{p,\Omega_{2\delta}}\leq C\|f\|_{p,\Omega}$. Then  $\|u_1\|_{p,\Omega_{2\delta}}\leq C\|f\|_{p,\Omega}$. 
Consider now $u_2(x)$; since $|x|\geq a$,
\begin{align*}
|u_2(x)|&\leq C\int_{\{|x||y|\leq |x-y|^2\}} \frac{(|x||y|)^{\sqrt{D} -\frac{N-2}{2}}}{|x-y|^{2\sqrt{D}}}|y|^{\frac{c}{2}-\alpha}|f(y)|\, dy.
\end{align*}
As before, consider separately
$$J_1(x):=\int_{\{|x||y|\leq |x-y|^2,\ |y|\geq\frac{a}{2}\}}  \frac{(|x||y|)^{\sqrt{D} -\frac{N-2}{2}}}{|x-y|^{2\sqrt{D}}}|y|^{\frac{c}{2}-\alpha}|f(y)|\, dy$$ and
$$J_2(x):=\int_{\{|x||y|\leq |x-y|^2,\ |y|\leq\frac{a}{2}\}} \frac{(|x||y|)^{\sqrt{D} -\frac{N-2}{2}}}{|x-y|^{2\sqrt{D}}}|y|^{\frac{c}{2}-\alpha}|f(y)|\, dy.$$ 
Concerning $J_1$, we have $\frac{(|x||y|)^{\sqrt{D}}}{|x-y|^{2\sqrt{D}}}\leq 1$ and $(|x||y|)^{-\frac{N-2}{2}+\frac{c}{2}-\alpha}\leq C$, therefore
$$J_1(x)\leq C\int_\Omega|f(y)|\, dy$$ and
$\|J_1\|_{p,\Omega_{2\delta}}\leq C\|f\|_{p,\Omega}$. Finally, for $J_2$ we have $|x-y|\geq \frac a 2$ and 
$$J_2(x)\leq 
C\int_{\{|x||y|\leq |x-y|^2,\ |y|\leq\frac{a}{2}\}} |y|^{\sqrt{D} -\frac{N-2}{2}+\frac{c}{2}-\alpha}|f(y)|\, dy\leq C\|f\|_{p,\Omega}\,\left\||y|^{\sqrt{D} -\frac{N-2}{2}+\frac{c}{2}-\alpha}\right\|_{p',\Omega}.$$ 
The last norm is finite if and only if
$$\left(\sqrt{D} -\frac{N-2}{2}+\frac{c}{2}-\alpha\right)p'>-N $$
which is equivalent to $\alpha<N\Bigl (\frac12-\frac1p\Bigr)+1+\frac{c}{2}+ \sqrt{D}$, our assumption.

Let us suppose now that $D=0$. Then, similarly, we write for $x\in\Omega_{2\delta}$,
\begin{align*}
|u(x)|&\leq C\int_\Omega G_0(x,y)|y|^{-\alpha}|f(y)|\, dy\\[1ex]&=
C\int_{\{|x-y|\leq 1\}}G_0(x,y)|y|^{-\alpha}|f(y)|\, dy\\[1ex]
&+ C\int_{\{|x-y|\geq 1\}}G_0(x,y)|y|^{-\alpha}|f(y)|\, dy=: u_1(x)+u_2(x).
\end{align*}
Concerning $u_1$  we get 
\begin{align*}
u_1(x)\leq C\int_\Omega |x-y|^{2-N}|f(y)|\, dy,\quad \text{if\, }N>2,\\[1ex]
u_1(x)\leq C\int_\Omega \left|\log|x-y|\right||f(y)|\, dy,\quad \text{if\, }N=2,
\end{align*}
which implies $\|u_1\|_{p,\Omega_{2\delta}}\leq C\|f\|_{p,\Omega}$ as before.
 Finally, for $u_2$ we have 
$$u_2(x)\leq 
C\int_{\Omega} |y|^{ -\frac{N-2}{2}+\frac{c}{2}-\alpha}|f(y)|\, dy\leq C\|f\|_{p,\Omega}\,\left\||y|^{\-\frac{N-2}{2}+\frac{c}{2}-\alpha}\right\|_{p',\Omega}$$ 
 which is finite if and only if  $\alpha<N\Bigl (\frac12-\frac1p\Bigr)+1+\frac{c}{2}$, our assumption when $D=0$ (note that in this case  $c-s_1= -\frac{N-2}{2}+\frac{c}{2}$).

Let us now show that the conditions on $\alpha$ are also necessary and here we do not need to distinguish between $D>0$ and $D=0$.
\\When $\alpha= N\Bigl (\frac12-\frac1p\Bigr)+1+\frac{c}{2}-\sqrt {D+\lambda_n}, \  n\in\N_0$, Rellich inequalities fail, by Example 2.3.
Let  $\alpha> N\Bigl (\frac12-\frac1p\Bigr)+1+\frac{c}{2}+ \sqrt{D}$ and assume, as above, that  $\Omega\subseteq B$. Let $s_2$ be  defined in \eqref{defs}  and 
\begin{align*}
u(x):=|x|^{-s_2} , \quad x\in \Omega.
\end{align*}
 Then (see Proposition \ref{Counterexample in B}),  $Lu=0$  and $|x|^{\alpha-2} u\in L^p\left(B\right)$ since
$\alpha-2-s_2>-\frac{N}{p}$. On the other hand $u\notin D_{p,\alpha}(\Omega)$ since $u$ does not vanish on $\partial\Omega$. We use Lemma \ref{Maximum principle for L} and, for $0<\epsilon<1$, let    $v_\epsilon\in C^2\left(\Omega\setminus\bar{B_\epsilon}\right)\cap C\left(\bar \Omega\setminus B_\epsilon\right)$ satisfy
\begin{align*}
\begin{cases}
-Lv_\epsilon=0,\quad & \text{in}\quad \Omega\setminus\bar{B_\epsilon},\\
v_\epsilon=u,\quad & \text{in}\quad \partial\Omega,\\
v_\epsilon=\epsilon^{-s_1},\quad & \text{in}\quad \partial B_\epsilon.
\end{cases}
\end{align*}
Since $s_2>s_1$, one has, by construction, $v_\epsilon(x)\leq |x|^{-s_1}$ for every $x\in \partial\Omega\cup\partial B_\epsilon$. It follows from Lemma \ref{Maximum principle for L} that $0\leq v_\epsilon(x)\leq |x|^{-s_1}$ in $\Omega\setminus B_\epsilon$. Using local elliptic regularity and a standard diagonal argument, we prove that  $v_\epsilon$ converges, up to  subsequences, to a function $v$ in $W^{2,p}_{loc}\left(\Omega\setminus\{0\}\right)$. By construction $v$ satisfies  $v=u$ in $\partial\Omega$ and $Lv=0$, $0\leq v\leq |x|^{-s_1}$ in $\Omega\setminus\{0\}$; in particular $|x|^{\alpha-2}v\in L^p(\Omega)$,  since $\alpha-2-s_1>-\frac{N}{p}$. Then  the function $w:=u-v$ satisfies $w=0$ in $\partial\Omega$ and $Lw=0$ in $\Omega\setminus\{0\}$. In particular $w\in D_{p,\alpha}(\Omega)$ but Rellich inequalities \eqref{Intr 1} fail for $w$.\\\qed

\subsection{Rellich inequalities in exterior domains}
Let $V\subseteq \R^N$ be  an  exterior domain  (that is the complement of a bounded set)  which is also  open, connected and does not not containing the origin. We also assume that   $\partial V$ is $C^{2,\beta}$. 
As  before, we define
\begin{align*}
D_{p,\alpha}(V):&=\left\{u:\ |x|^{\alpha-2}u,\ |x|^{\alpha}Lu\in L^p\left(\Omega\right),\ u=0 \text{ on } \partial V\right\}.
\end{align*} 

\begin{prop} \label{RellichExterior}
Let $N\geq 2$, $1< p <\infty$ and assume that (\ref{D}) holds. 
Rellich inequalities 
\begin{align*}
\||x|^\alpha Lu\|_p \geq C\||x|^{\alpha-2} u\|_p,\quad u\in D_{p,\alpha}(V)
 \end{align*} hold if and only if
\begin{align*}
\alpha&>N\Bigl (\frac12-\frac1p\Bigr)+1+\frac{c}{2}- \sqrt{D} \quad\text{and}\;\\[1ex]
\alpha&\neq N\Bigl (\frac12-\frac1p\Bigr)+1+\frac{c}{2}+\sqrt {D+\lambda_n}, \quad\forall\, n\in\N_0.
\end{align*}
 When $V=B_r^c$ the same result  holds  when $D<0$ (replacing the square roots  with their real parts) and for $p=1,\infty$.
\end{prop}

{\sc Proof.} 
 For $u\in D_{p,\alpha}(V)$, we use the Kelvin transform $u(x)=|x|^{2-N}v\left (\frac{x}{|x|^2}\right )$ where $ v $ is defined in the bounded domain  $\Omega=\left\{x\in\R^N: x/|x|^2\in V\right\}$,  which contains the origin. Then by elementary computation
$$L u(x)=|x|^{-N-2}\tilde{L} v\left (\frac{x}{|x|^2} \right )$$ where  
\begin{align*}
\tilde{L}=\Delta +\tilde c\frac{x}{|x|^2}\cdot\nabla -\frac{\tilde b}{|x|^2},\quad 
\tilde c:=-c,\quad \tilde b:=b+(N-2)c.
\end{align*}
In particular its discriminant $\tilde D$ satisfies $\tilde D=D$.
Setting $y=x/|x|^2$, $dx=|y|^{-2N} dy$ we see that the inequality
$$
\||x|^\alpha L u\|_{L^p\left(V\right)} \geq C\||x|^{\alpha-2} u\|_{L^p\left(V\right)}
$$
is equivalent to 
$$
\||x|^{\tilde \alpha} \tilde L v\|_{L^p\left(\Omega\right)} \geq C\||x|^{{\tilde \alpha}-2} v\|_{L^p\left(\Omega\right)}
$$
with the same constant $C$ and $\tilde \alpha:=-\alpha+N+2-2N/p$. The statements then follow from Theorems \ref{RellichFJ} and \ref{RellichOmega}. \qed

\section{Critical cases in $L^p(\R^N)$}

In this section we assume that $\Omega$ coincides with  $\R^N$ and prove that, when Rellich inequalities fail,  modified inequalities which include logarithmic terms are still valid. The situation is similar to Hardy inequality, when the classical one fails.
By Theorem \ref{RellichFJ} Rellich inequalities fail in $\R^N$ if and only if
\begin{align} \label{fail}
\alpha= N\Bigl (\frac12-\frac1p\Bigr)+1+\frac{c}{2}\pm \Rp \sqrt {D+\lambda_n}.
\end{align}
or equivalently when 
\begin{equation} \label{failmod}
b +\gamma_p(\alpha,c)+\lambda_n = 0 \quad {\rm  for\  some}\  n \in \N_0.
\end{equation}
To study these cases we need an unweighted one dimensional result for a general second order operator on the half line.

\begin{prop} \label{OneD}
Consider the operator with real constant coefficients 
$$
\Gamma = D^2+\beta D
$$
 in $(0,\infty)$ and fix $a>0$. If  $\beta \neq 0$, then    
for every $v \in C_c^2 (\R_+)$,
\begin{equation} \label{one-dim-Rellich}
\left\|\frac{v}{s}\right\|_{L^p(a,\infty)}\leq C \|\Gamma v\|_{L^p(0,\infty)}. 
\end{equation}
for $1<p \le \infty$ and 
\begin{equation} \label{one-dim-Rellich1}
\left\|\frac{v}{s^{1+\eps}}\right\|_{L^1(a,\infty)}\leq C_\eps \|\Gamma v\|_{L^1(0,\infty)}. 
\end{equation}
for  $\eps>0$. The weaker 
inequalities
\begin{equation} \label{one-dim-Rellichw}
\left\|\frac{v}{s^2}\right\|_{L^p(a,\infty)}\leq C \|\Gamma v\|_{L^p(0,\infty)}
\end{equation}
 and
\begin{equation} \label{one-dim-Rellichw1}
\left\|\frac{v}{s^{2+\eps}}\right\|_{L^1(a,\infty)}\leq C_\eps \|\Gamma v\|_{L^1(0,\infty)}
\end{equation}
hold when $\beta=0$.
\end{prop}

In the proof we need the following lemma. 
\begin{lem} \label{GreenFunc}
Let $v\in C_c^2(\R_+)$ and 
$f=\Gamma v$, with $\beta \neq 0$. Then 
\begin{align} \label{rappresentazione}
v(s)=-\frac{1}{\beta}\left (\int_0^s e^{-\beta(s-\sigma)}f(\sigma)\, d\sigma+\int_s^\infty f(\sigma)\,d\sigma \right ).
\end{align}
Moreover, one has 
\begin{equation} \label{orthogonal}
\int_0^\infty f(\sigma)\,d\sigma=
\int_0^\infty e^{\beta \sigma}f(\sigma)\,d\sigma=0.
\end{equation}
\end{lem}
{\sc Proof.} Identity (\ref{orthogonal})  holds since $1, e^{\beta s}$ are solution of the adjoint $\Gamma^*=D^2-\beta D$. If  $w$ is the right hand side of (\ref{rappresentazione}), by the variation of constants formula, $\Gamma w=f$ and, by (\ref{orthogonal}) and since $f$ has a compact support, $w$ has a compact support, too. On the other hand, $\Gamma (v-w)=0$, hence $v-w=c_1 +c_2 e^{\beta s}$. Since both have a compact support in $(0,\infty)$, then $c_1=c_2=0$.
 \qed

{\sc Proof.} (Proposition \ref{OneD})
Let $f:=\Gamma v$ and assume first that $\beta=0$. Then
$$
\frac{|v(s)|}{s^2} \le s^{-2}\int_0^s (s-\sigma)|f(\sigma)|\, d\sigma \le s^{-1} \int_0^s |f(\sigma)|\, d\sigma
$$ and (\ref{one-dim-Rellichw}) follows from Hardy inequality. When $p=1$ we write
$$
|v(s)|=\left |-\int_0^s d\sigma \int_\sigma ^\infty v''(\xi)\, d\xi\right | \le s\|v''\|_1
$$ and  (\ref{one-dim-Rellichw1}) is immediate.

We assume now that $\beta \neq 0$ and use (\ref{rappresentazione})
\begin{align*} 
v(s)=-\frac{1}{\beta}\left (\int_0^s e^{-\beta(s-\sigma)}f(\sigma)\, d\sigma+\int_s^\infty f(\sigma)\,d\sigma \right )=:v_1+v_2.
\end{align*}
Since (\ref{orthogonal}) holds, then 
$$
\frac{|v_2(s)|}{s} \le C\frac{1}{s}\int_0^s |f(\sigma)|\, d\sigma,  \quad \left \|\frac{v_2}{s}\right \|_{L^p(a, \infty)} \le C\|f\|_p, 
$$
 if $1<p \le \infty$, by Hardy inequality. When $p=1$, then $|v_2(s)| \le C\|f\|_1$. This shows that (\ref{one-dim-Rellich}), (\ref{one-dim-Rellich1}) hold for $v_2$. 
Since by (\ref{orthogonal})
$$
-\beta v_1(s)=\int_0^s e^{-\beta (s-\sigma)}f(\sigma)\, d\sigma=-\int_s^\infty e^{-\beta(s-\sigma)}f(\sigma)\, d\sigma=e^{-\beta \cdot} \chi_{(0, \infty)}*f (s)=-e^{-\beta \cdot} \chi_{(- \infty,0)}*f (s),
$$
the estimate
$$
\|v_1\|_{L^p(0, \infty)} \le C\|f\|_{L^p(0,\infty)}
$$ follows from Young's inequality for every $1 \le p \le \infty$ and concludes the proof.
\qed

In the following theorem we concentrate on the singularity at 0, hence we consider only $C^2$-functions vanishing on a neighbourhood of the origin and with a fixed common support which can be assumed to be $ B_{R/2}$.  
We set 
$$
{\cal D}_R=\{u \in C^2(\R^N): u=0 \ {\rm   in\  a\  neighborhood\  of  }\ 0,  \ \ {\rm spt }\ u \subset B_{R/2} \}.
$$

\begin{teo}  \label{critical}
Assume that 
$$\alpha= N\Bigl (\frac12-\frac1p\Bigr)+1+\frac{c}{2}\pm \Rp \sqrt {D+\lambda_n}$$ for some $n \in \N_0$
 or, equivalently, that (\ref{failmod}) holds. 

Then for $1<p\le \infty$
there exists a positive constant $C$, independent of $R$,  such that for every $u\in {\cal D}_R$
\begin{equation}  \label{Rellichlog}
\||x|^\alpha Lu\|_p \geq C\Big\||x|^{\alpha-2}\left|\log |R^{-1}x|\right|^{-2} u\Big\|_p \quad {\rm when }\ D+\lambda_n \le 0
\end{equation}
\begin{equation}  \label{Rellichlog1}
\||x|^\alpha Lu\|_p \geq C\Big\||x|^{\alpha-2}\left|\log |R^{-1}x|\right|^{-1} u\Big\|_p \quad {\rm when }\ D+\lambda_n >0.
\end{equation} \\

When $p=1$, inequalities (\ref{Rellichlog}) and (\ref{Rellichlog1}) hold with $|\log |R^{-1}x||^{-2}$ and  $|\log |R^{-1}x||^{-1}$ replaced by $|\log |R^{-1}x||^{-2-\eps}$ and  $|\log |R^{-1}x||^{-1-\eps}$, respectively.
\end{teo}

{\sc Proof.} By scaling we may assume that $R=1$.
By   Theorem \ref{RellichFJ},  Rellich inequalities hold in $D_{p,\alpha}(\R^N) \cap L^p_{ \neq  n}$. Then (\ref{Rellichlog1}) hold in 
 ${\cal D}_1 \cap L^p_{ \neq n}$, since the singularity at $0$ is weaker and $u$ has support in $ B_{1/2}$.
Since, by Lemma \ref{projection}
$$
L^p(\R^N)= L^p_{n}(\R^N) \oplus L^p_{\neq n}(\R^N)
$$
and $L$ preserves both $L^p_{n}(\R^N)$ and $ L^p_{\neq n}(\R^N)$, 
we have to show that (\ref{Rellichlog}) or (\ref{Rellichlog1}) or their variants for $p=1$ hold in ${\cal D}_1 \cap L^p_{n}(\R^N)$.

Let  $u(\rho, \omega)=c(\rho)P(\omega)$, where $P$ is a fixed spherical harmonic of order $n$. Using the transformation $c(\rho)=\rho^{-\alpha+2-\frac{N}{p}} v(-\log \rho)$ we have
\begin{align*}
&\||x|^\alpha Lu\|^p_p =\int_{\R^N}|x|^{\alpha p}\left|\Delta u+c\frac{x}{|x|^2}\nabla u-\frac{b}{|x|^2}u\right|^p\ dx\\[1ex]
&=\int_{S^{N-1}}|P(\omega)|^p\int_0^{\frac12} \rho^{\alpha p+N-1}\left|\frac{\partial^2 c(\rho)}{\partial\rho^2}+\frac{(N-1+c)}{\rho} \frac{\partial c(\rho)}{\partial \rho}-\frac{\lambda_n+b}{\rho^2}c(\rho)\right|^pd\rho\ d\omega\\[1ex]
&=\int_{S^{N-1}}|P(\omega)|^p\int_{\log 2}^\infty \left|\frac{\partial^2 v(s)}{\partial s^2}+\left(2\alpha-2-N +\frac{2N}{p}-c\right)\frac{\partial v(s)}{\partial s}-(\gamma_p(\alpha,c)+b+\lambda_n)v(s)\right|^p ds\ d\omega \\[1ex]
&=\int_{S^{N-1}}|P(\omega)|^p\int_{\log 2}^\infty \left|\frac{\partial^2 v(s)}{\partial s^2}+\left(2\alpha-2-N +\frac{2N}{p}-c\right)\frac{\partial v(s)}{\partial s} \right|^p ds\ d\omega.
\end{align*}
since $\gamma_p(\alpha,c)+b+\lambda_n=0$. At this point we apply Proposition \ref{OneD} with $\beta=2\alpha-2-N+\frac{2N}{p}-c$ after noticing that
\begin{align*}
\int_{\R^N}|x|^{\alpha p}\left|\frac{u(x)}{|x|^2|\log |x||^\gamma}\right|^p\ dx=
\int_{S^{N-1}}|P(\omega)|^p\int_{\log 2}^\infty \left|\frac{v(s)}{s^\gamma}\right|^p ds\ d\omega.
\end{align*}
Observe that, since 
$\alpha= N\Bigl (\frac12-\frac1p\Bigr)+1+\frac{c}{2}\pm \Rp \sqrt {D+\lambda_n}$, then $\beta\neq 0$ if and only if  $D+\lambda_n>0$. 
\qed

\section{Best constants and remainder terms} 

When 
$D:=b+\left(\frac{N-2+c}{2}\right)^2>0$ and 
\begin{equation*} 
 N\Bigl (\frac12-\frac1p\Bigr)+1+\frac{c}{2}-\sqrt {D}<\alpha < N\Bigl (\frac12-\frac1p\Bigr)+1+\frac{c}{2}+\sqrt {D}. 
\end{equation*}
we have seen in Proposition \ref{easy} that 
Rellich inequalities (\ref{Intr 1}) hold in $D_{p, \alpha}(\Omega)$ with the best constant 
\begin{equation} \label{best}
C:=b+\Bigl(\frac{N}{p}-2+\alpha\Bigr)\Bigl(\frac{N}{p'}-\alpha+c\Bigr).
\end{equation}
 As usual, $\Omega$ is an open bounded and connected  set containing $0$ and with a smooth boundary, or $\Omega =\R^N$. Best constants are not known in other cases, except for $p=2$ or in special subspaces, see \cite{rellich}.

A direct proof that, in the above range, the constant $C$ is optimal can be achieved by truncating  the function $u(x)=|x|^{2-\alpha-N/p}$ as in Example (\ref{esem1}).

\begin{lem} \label{c1}
Assume $1<p<\infty$. Under the above assumption on $\alpha$, Rellich inequalities hold in $D_{p,\alpha}(\R^N)\cap L^p_{\ge 1}(\R^N)$ with a constant $C_1>C$.
\end{lem}
{\sc{Proof.}} According to equation \eqref{disv} of Section 2, we have to show that the inequality
$$
\|\mu v-Av\|_p \ge C_1 \|v\|_p, \quad v \in D_{p,max}(\R^N)
$$ holds for a suitable $C_1>C$.
We revisit Theorem \ref{dissipativity} where,  we recall, $\mu=\lambda-\omega_p$ and $\lambda=\mu+\omega_p=b+\gamma_p(\alpha, c)=C$ (see also Lemma \ref{Parameters}). 
Theorem \ref{contractivity} holds in $D_{p,max}(\R^N)$ with a suitable $\omega_p^1 >\omega$, by the results in Section 2 of \cite{met-soba-spi3}, see in particular Proposition 2.8 and Remark 2.9. with $\beta=0$ therein. It follows that $\mu=\lambda-\omega=\lambda_1-\omega_p^1$ and $\lambda_1>\lambda$. Then estimate (\ref{contractivity}) holds with $\lambda_1 >\lambda=C$.
\qed
The remainder term can arise, therefore, when considering radial functions. To deal with them, we need the following auxiliary result.

\begin{lem} \label{aux}
Let $1<p<\infty$ and $\Gamma=D^2+\beta D-\lambda$ be an operator with real constant coefficients on $(0, \infty)$. Then for every $v \in C_c^2(0,\infty)$ and $\lambda >0$
$$
\|\Gamma v\|_p^p -\lambda^p \|v\|_p^p \ge \lambda^{p-1}\frac{p-1}{p^2} \int_0^\infty \frac{|v(s)|^p}{s^2}\, ds.
$$
\end{lem}
{\sc Proof. } We have
$$
\int_0^\infty (\lambda v-v''-\beta v')v|v|^{p-2}=\int_0^\infty \left (\lambda |v|^p +(p-1) |v'|^2 |v|^{p-2}-\beta v'v|v|^{p-2} \right ).
$$
Since $v'v|v|^{p-2}$ is the derivative of $p^{-1}|v|^p$, the last integral vanishes. By Hardy inequality of Proposition \ref{hardy2}, with $N=1, \beta =0$ we have
$$
\int_0^\infty |v'|^2|v|^{p-2} \ge \frac{1}{p^2} \int_0^\infty \frac{|v(s)|^p}{s^2}\, ds
$$
and therefore
$$\lambda \|v\|_p^p +\frac{p-1}{p^2} \int_0^\infty \frac{|v(s)|^p}{s^2}\, ds \le \|\Gamma v\|_p \|v\|_p^{p-1}.$$
Let 
$$
A^p=\|V\|_p^P, \quad B^p=\frac{p-1}{p^2} \int_0^\infty \frac{|v(s)|^p}{s^2}\, ds, \quad C=\|\Gamma v\|_p.
$$
then from $\lambda A^p+B^p \le CA^{p-1}$ we get $\lambda A \le C$ and 
$$
C^p-\lambda^p A^p \ge C^p-C\lambda^{p-1}A^{p-1}+\lambda^{p-1}B^p=C(C^{p-1}-\lambda^{p-1}A^{p-1})+\lambda^{p-1}B^p \ge \lambda^{p-1}B^p.$$
\qed

The main result of this section is stated below. As in the previous section we formulate it for functions belonging to 
$$
{\cal D}_R=\{u \in C^2(\R^N): u=0 \ {\rm   in\  a\  neighborhood\  of  }\ 0,  \ \ {\rm spt }\ u \subset B_{R/2} \}.
$$

\begin{teo} \label{Rellichremainder}
Let $1<p<\infty$, $D:=b+\left(\frac{N-2+c}{2}\right)^2>0$ and 
\begin{equation*} 
 N\Bigl (\frac12-\frac1p\Bigr)+1+\frac{c}{2}-\sqrt {D}<\alpha < N\Bigl (\frac12-\frac1p\Bigr)+1+\frac{c}{2}+\sqrt {D}. 
\end{equation*} 
If $C$ is the best constant defined in (\ref{best}), then there exists $c>0$, independent of $R$,  such that for every $u \in {\cal D}_R$
\begin{equation} \label{RR}
\Big\||x|^\alpha Lu\Big\|_p^p -C^p \Big\||x|^{\alpha-2} u\Big\|_p^p \ge c \Big\||x|^{\alpha-2}\left|\log |R^{-1}x|\right|^{-\frac{2}{p}} u\Big\|_p^p.
\end{equation}
\end{teo}
{\sc Proof.} By scaling we may assume that $R=1$. If $u \in {\cal D}_1$, we split $u=u_0+u_1$, where $u_0$ is radial and $u_1 \in L^p_{\ge 1}(\R^N)$. By Lemma \ref{c1}, inequality (\ref{RR}) holds for $u_1$.

For $u_0$ we proceed as in Theorem \ref{critical} writing
$u_0(\rho)=\rho^{-\alpha+2-\frac{N}{p}} v(-\log \rho)$. Then 
\begin{align*}
\||x|^\alpha Lu_0\|^p_p =N\omega_N\int_{\log 2}^\infty \left|\frac{\partial^2 v(s)}{\partial s^2}+\left(2\alpha-2-N +\frac{2N}{p}-c\right)\frac{\partial v(s)}{\partial s}-(\gamma_p(\alpha,c)+b)v(s)\right|^p ds.
\end{align*}
Next we use Lemma \ref{aux} with $\lambda=\gamma(\alpha,c)+b=C$ to obtain
\begin{align*}
\||x|^\alpha Lu_0\|^p_p-C^p\||x|^{\alpha-2}u\|_p^p&= N\omega_N (\|\Gamma v\|_p^p-C^p \|v\|_p) \ge N\omega_N C^{p-1}\frac{p-1}{p^2}\int_0^\infty \frac{|v(s)|^p}{s^2}\, ds \\[1ex]
&=C^{p-1}\frac{p-1}{p^2}\Big\||x|^{\alpha-2}\big|\log |x|\big|^{-\frac{2}{p}} u\Big\|_p^p.
\end{align*}

The general case now follows, since $L_0(\R^N), L_{\ge 1}(\R^N)$ are invariant under $L$ and under multiplication by radial weights and since $|u|^p_p:=\|u_0\|_p^p +\|u_1\|_p^p$ is an equivalent norm on $L^p(\R^N)$.
\qed

\section{Appendix}

\subsection{Approximation on Sobolev spaces on domains}

Let  $V$ be a $C^{2,\beta}$ bounded connected open subset of $\R^N$ and let
$A$ be a uniformly elliptic operator $A=\mbox{tr}(A(x)D^2)+c(x)\cdot\nabla-b(x)$, with $C^\beta$ coefficients, endowed with Dirichlet boundary conditions. We  recall that for $1<p<\infty$
\begin{align*}
D_p(\Omega)=W^{2,p}(\Omega)\cap W^{1,p}_0(\Omega),
\end{align*}
whereas for $p=1$ 
 \begin{align*}
D_1(\Omega)=\left\{u\in W^{1,p}_0(\Omega): \mbox{tr}(A(x)D^2u)\in L^1(\Omega)\right\},
\end{align*}
and for  $p=\infty$ 
\begin{align*}
D_\infty(\Omega)=\left\{u\in C^1(\Omega)\cap C_0(\overline{\Omega}): \mbox{tr}(A(x)D^2u)\in C(\overline{\Omega})\right\},
\end{align*}
both endowed with the graph norm.

\begin{prop}\label{Sobolev approximation 1,infty}
Under the above assumptions the set
$$
C^2_0(\Omega)=\{u \in C^2(\overline{\Omega}): u=0 \ {\rm on}\ \partial \Omega \}
$$ is dense in $D_p(\Omega)$ for every $1 \le p\le \infty$.
\end{prop}
{\sc Proof. } Let $\lambda>0$ such that $\lambda-A$ is invertible from $D_p(\Omega)$ to $L^p(\Omega)$. If $u \in D_p(\Omega)$, $f=\lambda u-Au$ and $(f_n) \subset C^\beta (\Omega)$ tends to $f$ in $L^p(\Omega)$, then $u_n=(\lambda-A)^{-1}f_n$ belongs to $C^{2,\beta}(\Omega)$, by the Schauder  theory, vanishes at $\partial \Omega$ and approximates $u$ in the graph norm.
\qed
The following partition of unity of $\Omega$ has been used several times. 
\begin{prop}\label{Partition unity}
Let $0\leq\beta\leq 1$ and let $\Omega$ be a bounded connected open subset of $\R^N$ whose boundary $\partial \Omega$ is of class $C^{2,\beta}$. Then there exist $\delta>0$ such that the distance function $x\mapsto \mbox{dist}(x,\partial\Omega)$ is $C^{2,\beta}$ over the set 
\begin{align*}
K_{\delta}:=\left\{x\in\R^N:\ \mbox{dist}(x,\partial\Omega)<\delta\right\}.
\end{align*}
In particular $K_{\delta}$ and the subset
\begin{align*}
\Omega_{\delta}:=K_{\delta}\cap\Omega
\end{align*}
have $C^{2,\beta}$ boundary. Furthermore there exists an open subset $\Omega_0\subset\subset\Omega$ for which  $\overline\Omega=\overline\Omega_\delta\cup\Omega_0$ and there exists a partition of  unity $\{\eta_\delta^2,\eta_0^2\}$ such that
\begin{itemize}
\item[(i)] $\eta_\delta\in C_c^\infty(K_{\delta})$, $0\leq\eta_\delta\leq 1$, $\eta_\delta=1$ in $\overline\Omega_{\frac \delta 2}$;
\item[(ii)] $\eta_0\in C_c^\infty(\Omega_0)$, $0\leq\eta_\delta\leq 1$;
\item[(iii)] $\eta_\delta^2+\eta_0^2=1$ in $\overline\Omega$.
\end{itemize}
\end{prop}
{\sc{Proof.}}  \cite[Lemma 14.16]{gil-tru} proves the case $\beta=0$  and that, for sufficiently small $\delta$, for  every point  $x\in K_{\delta}$ there exist a unique $y\in\partial\Omega$ such that $|x-y|=d(x,\partial\Omega)$.  The result for $\beta>0$ then follows by \cite{Li-Nirenberg}. The existence of such a partition of unity is a standard result.\qed

\subsection{Some results on spectral theory } \label{spectral}
We collect some definitions and results from spectral theory which are used  throughout the paper. Let $X$ be a Banach space and let $A$ be a closed operator $A:D(A)\subseteq X\to X$. The spectrum of $A$ is denoted by $\sigma (A)$ and the resolvent set $\C \setminus \sigma (A)$ by $\rho (A)$.
\begin{defi}
The set
$$P\sigma(A):=\{\lambda\in\C: \lambda-A\ \textrm{is not injective}\}$$
is called the point spectrum of $A$. Moreover each $\lambda\in P\sigma(A)$ is called an eigenvalue and each $0\neq x\in D(A)$ satisfying $(\lambda-A)x=0$ is an eigenvector of $A$ (corresponding to $\lambda$).
\end{defi}

\begin{defi}\label{defi A spectrum}
The set 
$$A\sigma(A):=\{\lambda\in\C: \lambda-A\ \textrm{is not injective or}\  \rm rg(\lambda-A)\ \textrm{is not closed in X}\}$$
is called the approximate point spectrum of $A$. Obviously $P\sigma(A)\subseteq A\sigma(A)$.
\end{defi}

\begin{defi}\label{defi R spectrum}
The set 
$$R\sigma(A):=\{\lambda\in\C: \rm rg(\lambda-A)\ \textrm{is not dense in X}\}$$
is called the residual spectrum of $A$. 
\end{defi}

Note that $P\sigma (A) \subset A\sigma (A)$, that $P\sigma (A)$ and  $R\sigma (A)$, as well as $A\sigma (A)$ and $R\sigma (A)$ may overlap and that $\sigma (A)=A\sigma (A) \cup R\sigma (A)$.

\begin{lem} (\cite[Lemma 1.9, Chapter IV]{engel-nagel})\label{Char Aspectrum}
A number $\lambda\in\C$ belongs to $A\sigma(A)$ if and only if there exists a sequence $(x_n)_{n\in\N}\subset D(A)$, called an approximate eigenvector, such that $\|x_n\|=1$ and $\lim_{n\to\infty}\|Ax_n-\lambda x_n\|=0$.
\end{lem}
The following result is an elementary consequence of the previous Lemma.
\begin{prop}  \label{Rellich-spectrum}
The following properties are equivalent
\begin{itemize}
\item[(i)] There exists $C>0$ such that
\begin{align*}
\|x\|\leq C\|\lambda x-Ax\|,\quad \forall x\in D(A);
\end{align*}
\item[(ii)] $\lambda$ does not belong to the approximate point spectrum of $A$.
\end{itemize}
\end{prop}
The next Proposition implies that $A\sigma(A)$ is never empty.
\begin{prop} \label{boundary-spectrum} \cite[Proposition 1.10, Chapter IV]{engel-nagel}
The topological boundary of the spectrum is contained in the approximate point spectrum. 
\end{prop}

\subsection{Spectrum of a second order ordinary differential operator} 
We present the following elementary result on the spectrum of  the second order ordinary differential operator $B=D^2+\beta D$ in $L^p([0, \infty[)$, endowed with Dirichlet boundary condition at 0, that is 
$$
D(B)=\{u \in W^{2,p}([0, \infty[): u(0)=0\}.
$$

As usually $L^\infty([0,+\infty))$ stands for  $C_0^0([0,+\infty)$. Here $\beta\in\R$ and we recall that 
\begin{align*}
{\cal Q}&=\left\{\lambda\in \C: \  ({\rm Im}\lambda)^2 \leq -\beta^2 {\rm Re }\lambda \right\},\quad
{\cal P}=\left\{\lambda\in \C: \  ({\rm Im}\lambda)^2 = -\beta^2 {\rm Re }\lambda \right\}.
\end{align*}

Note that 
$$
{\cal P}(\kappa):= \{-\xi^2+i \beta \xi\;;\;\xi\in \R\}
$$
and that
\begin{equation} \label{distance parabola}
{\rm dist}(\lambda,{\cal P})^2
=
\begin{cases}
\lambda^2 & {\rm if}\ \lambda\geq-\frac{\beta^2}{2}, 
\\[1ex]
\beta^2(-\lambda-\frac{\beta^2}{4}) & { \rm if}\ \lambda<-\frac{\beta^2}{2}. 
\end{cases}
\end{equation}

 Observe that the spectrum of $B$ in $L^p(\R)$ is given by $\cal P$ and consists of approximate eigenvalues. This can be seen by noticing that the spectrum is independent of $p$ and using the Fourier transform in $L^2(\R)$. \\

For $\lambda \in \C$, we consider the solutions of the homogeneous equation $\lambda u -Bu=0$ given by $e^{\mu_it}$, $i=1,2$ where
$$\mu_1=\frac{-\beta-\sqrt{\beta^2+4\lambda}}{2},\quad \mu_2=\frac{-\beta+\sqrt{\beta^2+4\lambda}}{2}.$$ When $\lambda=-\beta^2/4$ then  $\mu_1=\mu_2=-b/2$ and  we substitute $e^{\mu_2 t}$ with 
 $te^{-\frac{\beta}{2}t}$.

\begin{lem} \label{sqrt}
The inequality ${\rm Re }\sqrt{\beta^2+4\lambda}<|\beta|$ holds if and only if $\lambda \in \overset{\mathrm{o}}{\cal Q}$. Similarly, ${\rm Re }\sqrt{\beta^2+4\lambda}>|\beta|$ if and only if $\lambda \not \in {\cal Q}$ and ${\rm Re }\sqrt{\beta^2+4\lambda}=|\beta|$ if and only if $\lambda \in {\cal P}$. Here $\sqrt z$ denotes any square root of $z$ with non negative real part.
\end{lem}
{\sc Proof.} If $\sqrt{\beta^2+4\lambda}=x+iy$, with $x \ge 0$, then $4\lambda=(x^2-y^2-\beta^2)+2i xy$ and $x=|\beta|$ if and only if  $({\rm Im}\lambda)^2 =-\beta^2 {\rm Re }\lambda$. The other cases are similar.
 \qed

\begin{prop} \label{ODE}
The spectrum of $B=D^2+\beta D$ in $L^p([0,+\infty))$, with Dirichlet boundary condition at 0, is given by $\sigma (B)={\cal Q}$.
More specifically we have 
\begin{itemize}
\item[(i)] if $\beta>0$, then \quad $\sigma(B)=A\sigma(B)={\cal Q}$,\ $P\sigma(B) \supset \overset{\mathrm{o}}{\cal Q}$;
\item[(ii)] if $\beta=0$, then \quad  $\sigma(B)=A\sigma(B)=(-\infty, 0]$;
\item[(iii)]  if $\beta<0$, then \quad $ A\sigma(B)={\cal P}$,\; $R\sigma(B)\setminus A\sigma(B)=\overset{\mathrm{o}}{\cal Q}$.
\end{itemize} 
\end{prop} 
{\sc Proof.} Let us prove preliminarily that $\cal{Q}^{\mbox{c}}\subseteq\rho(B)$ in all cases.  If $\lambda \notin \cal{Q}$ by the lemma above  
 ${\rm Re }\sqrt{\beta^2+4\lambda}>|\beta|$, hence ${\rm Re}\,\mu_1<0<{\rm Re}\,\mu_2$.
 It is then easy to see that $\lambda-B$ is invertible and that its inverse is given by the Green function
\begin{equation*}
G(t,s)=\left\{
\begin{array}{ll}
\displaystyle\frac{u_1(t)u_2(s)}{W(s)}  & \quad t\leq s,\\[3ex]
\displaystyle\frac{u_1(s)u_2(t)}{W(s)}   & \quad t\geq s;
\end{array}\right.\\
\end{equation*}
where   $u_1(t)=e^{\mu_2t}-e^{\mu_1t}$, $u_2(t)=e^{\mu_1t}$ and $W(t)=(\mu_1-\mu_2)e^{(\mu_1+\mu_2)t}=(\mu_1-\mu_2)e^{-\beta t}$ is their Wronskian.
 
Let us suppose now that $\lambda\in\overset{\mathrm{o}}{\cal Q}$ and assume first $\beta>0$. Then 
${\rm Re }\sqrt{\beta^2+4\lambda}<\beta$ and ${\rm Re } \mu_1 \leq {\rm Re } \mu_2<0$.
It follows that $\lambda $ is an eigenvalue with
 eigenfunction $u(t)=e^{\mu_1t}-e^{\mu_2t}$ (or $te^{-\frac{\beta}{2}t}$ when $\lambda=-\beta^2/4$). This proves that $\overset{\mathrm{o}}{\cal Q}\subseteq P\sigma(B)$ and case (i) is done, since the boundary of the spectrum is always contained in the approximate point spectrum, see Proposition \ref{boundary-spectrum}.\\
Assume now $\beta<0$ and still that $\lambda\in\overset{\mathrm{o}}{\cal Q}$. Then ${\rm Re }\sqrt{\beta^2+4\lambda}<-\beta$ and $0 <{\rm Re }\, \mu_1 \le  {\rm Re }\, \mu_2$, hence $\lambda-B$ is injective. Moreover, 
$\lambda-B$ is invertible with a continuous inverse from its domain  onto the closed subspace
$$X=\left\{f\in L^p\left([0,+\infty)\right):\ \int_0^\infty f(e^{-\mu_1s}-e^{-\mu_2s})ds=0\right\}$$
(with the usual change here and in what follows if $\lambda=-\beta^2/4$). 

Indeed if $u\in D(B)$ set  $f= (\lambda-B)u$ and  $B^*u=u''-\beta u'$. Since   $(e^{-\mu_1s}-e^{-\mu_2s})(0)=0$ and $(\lambda -B^*)(e^{-\mu_1s}-e^{-\mu_2s})=0$ , one has 
\begin{align*}
\int_0^\infty f(e^{-\mu_1s}-e^{-\mu_2s})ds&=\int_0^\infty (\lambda-B)u\,(e^{-\mu_1s}-e^{-\mu_2s})ds\\
&=\int_0^\infty u(\lambda -L^*)(e^{-\mu_1s}-e^{-\mu_2s})ds=0.
\end{align*}
On the other hand, if $f\in L^p([0,+\infty[)$ satisfies $\int_0^\infty f(e^{-\mu_1s}-e^{-\mu_2s})ds=0$, by the variation of constants method, one finds that 
$$u(t)=\frac{1}{\mu_1-\mu_2}e^{\mu_2 t}\int_t^\infty e^{-\mu_2 s}f(s)ds+\frac{1}{\mu_2-\mu_2}e^{\mu_1 t}\int_t^\infty e^{-\mu_1 s}f(s)ds$$
satisfies $u(0)=0$, $u \in D(B)$  and $(\lambda-B)u=f$. 

This proves that $\lambda-B$ is injective and that $rg(\lambda-B)=X\subset L^p([0,\infty[)$ which, recalling Definitions  \ref{defi A spectrum},  \ref{defi R spectrum}, gives $\overset{\mathrm{o}}{\cal Q}\subseteq R\sigma(B)\setminus A\sigma(B)$. Using again Proposition \ref{boundary-spectrum}, (iii) is proved.

When $\beta=0$ one sees that 
$A\sigma(B) \supset(-\infty,0]$ by truncating the functions $\sin (\sqrt {-\lambda}\, t)$.\\
\qed

An analogous result can be obviously proved in $L^p(]-\infty,0])$ using the isometry 
$$S:L^p([0,\infty[)\to L^p(]-\infty,0]),\quad Su(t)=u(-t).$$
\begin{prop} \label{ODE2}
The spectrum of $B=D^2+\beta D$ in $L^p(]-\infty,0])$, with Dirichlet boundary condition at 0, is given by $\sigma (B)={\cal Q}$.
More specifically we have 
\begin{itemize}
\item[(i)] if $\beta<0$, then \quad $\sigma(B)=A\sigma(B)={\cal Q}$,\ $P\sigma(B) \supset \overset{\mathrm{o}}{\cal Q}$;
\item[(ii)] if $\beta=0$, then \quad  $\sigma(B)=A\sigma(B)=(-\infty, 0]$;
\item[(iii)]  if $\beta>0$, then \quad $ A\sigma(B)={\cal P}$,\; $R\sigma(B)\setminus A\sigma(B)=\overset{\mathrm{o}}{\cal Q}$.
\end{itemize} 
\end{prop}

\end{document}